\documentclass[3p]{elsarticle}

\usepackage{amsmath}
\usepackage{amssymb}
\allowdisplaybreaks

\usepackage{fix-cm}
\DeclareMathAlphabet{\mathsfit}{T1}{\sfdefault}{\mddefault}{\sldefault}
\SetMathAlphabet{\mathsfit}{bold}{T1}{\sfdefault}{\bfdefault}{\sldefault}

\newcommand*{\bvec}[1]{\boldsymbol{#1}}
\newcommand*{\bmat}[1]{\boldsymbol{\mathsfit{#1}}}
\newcommand*{\bmatu}[1]{\boldsymbol{\mathsf{#1}}}

\DeclareMathOperator{\erf}{erf}
\DeclareMathOperator{\erfc}{erfc}
\DeclareMathOperator{\sgn}{sgn}
\DeclareMathOperator{\rect}{rect}
\DeclareMathOperator{\diam}{diam}

\newcommand*{\E}{\mathrm{e}}
\newcommand*{\I}{\mathrm{i}}

\newcommand*{\stokeslet}{S}
\newcommand*{\rotlet}{\Omega}
\newcommand*{\rotlett}{\boldsymbol{\Omega}}
\newcommand*{\stresslet}{T}
\newcommand*{\harmonic}{H}
\newcommand*{\biharmonic}{B}
\newcommand*{\pressurelet}{\Pi}
\newcommand*{\aB}{a_\biharmonic}
\newcommand*{\bB}{b_\biharmonic}
\newcommand*{\cB}{c_\biharmonic}

\newcommand*{\Kop}{\mathrm{K}}
\newcommand*{\Kopt}{\bmatu{K}}
\newcommand*{\kernel}{\bmat{G}}
\newcommand*{\matop}{\cdot}

\newcommand*{\per}{\mathcal{P}}
\newcommand*{\boxvar}{\mathcal{B}}
\newcommand*{\wavenumset}{\mathcal{K}}
\newcommand*{\fourier}{\mathcal{F}}
\newcommand*{\extended}[1]{#1'}
\newcommand*{\extbox}{\extended{\boxvar}}
\newcommand*{\extL}{\extended{L}}
\newcommand*{\extM}{\extended{M}}
\newcommand*{\gridvar}{\mathcal{G}}
\newcommand*{\isc}[1]{\bar{#1}}

\newcommand*{\rp}{\mathrm{R}}
\newcommand*{\fp}{\mathrm{F}}
\newcommand*{\selfterm}{\mathrm{self}}

\newcommand*{\ftpair}{\rightleftharpoons}
\newcommand*{\wh}[1]{\widehat{#1}}
\newcommand*{\pft}[1]{\breve{#1}}

\newcommand*{\gE}{\gamma_\mathrm{E}}
\newcommand*{\gH}{\gamma_\mathrm{H}}
\newcommand*{\hgE}{\wh{\gamma}_\mathrm{E}}
\newcommand*{\hgH}{\wh{\gamma}_\mathrm{H}}

\newcommand*{\window}{w}
\newcommand*{\gridh}{h}
\newcommand*{\rc}{r_\mathrm{c}}
\newcommand*{\gridfactor}{f_M}

\newcommand*{\Erms}{E_\mathrm{rms}}
\newcommand*{\Erel}{E_\mathrm{rms,rel}}
\newcommand*{\abstol}{\tau_\mathrm{rms}}
\newcommand*{\reltol}{\tau_\mathrm{rms,rel}}
\newcommand*{\Etrunc}{E_\mathrm{rms,trunc}}
\newcommand*{\Ewindow}{E_\mathrm{rms,window}}
\newcommand*{\Ewindowrel}{E_\mathrm{rms,window,rel}}
\newcommand*{\potmag}{U}

\newcommand*{\Matlab}{\textsc{Matlab}}

\usepackage{algorithm,algpseudocode}

\algnewcommand\algorithmicinput{\textbf{Input:}}
\algnewcommand\INPUT{\item[\algorithmicinput]}
\algnewcommand\algorithmicoutput{\textbf{Output:}}
\algnewcommand\OUTPUT{\item[\algorithmicoutput]}

\makeatletter
\newenvironment{breakablealgorithm}
  {
   \begin{center}
     \refstepcounter{algorithm}
     \hrule height.8pt depth0pt \kern2pt
     \renewcommand{\caption}[2][\relax]{
       {\raggedright\textbf{\ALG@name~\thealgorithm} ##2\par}%
       \ifx\relax##1\relax 
         \addcontentsline{loa}{algorithm}{\protect\numberline{\thealgorithm}##2}%
       \else 
         \addcontentsline{loa}{algorithm}{\protect\numberline{\thealgorithm}##1}%
       \fi
       \kern2pt\hrule\kern2pt
     }
  }{
     \kern2pt\hrule\relax
   \end{center}
  }
\makeatother

\usepackage{booktabs}
\usepackage[colorlinks]{hyperref}
\usepackage{placeins}

\journal{Journal of Computational Physics}

\begin{document}

\begin{frontmatter}

\title{Fast Ewald summation for Stokes flow with arbitrary periodicity}

\author{Joar Bagge\corref{corresponding}}
\ead{joarb@kth.se}
\author{Anna-Karin Tornberg}
\ead{akto@kth.se}
\address{KTH Mathematics, Linn\'{e} FLOW Centre/Swedish e-Science Research Centre,\\
Royal Institute of Technology, SE-100 44 Stockholm, Sweden}
\cortext[corresponding]{Corresponding author}

\begin{abstract}
  A fast and spectrally accurate Ewald summation method for the
  evaluation of stokeslet, stresslet and rotlet potentials of
  three-dimensional Stokes flow is presented. This work extends
  the previously developed Spectral Ewald method for Stokes flow
  to periodic boundary conditions in any number (three, two, one,
  or none) of the spatial directions, in a unified framework. The
  periodic potential is split into a short-range and a long-range
  part, where the latter is treated in Fourier space using the
  fast Fourier transform.
  A crucial component of the method is the modified kernels used
  to treat singular integration. We derive new modified kernels,
  and new improved truncation error estimates for the stokeslet
  and stresslet. An automated procedure for selecting parameters
  based on a given error tolerance is designed and tested.
  Analytical formulas for validation in the doubly and singly
  periodic cases are presented. We show that the computational
  time of the method scales like $O(N \log N)$ for $N$ sources and targets, and
  investigate how the time depends on the error tolerance and
  window function, i.e.\ the function used to smoothly spread
  irregular point data to a uniform grid. The method is fastest in the fully periodic
  case, while the run time in the free-space case is around three
  times as large.
  Furthermore, the highest efficiency is reached when applying
  the method to a uniform source distribution in a primary cell
  with low aspect ratio. The work presented in this paper enables
  efficient and accurate simulations of three-dimensional Stokes
  flow with arbitrary periodicity using e.g.\ boundary integral
  and potential methods.
\end{abstract}


\end{frontmatter}


\section{Introduction}
\label{sec:introduction}

Stokes flow, also known as creeping flow or viscous flow, is a
model of fluid flow in which inertial forces are assumed to be
negligible in comparison to viscous forces (i.e., the Reynolds
number is very small). This is often a valid assumption for
phenomena involving suspension flows on the micro- and nanoscales
(often in combination with Brownian motion), such as swimming
microorganisms \citep{Guasto2012}, cell dynamics
\citep{Maxian2021c}, microfluidic devices \citep{Squires2005},
gels \citep{Whitaker2019,Das2021}, dynamics of nanoparticles and
nanofibrils \citep{Mittal2018,Sherman2019,Turetta2022,Alcazar-Cano2022},
electrolytes \citep{Ladiges2021}, and antibodies \citep{Lai2021}.
For Stokes flow, the Navier--Stokes equations reduce to the
Stokes equations, which for an incompressible Newtonian fluid are
given by
\begin{gather}
  \label{eq:stokes-equations-1}
  - \nabla p(\bvec{x}) + \mu \nabla^2 \bvec{u}(\bvec{x}) + \bvec{g}(\bvec{x}) = \bvec{0},
  \\
  \nabla \cdot \bvec{u}(\bvec{x}) = 0.
\end{gather}
Here, $p$ is the pressure, $\bvec{u}$ is the fluid velocity,
$\bvec{g}$ is the body force per unit volume acting on the fluid,
and $\mu$ is the viscosity of the fluid. We will here consider
the nondimensionalized Stokes equations, which is equivalent to
setting $\mu=1$.

In boundary integral and potential methods for Stokes flow,
the fundamental solution (Green's function) of the Stokes equations
appears, namely the stokeslet kernel. In the three-dimensional
case, the stokeslet is a $3 \times 3$ tensor given by
\begin{equation}
  \stokeslet_{jl}(\bvec{r})
  = \frac{\delta_{jl}}{\lvert \bvec{r} \rvert}
  + \frac{r_j r_l}{\lvert \bvec{r} \rvert^3},
\end{equation}
where $\delta_{jl}$ is the Kronecker delta. Also commonly used
are the stresslet and rotlet kernels, which can be seen as
derivatives of the stokeslet and will be introduced in
section~\ref{sec:ewald-summation}. In this paper, we are
interested in problems with periodic boundary conditions, in
which $D$ of the three spatial directions will be periodic, and
the remaining $3-D$ directions will be free, in this context
meaning that the domain extends to infinity and that no boundary
conditions are enforced in these directions in the summation procedure.
We will call the problem triply, doubly and singly periodic if $D=3,2,1$,
respectively, and free-space if $D=0$. We will consider a system
of $N$ point sources of strengths $\bvec{f}(\bvec{x}_n)$ located
at positions $\bvec{x}_n$. The velocity field generated by this
system is given as a periodic sum over the point sources, i.e.
\begin{equation}
  \label{eq:periodic-sum}
  \bvec{u}^{D\per}(\bvec{x}) =
  \sum_{n=1}^N \sum_{\bvec{p} \in P_{D\per}}
  \bmat{\stokeslet}(\bvec{x} - \bvec{x}_n + \bvec{p})
  \bvec{f}(\bvec{x}_n),
\end{equation}
where the set $P_{D\per}$ is a $D$-dimensional lattice containing
all periodic images, to be properly defined in section~\ref{sec:ewald-summation}.
We assume that we want to evaluate \eqref{eq:periodic-sum} also
in $N$ target points, which may or may not be the same as the source
points. The periodic sum can be computed using Ewald summation,
which splits the sum into a short-range part, to be summed
directly, and a smooth long-range part, to be treated in Fourier
space. For the stokeslet, the split was derived by Hasimoto
\cite{Hasimoto1959}. A decomposition parameter $\xi$, to be introduced in
section~\ref{sec:ewald-summation}, controls the Ewald
summation split. For fixed $\xi$, computing the Ewald sums
directly leads to a method with time complexity $O(N^2)$.
Adjusting $\xi$ properly reduces the complexity of the direct
summation to $O(N^{3/2})$, see e.g.\ \cite{Darden1993,Karasawa1989}.
Fast Ewald summation methods such as
the Particle--Mesh--Ewald (PME) method \citep{Darden1993} and
smooth Particle--Mesh--Ewald (SPME) method
\citep{Essmann1995,Saintillan2005} reduce
the complexity further to $O(N \log N)$, a significant improvement.

The Spectral Ewald (SE) method is a PME-type method with spectral
accuracy. As in all PME methods, the interactions between source
points and target points are computed via a uniform grid and the
fast Fourier transform (FFT). A so-called window function is used
to interpolate between source/target points and the uniform grid,
much like in the nonuniform fast Fourier transform (NUFFT)
\citep{Dutt1993,Greengard2004}.
The spectral accuracy of the SE method comes from the choice of
window function, which was originally a truncated Gaussian. In
contrast, other methods such as e.g.\ SPME
\citep{Essmann1995,Saintillan2005} use Cardinal
$B$-splines for interpolation, which leads to algebraic accuracy.

The SE method for Stokes flow has been developed in a series of
papers, first for the triply periodic stokeslet \citep{Lindbo2010},
doubly periodic stokeslet \citep{Lindbo2011b}, triply periodic
stresslet \citep{afKlinteberg2014}, and triply periodic rotlet
\citep{afKlinteberg2016}. The method was extended to the
free-space case for all three kernels by \cite{afKlinteberg2017}.
The current paper serves to complete this development, much in
the same way as was recently done for electrostatics by
\cite{Shamshirgar2021}, by adding the missing pieces (singly
periodic case for all three kernels, and doubly periodic case for
the stresslet and rotlet), and unifying all periodic cases within
a single framework. This opens up the possibility to perform
efficient simulations of three-dimensional Stokes flow with
arbitrary periodicity ($D=3,2,1,0$), using boundary integral and
potential methods.

The SE method has also been adapted to related models, such as
Brinkman flow \citep{Nguyen2016}, and Stokesian dynamics
\citep{Wang2016}. To facilitate the inclusion of Brownian motion,
\cite{Fiore2017} proposed the Positively Split Ewald (PSE) method
for the Rotne--Prager--Yamakawa (RPY) tensor of Stokesian
dynamics, thus ensuring that both the short-range and long-range
parts are symmetric positive definite. The PSE method has found
much use in Brownian dynamics
\citep{Sprinkle2017,Fiore2018,Bao2018,Whitaker2019,Sherman2019,Das2021,Lai2021,Maxian2021a,Maxian2021c,Turetta2022,Alcazar-Cano2022}.

In this paper, the aperiodic directions are assumed to be free
and extend to infinity. If they are instead bounded, one
possibility is to explicitly discretize the boundary, and enforce
boundary conditions on it, as done e.g.\ in \cite{Bagge2021},
without modifying the underlying SE method. Another option is the
general geometry Ewald-like method (GGEM)
\citep{HernandezOrtiz2007,Zhao2017}, which uses the same split
into a short-range and long-range part, but treats the long-range
part in real space using a mesh-based solution method. GGEM can
in principle handle nontrivial boundary conditions, but involves
an expensive correction solve if high accuracy is desired. A more
recent solution is presented by \cite{Maxian2021b}, in which a
Chebyshev method is used for the boundary value problem in the
aperiodic direction; this method has been demonstrated for
electrostatics but is expected to generalize to Stokes flow.
Related work on Ewald-type methods for reduced periodicity in
electrostatics includes \cite{Nestler2015,Nestler2016,Weeber2019}.
For half-space Stokes problems, the methods by \cite{Srinivasan2018} or
\cite{Yan2018b} can be used. The SE method has also been
implemented for two-dimensional Stokes flow \citep{Paalsson2020}.

An alternative to Ewald-like summation methods is the fast
multipole method (FMM) \citep{Greengard1987,Fu2000,Wang2007,Tornberg2008}, which in
general achieves $O(N)$ complexity. Unlike Ewald-type methods,
which reach their highest efficiency for fully periodic problems,
the FMM is most efficient and most natural to formulate in the
fully aperiodic (free-space) setting. Nevertheless, the FMM has
also been generalized to arbitrary periodicity
\citep{Kabadshow2012,Yan2018,Yan2021}. An advantage of the FMM is
that it is spatially adaptive, while the SE method requires a
uniform spatial grid due to the FFT. Thus, the FMM will typically
be faster for highly nonuniform point distributions (especially
in free space), while the SE method may be faster for uniform
distributions, as shown e.g.\ by \cite{afKlinteberg2017,Shamshirgar2021}.
Yet another alternative method is found in \cite{Barnett2018}, in
which the long-range interaction is represented by auxiliary
sources. The principle is the same in both two and three
dimensions, but it has not been demonstrated that this method
would be competitive with Ewald-type methods in three dimensions.

The contribution of the current paper is, as mentioned above, to
complete and unify the previous work on the SE method for the
kernels of Stokes flow (stokeslet, stresslet, rotlet), by adding
the singly and doubly periodic cases, and treat
all periodic cases within the same framework. An important part
of the unified SE method is the modified kernels that are used to
treat the singular integration in cases with reduced periodicity,
based on an idea by \cite{Vico2016}. In this paper, we improve
the convergence of the modified kernels used by
\cite{afKlinteberg2017} in the free-space case, and derive new
modified kernels for the singly and doubly periodic cases. We
also derive a new improved truncation error estimate for the
stokeslet and stresslet, valid in all periodic cases, based on
techniques from \cite{afKlinteberg2017}. Analytical formulas
useful for validation are derived in the singly and doubly
periodic cases, completing the formulas previously derived by
\cite{Lindbo2011b} for the doubly periodic stokeslet. The SE
method presented in this paper furthermore uses the polynomial
Kaiser--Bessel (PKB) window introduced by \cite{Barnett2019} and \cite{Shamshirgar2021},
and the adaptive Fourier transform (AFT) introduced by
\cite{Shamshirgar2017} for the singly and doubly periodic cases.

The paper is organized as follows. In
section~\ref{sec:ewald-summation}, the Ewald splits and Ewald
sums are given for the three kernels. In
section~\ref{sec:spectral-ewald}, the fast method to compute the
Fourier-space Ewald sums is presented, i.e.\ the SE method; this
section also describes the modified kernels, PKB window, and AFT.
In section~\ref{sec:estimates}, we give error estimates for the
SE method and describe an automated procedure to select the
parameters of the method given an error tolerance; this procedure
is also tested by numerical examples. In
section~\ref{sec:results}, more numerical results follow,
regarding the pointwise error, computational time and complexity,
and the window function. Finally, conclusions are drawn in
section~\ref{sec:conclusions}. In the appendices, we derive
analytical formulas for validation in \ref{app:2p-integrals} and
\ref{app:1p-integrals}, and the improved truncation error
estimate for the stokeslet and stresslet in
\ref{app:truncation-estimates}.

\section{Ewald summation for Stokes flow}
\label{sec:ewald-summation}

We consider three fundamental solutions of Stokes flow, namely
the stokeslet~$\bmat{\stokeslet}$, rotlet~$\rotlett$, and
stresslet~$\bmat{\stresslet}$, which are tensorial kernels given by
\begin{equation}
  \label{eq:stokes-kernels}
  \stokeslet_{jl}(\bvec{r})
  = \frac{\delta_{jl}}{\lvert \bvec{r} \rvert}
  + \frac{r_j r_l}{\lvert \bvec{r} \rvert^3},
  \qquad
  \rotlet_{jl}(\bvec{r})
  = \epsilon_{jlm} \frac{r_m}{\lvert \bvec{r} \rvert^3},
  \qquad
  \stresslet_{jlm}(\bvec{r})
  = -6 \frac{r_j r_l r_m}{\lvert \bvec{r} \rvert^5}
  ,
\end{equation}
where $j,l,m \in \{1,2,3\}$.
Here, $\delta_{jl}$ denotes the Kronecker delta and
$\epsilon_{jlm}$ denotes the Levi-Civita symbol; Einstein's
summation convention is used, with repeated indices implicitly
summed over $\{1,2,3\}$. Given a point force $8\pi\bvec{f}$ acting on
the fluid at $\bvec{x}_0$, the Stokesian velocity field is given
by $u_j(\bvec{x}) = \stokeslet_{jl}(\bvec{x}-\bvec{x}_0) f_l$, the vorticity is
given by $\omega_j(\bvec{x}) = 2\rotlet_{jl}(\bvec{x}-\bvec{x}_0) f_l$, and
the stress field is given by $\sigma_{jl}(\bvec{x}) =
\stresslet_{jlm}(\bvec{x}-\bvec{x}_0) f_m$
\citep[ch.~2.2]{Pozrikidis1992}.
In boundary integral equations, the kernels $\bmat{\stokeslet}$,
$\rotlett$ and $\bmat{\stresslet}$ are multiplied by different source
quantities which all give rise to velocity fields, namely
\begin{align}
  u_j^\stokeslet(\bvec{x})
  &= \stokeslet_{jl}(\bvec{x} - \bvec{x}_0)
  f_l^\stokeslet
  , \\
  u_j^\rotlet(\bvec{x})
  &= \rotlet_{jl}(\bvec{x} - \bvec{x}_0)
  f_l^\rotlet
  , \\
  u_j^\stresslet(\bvec{x})
  &= \stresslet_{jlm}(\bvec{x} - \bvec{x}_0)
  f_{lm}^\stresslet
  ,
\end{align}
where notably the stresslet source $f_{lm}^\stresslet$ has two
indices. Commonly, the stresslet source is of the form
$f_{lm}^\stresslet = q_l \nu_m$, where $\bvec{q}$ and
$\bvec{\nu}$ are vectors, and this will be assumed in this paper.

In the following, we will make frequent use of the fact that the
fundamental solutions of Stokes flow can be related to the
harmonic (Laplace) Green's function $\harmonic(\bvec{r}) = 1/\lvert
\bvec{r} \rvert$ or the biharmonic
Green's function $\biharmonic(\bvec{r}) = \lvert \bvec{r} \rvert$, via
\citep[Appendix~E, p.~113]{Claeys1991}\citep{afKlinteberg2017}\citep{Fan1998}
\begin{align}
  \label{eq:diffop-rel-1}
  \stokeslet_{jl}(\bvec{r}) &=
  (\delta_{jl} \nabla^2 - \nabla_j \nabla_l) \biharmonic(\bvec{r})
  =: \Kop^{\stokeslet}_{jl} \biharmonic(\bvec{r})
  , \\
  \label{eq:diffop-rel-2}
  \rotlet_{jl}(\bvec{r}) &=
  - \epsilon_{jlm} \nabla_m \harmonic(\bvec{r})
  =: \Kop^{\rotlet}_{jl} \harmonic(\bvec{r})
  , \\
  \label{eq:diffop-rel-3}
  \stresslet_{jlm}(\bvec{r}) &=
  [(\delta_{jl} \nabla_m + \delta_{mj} \nabla_l + \delta_{lm}\nabla_j)
  \nabla^2 - 2 \nabla_j \nabla_l \nabla_m]
  \biharmonic(\bvec{r})
  =: \Kop^{\stresslet}_{jlm} \biharmonic(\bvec{r})
  .
\end{align}
We have here introduced a linear differential operator $\Kopt$ for
each kernel.
To be able to treat all kernels together where appropriate, we
now introduce some special notation following
\cite{afKlinteberg2017}; we will write
\begin{equation}
  \label{eq:generic-kernel}
  \bvec{u}(\bvec{x}) = \bmat{G}(\bvec{x}-\bvec{x}_0) \matop \bmat{f},
\end{equation}
where the kernel $\bmat{G}$ may be the stokeslet~$\bmat{\stokeslet}$,
rotlet~$\rotlett$ or stresslet~$\bmat{\stresslet}$, and the
notation $\bmat{G} \matop \bmat{f}$ will be understood to mean one of
\begin{equation}
  \label{eq:tensor-products}
  \stokeslet_{jl} f_l, \qquad
  \rotlet_{jl} f_l, \qquad
  \stresslet_{jlm} f_{lm},
\end{equation}
depending on the actual kernel. (Note that $\bmat{G} \matop
\bmat{f}$ is always a vector quantity.)

We are interested in efficiently evaluating periodic potentials
of the form
\begin{equation}
  \label{eq:periodic-potential}
  \bvec{u}^{D\per}(\bvec{x})
  = \sum_{n=1}^N \sum_{\bvec{p} \in P_{D\per}}
  \kernel(\bvec{x} - \bvec{x}_n + \bvec{p}) \matop \bmat{f}(\bvec{x}_n),
\end{equation}
generated by $N$ sources of strengths $\bmat{f}(\bvec{x}_n)$
located at points $\bvec{x}_n$ inside a box $\boxvar = [0,L_1)
\times [0,L_2) \times [0,L_3)$, called the primary cell.
The set $P_{D\per}$ of periodic images depends on the number of
periodic directions $D$ and is given by
\begin{equation}
  \label{eq:periodic-images}
  P_{D\per} = \begin{cases}
    \{ (\isc{p}_1 L_1, \isc{p}_2 L_2, \isc{p}_3 L_3) : \isc{p}_i
    \in \mathbb{Z} \}, & \text{if $D=3$ (triply periodic)}, \\
    \{ (\isc{p}_1 L_1, \isc{p}_2 L_2, 0) : \isc{p}_i
    \in \mathbb{Z} \}, & \text{if $D=2$ (doubly periodic)}, \\
    \{ (\isc{p}_1 L_1, 0, 0) : \isc{p}_i
    \in \mathbb{Z} \}, & \text{if $D=1$ (singly periodic)}, \\
    \{ (0, 0, 0) \}, & \text{if $D=0$ (free space)}.
  \end{cases}
\end{equation}
Note that the first $D$ coordinate directions are the periodic
ones; the remaining $3-D$ directions are called the free
directions. (The notation $D\per$ introduced here is to be read as
``$D$-periodic''.)
The sum \eqref{eq:periodic-potential} may be evaluated either at one of
the source locations $\bvec{x}_m$ ($m=1,\ldots,N$), in which case the term
corresponding to $m=n$, $\bvec{p}=\bvec{0}$ is omitted, or at
other arbitrary locations. In the former case, we write
\begin{equation}
  \label{eq:periodic-potential-star}
  \bvec{u}^{D\per\star}(\bvec{x}_m)
  = \sum_{n=1}^N \sum_{\bvec{p} \in P_{D\per}}^{\star}
  \kernel(\bvec{x}_m - \bvec{x}_n + \bvec{p}) \matop \bmat{f}(\bvec{x}_n),
\end{equation}
where the star above the second sum denotes precisely that the
term $\bvec{p}=\bvec{0}$ is omitted when $n=m$.

Since the kernels \eqref{eq:stokes-kernels} decay slowly as
$\lvert \bvec{r} \rvert \to \infty$ (the stokeslet decays like
$1/\lvert \bvec{r} \rvert$, and the rotlet and stresslet like
$1/\lvert \bvec{r} \rvert^2$), the periodic sums
\eqref{eq:periodic-potential} and
\eqref{eq:periodic-potential-star} are absolutely convergent only
in the case $D=0$ (for all kernels), and $D=1$ for the rotlet and
stresslet (but not the stokeslet). In all other cases, the sums
are either conditionally convergent (i.e.\ their values depend on
the order of summation) or even divergent. Special care must then
be taken when interpreting \eqref{eq:periodic-potential} and
\eqref{eq:periodic-potential-star}; an overview of these
considerations is given by \cite{Pozrikidis1996}. Ewald summation
corresponds to a spherical order of summation.

Ewald summation is based on the idea to split the periodic sum
into two parts, where the first part converges fast in real space
(the ``real-space part''), and the other part converges fast in
Fourier space (the ``Fourier-space part''). Ewald decompositions
for the kernels considered here are given in
section~\ref{sec:ewald-summation-decomp}. To aid the reader, we
then present the Ewald sums first in the triply periodic setting
($D=3$) in section~\ref{sec:ewald-summation-3p}, and then for
arbitrary periodicity ($D=3,2,1,0$) in section~\ref{sec:ewald-summation-Dp}.

\subsection{Ewald decompositions}
\label{sec:ewald-summation-decomp}

In general, there are two different but equivalent ways to derive
an Ewald decomposition, called screening and splitting
\citep{afKlinteberg2017}. The real-space part is easier to derive
using the splitting approach, while the Fourier-space part is
easier using the screening approach. In this paper, we use the
screening approach since the Fourier-space part is our main
focus.

In the screening approach, the kernel is convolved with a
screening function $\gamma(\bvec{r};\xi)$, which is a smooth function with the
property $\int_{\mathbb{R}^3} \gamma(\bvec{r};\xi) \,
\mathrm{d}\bvec{r} = 1$. The positive parameter $\xi$ is called
the Ewald decomposition parameter, and controls the decay of the
resulting decomposition
\begin{equation}
  \kernel = \kernel * (\delta - \gamma) + \kernel * \gamma,
\end{equation}
where $\kernel^\rp := \kernel * (\delta - \gamma)$ is the
real-space part and $\kernel^\fp := \kernel * \gamma$ is the
Fourier-space part; here, $\delta$ is the Dirac delta
distribution, and $*$ denotes convolution. In our case where $\kernel = \Kopt A$, with
$\Kopt$ being a linear differential operator and $A$ either the
harmonic~$\harmonic$ or biharmonic~$\biharmonic$, cf.~\eqref{eq:diffop-rel-1}--\eqref{eq:diffop-rel-3},
the decomposition can be written as $\kernel = \kernel^\rp +
\kernel^\fp$ with
\begin{align}
  \kernel^\rp &= \Kopt [A * (\delta - \gamma)] =: \Kopt A^\rp, \\
  \kernel^\fp &= \Kopt [A * \gamma] =: \Kopt A^\fp.
\end{align}
Thus, we start by writing down the Ewald decompositions for the
harmonic and biharmonic kernels, and may then apply the appropriate
differential operator $\Kopt$ to get the decompositions for the
stokeslet, rotlet, and stresslet.

For the harmonic, we select the classical Ewald screening
function
\begin{equation}
  \label{eq:ewald-screening}
  \gE(\bvec{r};\xi) = \xi^3 \pi^{-3/2} \E^{-\xi^2 \lvert \bvec{r} \rvert^2}
  \quad \ftpair \quad
  \hgE(\bvec{k};\xi) = \E^{-\lvert \bvec{k} \rvert^2 / (2\xi)^2},
\end{equation}
where the hat denotes the Fourier transform
\begin{equation}
  \wh{f}(\bvec{k}) = \fourier\{f\}(\bvec{k}) =
  \int_{\mathbb{R}^3} f(\bvec{r}) \E^{-\I \bvec{k} \cdot \bvec{r}}
  \, \mathrm{d}\bvec{r}
  .
\end{equation}
This leads to the decomposition originally derived by Ewald
\cite{Ewald1921}, namely
\begin{align}
  \harmonic^\rp(\bvec{r};\xi)
  &= [\harmonic * (\delta-\gE)](\bvec{r};\xi)
  = \frac{\erfc(\xi \lvert \bvec{r} \rvert)}{\lvert \bvec{r} \rvert}
  ,
  \\
  \harmonic^\fp(\bvec{r};\xi)
  &= [\harmonic * \gE](\bvec{r};\xi)
  = \frac{\erf(\xi \lvert \bvec{r} \rvert)}{\lvert \bvec{r} \rvert}
  ,
\end{align}
where $\erf(\cdot)$ is the error function and $\erfc(\cdot) = 1 - \erf(\cdot)$ is
the complementary error function.
Here, the singular behaviour of $1/\lvert \bvec{r}\rvert$ is
contained in $\harmonic^\rp$, while the long-range behaviour is
contained in $\harmonic^\fp$. Note that $\harmonic^\rp$ decays
fast as $\lvert \bvec{r} \rvert \to \infty$ in real space, while
$\harmonic^\fp$ is smooth and therefore its Fourier transform
decays fast. By the convolution theorem, the Fourier transform
(in the distributional sense) of $\harmonic^\fp$ is
\begin{equation}
  \wh{\harmonic}^\fp(\bvec{k};\xi)
  = \wh{\harmonic}(\bvec{k}) \hgE(\bvec{k};\xi)
  = \frac{4\pi}{\lvert \bvec{k} \rvert^2}
  \E^{-\lvert \bvec{k} \rvert^2 / (2\xi)^2},
\end{equation}
since $\wh{\harmonic}(\bvec{k}) = 4\pi / \lvert \bvec{k} \rvert^2$.

For the biharmonic, the classical Ewald screening function does
not yield a rapidly decaying real-space part, so a different
screening function is needed. We will use the Hasimoto screening
function \citep{HernandezOrtiz2007,Hasimoto1959}
\begin{equation}
  \label{eq:hasimoto-screening}
  \gH(\bvec{r};\xi) = \xi^3 \pi^{-3/2} \E^{-\xi^2 \lvert \bvec{r} \rvert^2}
  \left( \frac{5}{2} - \xi^2 \lvert \bvec{r} \rvert^2 \right)
  \quad \ftpair \quad
  \hgH(\bvec{k};\xi) = \E^{-\lvert \bvec{k} \rvert^2 / (2\xi)^2}
  \left( 1 + \frac{\lvert \bvec{k} \rvert^2}{(2\xi)^2} \right)
  .
\end{equation}
(Another option is the Beenakker screening function, see e.g.\
\cite[Table~1]{afKlinteberg2017}, but the Hasimoto screening
function yields a somewhat faster decay in both real and Fourier
space, and is therefore preferred.) The Hasimoto screening leads
to the decomposition
\begin{align}
  \biharmonic^\rp(\bvec{r};\xi)
  &= [\biharmonic * (\delta-\gH)](\bvec{r};\xi)
  = \lvert \bvec{r} \rvert \erfc(\xi \lvert \bvec{r} \rvert)
  - \frac{\E^{-\xi^2 \lvert \bvec{r} \rvert^2}}{\sqrt{\pi}\xi}
  ,
  \\
  \biharmonic^\fp(\bvec{r};\xi)
  &= [\biharmonic * \gH](\bvec{r};\xi)
  = \lvert \bvec{r} \rvert \erf(\xi \lvert \bvec{r} \rvert)
  + \frac{\E^{-\xi^2 \lvert \bvec{r} \rvert^2}}{\sqrt{\pi}\xi}
  .
\end{align}
By the convolution theorem, we have
\begin{equation}
  \wh{\biharmonic}^\fp(\bvec{k};\xi)
  = \wh{\biharmonic}(\bvec{k}) \hgH(\bvec{k};\xi)
  = -\frac{8\pi}{\lvert \bvec{k} \rvert^4}
  \E^{-\lvert \bvec{k} \rvert^2 / (2\xi)^2}
  \left( 1 + \frac{\lvert \bvec{k} \rvert^2}{(2\xi)^2} \right),
\end{equation}
since $\wh{\biharmonic}(\bvec{k}) = -8\pi / \lvert \bvec{k} \rvert^4$.
Again, $\biharmonic^\rp$ decays fast in real space, while the
Fourier transform of $\biharmonic^\fp$ decays fast.

By applying the relations \eqref{eq:diffop-rel-1}--\eqref{eq:diffop-rel-3},
one obtains Ewald decompositions also for the stokeslet, rotlet,
and stresslet. The real-space part kernels $\kernel^\rp$ become
\begin{align}
  \label{eq:stokeslet-real}
  \stokeslet_{jl}^\rp(\bvec{r};\xi) &=
  \Kop_{jl}^\stokeslet \biharmonic^\rp(\bvec{r};\xi)
  = \left( \delta_{jl} + \frac{r_j r_l}{\lvert \bvec{r} \rvert^2} \right)
  \left(
    \frac{\erfc(\xi \lvert \bvec{r} \rvert)}{\lvert \bvec{r} \rvert}
    + \frac{2 \xi \E^{-\xi^2 \lvert \bvec{r} \rvert^2}}{\sqrt{\pi}}
  \right)
  - \delta_{jl} \frac{4\xi \E^{-\xi^2 \lvert \bvec{r} \rvert^2}}{\sqrt{\pi}}
  , \\
  \label{eq:rotlet-real}
  \rotlet_{jl}^\rp(\bvec{r};\xi) &=
  \Kop_{jl}^\rotlet \harmonic^\rp(\bvec{r};\xi)
  = \epsilon_{jlm} \frac{r_m}{\lvert \bvec{r} \rvert^2}
  \left(
    \frac{\erfc(\xi \lvert \bvec{r} \rvert)}{\lvert \bvec{r} \rvert}
    + \frac{2 \xi \E^{-\xi^2 \lvert \bvec{r} \rvert^2}}{\sqrt{\pi}}
  \right)
  , \\
  \notag
  \stresslet_{jlm}^\rp(\bvec{r};\xi) &=
  \Kop_{jlm}^\stresslet \biharmonic^\rp(\bvec{r};\xi)
  = - \frac{2 r_j r_l r_m}{\lvert \bvec{r} \rvert^4}
  \left(
    \frac{3\erfc(\xi \lvert \bvec{r} \rvert)}{\lvert \bvec{r} \rvert}
    + (3 + 2\xi^2 \lvert \bvec{r} \rvert^2)
    \frac{2\xi \E^{-\xi^2 \lvert \bvec{r} \rvert^2}}{\sqrt{\pi}}
  \right)
  \\
  \label{eq:stresslet-real}
  &\hspace*{3cm}+ (\delta_{jl} r_m + \delta_{mj} r_l + \delta_{lm} r_j)
  \frac{4\xi^3 \E^{-\xi^2 \lvert \bvec{r} \rvert^2}}{\sqrt{\pi}}
  .
\end{align}
For the Fourier-space parts, we rewrite the relation
$\kernel = \Kopt A$, cf.~\eqref{eq:diffop-rel-1}--\eqref{eq:diffop-rel-3},
as $\wh{\kernel} = \wh{\Kopt} \wh{A}$, where the
Fourier-space operator $\wh{\Kopt}$ is found by letting the
differential operator $\Kopt$ act on $\E^{\I \bvec{k} \cdot \bvec{r}}$.
We have
\begin{align}
  \label{eq:diffop-rel-fourier-1}
  \wh{\Kop}^{\stokeslet}_{jl}(\bvec{k})
  &= -\delta_{jl} \lvert \bvec{k} \rvert^2 + k_j k_l
  , \\
  \label{eq:diffop-rel-fourier-2}
  \wh{\Kop}^{\rotlet}_{jl}(\bvec{k})
  &= -\epsilon_{jlm} \I k_m
  , \\
  \label{eq:diffop-rel-fourier-3}
  \wh{\Kop}^{\stresslet}_{jlm}(\bvec{k})
  &=
  -\I (\delta_{jl} k_m + \delta_{mj} k_l + \delta_{lm} k_j)
  \lvert \bvec{k} \rvert^2 + 2 \I k_j k_l k_m
  .
\end{align}
The Fourier transform of the Fourier-space part kernel then becomes
$\wh{\kernel}{}^\fp(\bvec{k};\xi) = \wh{\Kopt}(\bvec{k})
\wh{A}(\bvec{k}) \wh{\gamma}(\bvec{k};\xi)$. For the three
kernels of Stokes flow,
\begin{align}
  \label{eq:stokeslet-fourier}
  \wh{\stokeslet}_{jl}^\fp(\bvec{k};\xi) &=
  \wh{\Kop}_{jl}^\stokeslet(\bvec{k}) \wh{\biharmonic}(\bvec{k})
  \hgH(\bvec{k};\xi)
  =
  \left( -\delta_{jl} \lvert \bvec{k} \rvert^2 + k_j k_l \right)
  \left( -\frac{8\pi}{\lvert \bvec{k} \rvert^4} \right)
  \E^{-\lvert \bvec{k} \rvert^2 / (2\xi)^2}
  \left( 1 + \frac{\lvert \bvec{k} \rvert^2}{(2\xi)^2} \right)
  , \\
  \label{eq:rotlet-fourier}
  \wh{\rotlet}_{jl}^\fp(\bvec{k};\xi) &=
  \wh{\Kop}_{jl}^\rotlet(\bvec{k}) \wh{\harmonic}(\bvec{k})
  \hgE(\bvec{k};\xi)
  =
  -\epsilon_{jlm} \I k_m
  \frac{4\pi}{\lvert \bvec{k} \rvert^2}
  \E^{-\lvert \bvec{k} \rvert^2 / (2\xi)^2}
  , \\
  \notag
  \wh{\stresslet}_{jlm}^\fp(\bvec{k};\xi) &=
  \wh{\Kop}_{jlm}^\stresslet(\bvec{k}) \wh{\biharmonic}(\bvec{k})
  \hgH(\bvec{k};\xi)
  =
  \left( -\I (\delta_{jl} k_m + \delta_{mj} k_l + \delta_{lm} k_j)
  \lvert \bvec{k} \rvert^2 + 2 \I k_j k_l k_m \right)
  \\
  \label{eq:stresslet-fourier}
  &\hspace*{4.1cm}\times
  \left( -\frac{8\pi}{\lvert \bvec{k} \rvert^4} \right)
  \E^{-\lvert \bvec{k} \rvert^2 / (2\xi)^2}
  \left( 1 + \frac{\lvert \bvec{k} \rvert^2}{(2\xi)^2} \right)
  .
\end{align}
The Ewald sums are obtained by inserting $\kernel =
\kernel^\rp + \kernel^\fp$ into
\eqref{eq:periodic-potential} or
\eqref{eq:periodic-potential-star}, treating the sum associated
with $\kernel^\rp$ in real space and the one associated with
$\kernel^\fp$ in Fourier space. We do this below first in the
triply periodic case, and then for arbitrary periodicity.

\subsection{Triply periodic Ewald summation}
\label{sec:ewald-summation-3p}

Consider the case with $D=3$ periodic directions. Let us first
consider \eqref{eq:periodic-potential}, where the evaluation
point $\bvec{x}$ does not coincide with any of the source
locations $\bvec{x}_n$.
Inserting $\kernel = \kernel^\rp +
\kernel^\fp$ into \eqref{eq:periodic-potential} yields a
decomposition $\bvec{u}^{3\per} = \bvec{u}^{3\per,\rp} +
\bvec{u}^{3\per,\fp}$, where the real-space part is simply
\begin{equation}
  \label{eq:real-space-potential-3p}
  \bvec{u}^{3\per,\rp}(\bvec{x};\xi)
  = \sum_{n=1}^N \sum_{\bvec{p} \in P_{3\per}}
  \kernel^\rp(\bvec{x} - \bvec{x}_n + \bvec{p};\xi) \matop \bmat{f}(\bvec{x}_n).
\end{equation}
Note that since $\kernel^\rp$ decays fast, this periodic sum is
convergent.
For the Fourier-space part, we use the Poisson summation formula
to rewrite the periodic sum as
\begin{equation}
  \label{eq:fourier-space-potential-3p}
  \bvec{u}^{3\per,\fp,\bvec{k}\neq\bvec{0}}(\bvec{x};\xi)
  = \frac{1}{\lvert \boxvar \rvert} \sum_{n=1}^N
  \sum_{\substack{\bvec{k} \in \wavenumset^3 \\ \bvec{k} \neq \bvec{0}}}
  \wh{\kernel}{}^\fp(\bvec{k};\xi) \matop \bmat{f}(\bvec{x}_n)
  \E^{\I \bvec{k} \cdot (\bvec{x}-\bvec{x}_n)}
  ,
\end{equation}
where the box volume is $\lvert \boxvar \rvert = L_1 L_2 L_3$ and
the set of discrete wavenumbers is given by
\begin{equation}
  \label{eq:wavenumbers-3p}
  \wavenumset^3
  := \left\{ 2\pi
  \left(\frac{\isc{k}_1}{L_1}, \frac{\isc{k}_2}{L_2}, \frac{\isc{k}_3}{L_3}\right)
  : \isc{k}_i \in \mathbb{Z} \right\}
  .
\end{equation}
In \eqref{eq:fourier-space-potential-3p}, the term corresponding
to $\bvec{k}=\bvec{0}$ (the zero mode) has been omitted since
$\wh{\kernel}{}^\fp$ is singular at $\bvec{k}=\bvec{0}$. The zero
mode of a Fourier series corresponds to a constant
\begin{equation}
  \label{eq:fourier-space-potential-3p-k0}
  \bvec{u}^{3\per,\fp,\bvec{k}=\bvec{0}}
  = \frac{1}{\lvert \boxvar \rvert} \sum_{n=1}^N
  \bmat{C}^{(0)} \matop \bmat{f}(\bvec{x}_n)
  ,
\end{equation}
where $\bmat{C}^{(0)}$ is a constant tensor. The zero mode
\eqref{eq:fourier-space-potential-3p-k0} may
be added to \eqref{eq:fourier-space-potential-3p}, and the
freedom to select $\bmat{C}^{(0)}$ reflects that the fundamental
solutions are defined only up to a constant. To determine the
constant, one can e.g.\ impose a requirement of zero mean flow
\citep{afKlinteberg2014,afKlinteberg2016}, i.e.
\begin{equation}
  \label{eq:zero-mean-flow}
  \langle u_j^{3\per} \rangle :=
  \frac{1}{A_j} \int_{D_j} u_j^{3\per}(\bvec{x}) \,
  \mathrm{d}S(\bvec{x})
  = 0,
\end{equation}
where $D_j$ is the face of the box $\boxvar = [0,L_1) \times
[0,L_2) \times [0,L_3)$ lying in the plane $x_j=0$, and $A_j$ is
the area of $D_j$. (For instance, $D_1 = \{0\} \times [0,L_2)
\times [0,L_3)$ and $A_1=L_2 L_3$.) For the stokeslet and rotlet,
one can show \citep{afKlinteberg2016} that
\eqref{eq:zero-mean-flow} is satisfied when the
zero mode is zero, i.e.\ $\bvec{u}^{\stokeslet,3\per,\fp,\bvec{k}=\bvec{0}}
= \bvec{u}^{\rotlet,3\per,\fp,\bvec{k}=\bvec{0}} = \bvec{0}$.

For the stresslet, however, the situation is more complicated.
Setting the zero mode to zero will \emph{not} result in a zero
mean flow, as shown by \cite{afKlinteberg2014}. Yet another
perspective is provided by \cite[Paper~IV]{Marin2012}, where the
Ewald decomposition for the stresslet is derived using the
relation
\begin{equation}
  \label{eq:stresslet-physical-relation}
  \stresslet_{jlm}(\bvec{r}) =
  - \pressurelet_j(\bvec{r}) \delta_{lm}
  + \nabla_m \stokeslet_{jl}(\bvec{r})
  + \nabla_l \stokeslet_{mj}(\bvec{r}),
\end{equation}
where $\pressurelet_j(\bvec{r}) = 2 r_j / \lvert \bvec{r}
\rvert^3$ is the fundamental solution for the Stokesian pressure.
(The corresponding derivation in two dimensions was carried out
by \cite{VanDeVorst1996}.)
The resulting decomposition is the same as the one given here,
i.e.\ \eqref{eq:stresslet-real} and
\eqref{eq:stresslet-fourier}, for $\bvec{k} \neq \bvec{0}$.
However, it is noted that by setting $\bvec{u}^{\stokeslet,3\per,\fp,\bvec{k}=\bvec{0}}
= \bvec{0}$ for the stokeslet, a mean pressure gradient will appear (in
the Ewald sum for the pressure) to balance the point forces, represented by
$\nabla_l \pressurelet_j^{(0)}(\bvec{r}) = 8\pi \delta_{jl}/
\lvert \boxvar \rvert$. This implies that the pressure itself
gains a nonperiodic part
\begin{equation}
  \pressurelet_j^{(0)}(\bvec{r}) = \frac{8\pi}{\lvert \boxvar
  \rvert} r_j,
\end{equation}
which via \eqref{eq:stresslet-physical-relation} carries over to
the stresslet as
\begin{equation}
  \stresslet_{jlm}^{(0)}(\bvec{r}) = -\frac{8\pi}{\lvert \boxvar
  \rvert} r_j \delta_{lm}.
\end{equation}
Thus, the stresslet decomposition will gain the additional
nonperiodic term (here we abuse the notation for the zero mode
to represent this extra term, which is of course not constant)
\begin{equation}
  u_j^{\stresslet,3\per,\fp,\bvec{k}=\bvec{0}}(\bvec{x})
  = \sum_{n=1}^N \stresslet_{jlm}^{(0)}(\bvec{x}-\bvec{x}_n)
  f^\stresslet_{lm}(\bvec{x}_n)
  = -\frac{8\pi}{\lvert \boxvar \rvert}
  \sum_{n=1}^N (\bvec{x} - \bvec{x}_n)_j
  f^\stresslet_{ll}(\bvec{x}_n).
\end{equation}
Using the assumption that $f^\stresslet_{lm} = q_l \nu_m$, this
can be written as
\begin{equation}
  \label{eq:k0-stresslet-3p}
  \bvec{u}^{\stresslet,3\per,\fp,\bvec{k}=\bvec{0}}(\bvec{x})
  = -\frac{8\pi}{\lvert \boxvar \rvert}
  \sum_{n=1}^N (\bvec{x} - \bvec{x}_n)
  [\bvec{q}(\bvec{x}_n) \cdot \bvec{\nu}(\bvec{x}_n)].
\end{equation}
af Klinteberg and Tornberg \cite{afKlinteberg2014} derived a similar term under the
so-called rigid body assumption
\begin{equation}
  \label{eq:rigid-body}
  \sum_{n=1}^N [\bvec{q}(\bvec{x}_n) \cdot \bvec{\nu}(\bvec{x}_n)]=0,
\end{equation}
and our term \eqref{eq:k0-stresslet-3p} reduces to theirs under
the rigid body assumption \eqref{eq:rigid-body}.

To summarize, the full Ewald decomposition, for an evaluation
point $\bvec{x}$ which does not coincide with any of the source
locations, is given by
\begin{equation}
  \label{eq:3p-ewald-decomp}
  \bvec{u}^{3\per}(\bvec{x}) = \bvec{u}^{3\per,\rp}(\bvec{x};\xi)
  + \bvec{u}^{3\per,\fp,\bvec{k}\neq\bvec{0}}(\bvec{x};\xi)
  + \bvec{u}^{3\per,\fp,\bvec{k}=\bvec{0}}(\bvec{x}),
\end{equation}
where $\bvec{u}^{3\per,\rp}$ is given by
\eqref{eq:real-space-potential-3p},
$\bvec{u}^{3\per,\fp,\bvec{k}\neq\bvec{0}}$ is given by
\eqref{eq:fourier-space-potential-3p}, and
$\bvec{u}^{3\per,\fp,\bvec{k}=\bvec{0}}$ is zero for the
stokeslet and rotlet but given by \eqref{eq:k0-stresslet-3p} for
the stresslet.

When the evaluation point is one of the source locations
$\bvec{x}_m$ ($m=1,\ldots,N$), the term corresponding to $n=m$
for $\bvec{p} = \bvec{0}$ (the ``self interaction'') should be
omitted, cf.\ \eqref{eq:periodic-potential-star}. For the
real-space part, the self interaction can be omitted directly,
and the real-space part is thus given by
\eqref{eq:real-space-potential-3p} but with a star above the
second summation sign (denoted by $\bvec{u}^{3\per,\rp\star}$).
For the Fourier-space part, the self
interaction must be removed explicitly, which is done by adding
the term
\begin{equation}
  \label{eq:self-term}
  \bvec{u}^{\selfterm}(\bvec{x}_m;\xi) = -\kernel^\fp(\bvec{0};\xi)
  \matop \bmat{f}(\bvec{x}_m)
  = \lim_{\bvec{r} \to \bvec{0}} \left(
    \kernel^\rp(\bvec{r};\xi) - \kernel^\rp(\bvec{r})
  \right)
  \matop \bmat{f}(\bvec{x}_m)
  .
\end{equation}
For the stokeslet, we get
\begin{equation}
  \label{eq:3p-stokeslet-self-term}
  \bvec{u}^{\stokeslet,\selfterm}(\bvec{x}_m;\xi) =
  - \frac{4\xi}{\sqrt{\pi}} \bmat{f}(\bvec{x}_m),
\end{equation}
while for the rotlet and stresslet,
$\bvec{u}^{\rotlet,\selfterm}(\bvec{x}_m;\xi) =
\bvec{u}^{\stresslet,\selfterm}(\bvec{x}_m;\xi) =
\bvec{0}$. The full Ewald decomposition is then given by
\begin{equation}
  \label{eq:3p-ewald-decomp-self}
  \bvec{u}^{3\per\star}(\bvec{x}_m) = \bvec{u}^{3\per,\rp\star}(\bvec{x}_m;\xi)
  + \bvec{u}^{3\per,\fp,\bvec{k}\neq\bvec{0}}(\bvec{x}_m;\xi)
  + \bvec{u}^{3\per,\fp,\bvec{k}=\bvec{0}}(\bvec{x}_m)
  + \bvec{u}^{\selfterm}(\bvec{x}_m;\xi),
  \quad m=1,\ldots,N
  .
\end{equation}
Here, $\bvec{u}^{3\per,\rp\star}$ is given by
\eqref{eq:real-space-potential-3p} but with a star above the
second sum (to skip the term $\bvec{p}=\bvec{0}$ for $n=m$),
$\bvec{u}^{3\per,\fp,\bvec{k}\neq\bvec{0}}$ is given by
\eqref{eq:fourier-space-potential-3p},
$\bvec{u}^{3\per,\fp,\bvec{k}=\bvec{0}}$ is zero for the
stokeslet and rotlet but given by \eqref{eq:k0-stresslet-3p} for
the stresslet, and $\bvec{u}^{\selfterm}$ is zero for the rotlet
and stresslet but given by \eqref{eq:3p-stokeslet-self-term}
for the stokeslet. The expressions for $\kernel^\rp$ and
$\wh{\kernel}{}^\fp$ for the different kernels are found in
\eqref{eq:stokeslet-real}--\eqref{eq:stresslet-real} and
\eqref{eq:stokeslet-fourier}--\eqref{eq:stresslet-fourier},
respectively.
In summary, \eqref{eq:3p-ewald-decomp-self} is
used if the evaluation point coincides with one of the source
locations, while \eqref{eq:3p-ewald-decomp} is used otherwise.

\subsection{Ewald summation in arbitrary periodicity}
\label{sec:ewald-summation-Dp}

Let us now consider $D=3,2,1,0$ periodic directions.
When reducing the number of periodic directions, the real-space
part is straightforward to write down, and is given by
\begin{equation}
  \label{eq:real-space-potential-Dp}
  \bvec{u}^{D\per,\rp}(\bvec{x};\xi)
  = \sum_{n=1}^N \sum_{\bvec{p} \in P_{D\per}}
  \kernel^\rp(\bvec{x} - \bvec{x}_n + \bvec{p};\xi) \matop \bmat{f}(\bvec{x}_n),
\end{equation}
with a star above the second summation sign if $\bvec{x}$
coincides with one of the source locations $\bvec{x}_m$. Recall
that the set $P_{D\per}$, given by \eqref{eq:periodic-images},
ensures that the potential is periodically summed only in the
first $D$ directions. Since $\kernel^\rp$ decays fast, the sum
\eqref{eq:real-space-potential-Dp} is convergent, and it can be
computed efficiently for example using a cell list, as described
in section~\ref{sec:estimates-truncation-real}.
The term $\bvec{u}^\selfterm$, which removes
self interaction from the Fourier-space part and is defined by
\eqref{eq:self-term}, is independent of $D$.

In the Fourier-space part, cf.\ \eqref{eq:fourier-space-potential-3p},
the summation over discrete wavenumbers will in the free
directions be replaced by integration over continuous
wavenumbers (as the Fourier transform must be used instead of a
discrete Fourier series). Writing down the formulas explicitly
for all $D$ (at the moment ignoring the singularity at
$\bvec{k}=\bvec{0}$), we have
\begin{align}
  \label{eq:explicit-fourier-3p}
  \bvec{u}^{3\per,\fp}(\bvec{x};\xi)
  &= \frac{1}{L_1 L_2 L_3} \sum_{n=1}^N
  \sum_{\bvec{k} \in \wavenumset^3}
  \wh{\kernel}{}^\fp(\bvec{k};\xi) \matop \bmat{f}(\bvec{x}_n)
  \E^{\I \bvec{k} \cdot (\bvec{x}-\bvec{x}_n)}
  , \\
  \label{eq:explicit-fourier-2p}
  \bvec{u}^{2\per,\fp}(\bvec{x};\xi)
  &= \frac{1}{L_1 L_2 2\pi} \sum_{n=1}^N
  \sum_{(k_1,k_2) \in \wavenumset^2}
  \int_{\mathbb{R}}
  \wh{\kernel}{}^\fp(k_1,k_2,\kappa_3;\xi) \matop \bmat{f}(\bvec{x}_n)
  \E^{\I (k_1,k_2,\kappa_3) \cdot (\bvec{x}-\bvec{x}_n)}
  \, \mathrm{d}\kappa_3
  , \\
  \label{eq:explicit-fourier-1p}
  \bvec{u}^{1\per,\fp}(\bvec{x};\xi)
  &= \frac{1}{L_1 (2\pi)^2} \sum_{n=1}^N
  \sum_{k_1 \in \wavenumset^1}
  \int_{\mathbb{R}^2}
  \wh{\kernel}{}^\fp(k_1,\kappa_2,\kappa_3;\xi) \matop \bmat{f}(\bvec{x}_n)
  \E^{\I (k_1,\kappa_2,\kappa_3) \cdot (\bvec{x}-\bvec{x}_n)}
  \, \mathrm{d}\kappa_2 \, \mathrm{d}\kappa_3
  , \\
  \label{eq:explicit-fourier-0p}
  \bvec{u}^{0\per,\fp}(\bvec{x};\xi)
  &= \frac{1}{(2\pi)^3} \sum_{n=1}^N
  \int_{\mathbb{R}^3}
  \wh{\kernel}{}^\fp(\kappa_1,\kappa_2,\kappa_3;\xi) \matop \bmat{f}(\bvec{x}_n)
  \E^{\I (\kappa_1,\kappa_2,\kappa_3) \cdot (\bvec{x}-\bvec{x}_n)}
  \, \mathrm{d}\kappa_1 \, \mathrm{d}\kappa_2 \, \mathrm{d}\kappa_3
  ,
\end{align}
where the set of discrete wavenumbers is given by, cf.\
\eqref{eq:wavenumbers-3p},
\begin{equation}
  \label{eq:wavenumbers-Dp}
  \wavenumset^D
  := \left\{ 2\pi
  \left(\frac{\isc{k}_1}{L_1}, \cdots, \frac{\isc{k}_D}{L_D}\right)
  : \isc{k}_i \in \mathbb{Z} \right\}
  \subset \mathbb{R}^D
  .
\end{equation}
In \eqref{eq:explicit-fourier-3p}--\eqref{eq:explicit-fourier-0p},
we have written $\kappa_i$ instead of $k_i$ in the free directions
to emphasize that these wavenumbers are continuous.
Occasionally, we will use $\bvec{k}^\per$ to denote the vector of
discrete wavenumbers in the periodic directions, i.e.\
$\bvec{k}^\per \in \wavenumset^D$. Thus, $\bvec{k}^\per$ means
$(k_1,k_2,k_3)$ for $D=3$, $(k_1,k_2)$ for $D=2$, and $k_1$ for
$D=1$. Similarly, we may use $\bvec{\kappa} \in \mathbb{R}^{3-D}$
to mean $\kappa_3$ for $D=2$, $(\kappa_2,\kappa_3)$ for $D=1$ and
$(\kappa_1,\kappa_2,\kappa_3)$ for $D=0$. The notation
$\bvec{k}$ without superscript always refers to the full vector
of three wavenumbers regardless of $D$.

We can unify \eqref{eq:explicit-fourier-3p}--\eqref{eq:explicit-fourier-0p} by
introducing a ``mixed'' Fourier transform $\fourier_{D\per}$, as
follows. Let $f : \mathbb{R}^3 \to \mathbb{C}$ be a function that
is periodic in the first $D$ coordinate directions, and
nonperiodic in the remaining $3-D$ directions. We define the
mixed Fourier transform of $f$ by
\begin{equation}
  \label{eq:mixed-fourier-transform}
  \fourier_{D\per}\{f\}(\bvec{k})
  = \int_{\boxvar_{D\per}}
  f(\bvec{r}) \E^{-\I \bvec{k} \cdot \bvec{r}}
  \, \mathrm{d}\bvec{r},
\end{equation}
where the integration domain is
\begin{equation}
  \label{eq:mixed-fourier-integration-domain}
  \boxvar_{D\per} = \begin{cases}
    [0,L_1) \times [0,L_2) \times [0,L_3), & \text{if $D=3$}, \\
    [0,L_1) \times [0,L_2) \times \mathbb{R}, & \text{if $D=2$}, \\
    [0,L_1) \times \mathbb{R} \times \mathbb{R}, & \text{if $D=1$}, \\
    \mathbb{R} \times \mathbb{R} \times \mathbb{R}, & \text{if $D=0$}.
  \end{cases}
\end{equation}
The inverse transform is given by
\begin{equation}
  \label{eq:inverse-mixed-ft}
  \fourier_{D\per}^{-1} \{ g \}(\bvec{r})
  = \begin{cases}
    \displaystyle
    \frac{1}{L_1L_2L_3}
    \sum_{(k_1,k_2,k_3) \in \wavenumset^3}
    g(k_1,k_2,k_3)
    \E^{\I (k_1,k_2,k_3) \cdot \bvec{r}},
    & \text{if $D=3$},
    \\
    \displaystyle
    \frac{1}{L_1L_2 2\pi}
    \sum_{(k_1,k_2) \in \wavenumset^2}
    \int_{\mathbb{R}}
    g(k_1,k_2,\kappa_3)
    \E^{\I (k_1,k_2,\kappa_3) \cdot \bvec{r}}
    \, \mathrm{d}\kappa_3,
    & \text{if $D=2$},
    \\
    \displaystyle
    \frac{1}{L_1 (2\pi)^2}
    \sum_{k_1 \in \wavenumset^1}
    \int_{\mathbb{R}^2}
    g(k_1,\kappa_2,\kappa_3)
    \E^{\I (k_1,\kappa_2,\kappa_3) \cdot \bvec{r}}
    \, \mathrm{d}\kappa_2 \, \mathrm{d}\kappa_3,
    & \text{if $D=1$},
    \\
    \displaystyle
    \frac{1}{(2\pi)^3}
    \int_{\mathbb{R}^3}
    g(\kappa_1,\kappa_2,\kappa_3)
    \E^{\I (\kappa_1,\kappa_2,\kappa_3) \cdot \bvec{r}}
    \, \mathrm{d}\kappa_1 \, \mathrm{d}\kappa_2 \, \mathrm{d}\kappa_3,
    & \text{if $D=0$},
  \end{cases}
\end{equation}
and it holds that $f(\bvec{r}) = \fourier_{D\per}^{-1}\{
\fourier_{D\per}\{f \}\}(\bvec{r})$. Note that this is simply a
Fourier series in each periodic direction, and a Fourier
transform in each free direction. We can now write
\eqref{eq:explicit-fourier-3p}--\eqref{eq:explicit-fourier-0p}
compactly as
\begin{equation}
  \label{eq:fourier-potential-Dp}
  \bvec{u}^{D\per,\fp}(\bvec{x};\xi)
  =
  \sum_{n=1}^N
  \fourier_{D\per}^{-1}\left\{
  \wh{\kernel}{}^\fp \matop \bmat{f}(\bvec{x}_n)
  \right\} (\bvec{x}-\bvec{x}_n)
  .
\end{equation}

Let us consider the integrals that appear in the operator
$\fourier_{D\per}^{-1}$ in \eqref{eq:fourier-potential-Dp} for
$D=2,1,0$. Note that $\wh{\kernel}{}^\fp$, cf.\
\eqref{eq:stokeslet-fourier}--\eqref{eq:stresslet-fourier},
and thus the integrand, has a singularity at $\bvec{k}=\bvec{0}$.
This means that the integrand is singular for $D=2,1$ when
$\bvec{k}^\per=\bvec{0}$ (but not when $\bvec{k}^\per\neq\bvec{0}$),
as well as for $D=0$ (always). In fact, the integrals may not
even exist in the Lebesgue sense, but they can be interpreted as
inverse Fourier transforms in the distributional sense, as seen
in \ref{app:2p-integrals-k0} and \ref{app:1p-integrals-k0}.
For $D=2$ (see \ref{app:2p-integrals}), the integrals can be
evaluated analytically in both the nonsingular and singular
cases (in the distributional sense in the latter case), for all
three kernels. The results are of the form (with $\bvec{x}_n =
(x_n,y_n,z_n)$)
\begin{align}
  \label{eq:2p-evaluated-integrals-nonzero}
  \bvec{u}^{2\per,\fp,\bvec{k}^\per \neq \bvec{0}}
  (x,y,z;\xi)
  &=
  \frac{1}{L_1 L_2}
  \sum_{n=1}^N
  \sum_{\substack{(k_1,k_2)\in\wavenumset^2 \\
  (k_1,k_2) \neq (0,0)}}
  \bmat{Q}^{2\per}(k_1,k_2,z-z_n;\xi)
  \matop
  \bmat{f}(\bvec{x}_n)
  \E^{\I k_1 (x-x_n)}
  \E^{\I k_2 (y-y_n)}
  ,
  \\
  \label{eq:2p-evaluated-integrals-k0}
  \bvec{u}^{2\per,\fp,\bvec{k}^\per = \bvec{0}}
  (x,y,z;\xi)
  &=
  \frac{1}{L_1 L_2}
  \sum_{n=1}^N
  \bmat{Q}^{2\per,(0)}(z-z_n;\xi)
  \matop
  \bmat{f}(\bvec{x}_n)
  ,
\end{align}
where $\bmat{Q}^{2\per}$ and $\bmat{Q}^{2\per,(0)}$ are tensors
that depend on the kernel; $\bmat{Q}^{2\per}$ is given by
\eqref{eq:2p-rotlet-Q}, \eqref{eq:2p-stokeslet-Q}, \eqref{eq:2p-stresslet-Q},
and $\bmat{Q}^{2\per,(0)}$ by
\eqref{eq:2p-rotlet-Q0}, \eqref{eq:2p-stokeslet-Q0}, \eqref{eq:2p-stresslet-Q0}.

Also for $D=1$ (see \ref{app:1p-integrals}), the integrals can be
evaluated analytically. The results are of the form (again with
$\bvec{x}_n = (x_n, y_n, z_n)$)
\begin{align}
  \label{eq:1p-evaluated-integrals-nonzero}
  \bvec{u}^{1\per,\fp,\bvec{k}^\per \neq \bvec{0}}
  (x,y,z;\xi)
  &=
  \frac{1}{L_1}
  \sum_{n=1}^N
  \sum_{\substack{k_1\in\wavenumset^1 \\ k_1 \neq 0}}
  \bmat{Q}^{1\per}(k_1,y-y_n,z-z_n;\xi)
  \matop
  \bmat{f}(\bvec{x}_n)
  \E^{\I k_1 (x-x_n)}
  ,
  \\
  \label{eq:1p-evaluated-integrals-k0}
  \bvec{u}^{1\per,\fp,\bvec{k}^\per = \bvec{0}}
  (x,y,z;\xi)
  &=
  \frac{1}{L_1}
  \sum_{n=1}^N
  \bmat{Q}^{1\per,(0)}(y-y_n,z-z_n;\xi)
  \matop
  \bmat{f}(\bvec{x}_n)
  ,
\end{align}
where $\bmat{Q}^{1\per}$ is given by
\eqref{eq:1p-rotlet-Q}, \eqref{eq:1p-stokeslet-Q}, \eqref{eq:1p-stresslet-Q},
and $\bmat{Q}^{1\per,(0)}$ is given by
\eqref{eq:1p-rotlet-Q0}, \eqref{eq:1p-stokeslet-Q0}, \eqref{eq:1p-stresslet-Q0}.

The formulas given in
\eqref{eq:2p-evaluated-integrals-nonzero}--\eqref{eq:1p-evaluated-integrals-k0},
\ref{app:2p-integrals}, and \ref{app:1p-integrals}
for $D=2,1$ are not part of the fast method to be introduced in
section~\ref{sec:spectral-ewald}, but they serve to establish
that the Fourier integrals appearing in
\eqref{eq:fourier-potential-Dp} are well-defined (in the
distributional sense), even in the
singular case. Furthermore, they can be used to validate the fast
method, by truncating the periodic sums
\eqref{eq:2p-evaluated-integrals-nonzero}
and
\eqref{eq:1p-evaluated-integrals-nonzero}
at some maximum absolute wavenumber $k_{i,\text{max}}$ in each periodic
direction (so that only $k_i \in \{ 2\pi \isc{k}_i/L_i :
-\isc{k}_{i,\text{max}} \leq \isc{k}_i \leq \isc{k}_{i,\text{max}}-1,
\isc{k}_i \in \mathbb{Z}\}$ are included in the sum),
and then directly sum them.

For $D=0$, the integral in \eqref{eq:explicit-fourier-0p} can be
seen to exist for all three kernels (also in the Lebesgue sense),
for example by going to spherical coordinates. We do however not compute
it analytically here, since it is typically not needed;
validation is in the $D=0$ case most easily done by directly
summing the original sum \eqref{eq:periodic-potential}.

To summarize, the full Ewald decomposition in arbitrary
periodicity is given by
\begin{gather}
  \label{eq:Dp-ewald-decomp}
  \bvec{u}^{D\per}(\bvec{x}) = \bvec{u}^{D\per,\rp}(\bvec{x};\xi)
  + \bvec{u}^{D\per,\fp}(\bvec{x};\xi),
  \quad \bvec{x} \neq \bvec{x}_m,
  \\[5pt]
  \label{eq:Dp-ewald-decomp-self}
  \bvec{u}^{D\per\star}(\bvec{x}_m) = \bvec{u}^{D\per,\rp\star}(\bvec{x}_m;\xi)
  + \bvec{u}^{D\per,\fp}(\bvec{x}_m;\xi)
  + \bvec{u}^{\selfterm}(\bvec{x}_m;\xi),
  \quad m=1,\ldots,N,
\end{gather}
where \eqref{eq:Dp-ewald-decomp-self} is used if the evaluation
point coincides with one of the source locations, and
\eqref{eq:Dp-ewald-decomp} is used otherwise. Here, the
real-space part $\bvec{u}^{D\per,\rp(\star)}$ is given by
\eqref{eq:real-space-potential-Dp}, and the star signifies that
the term $\bvec{p}=\bvec{0}$ is skipped when $n=m$.
The Fourier-space part $\bvec{u}^{D\per,\fp}$
is given by \eqref{eq:fourier-potential-Dp} and can for $D=3,2,1$
be further decomposed as
\begin{equation}
  \label{eq:Dp-ewald-Fourier-decomp}
  \bvec{u}^{D\per,\fp}(\bvec{x};\xi)
  =
  \bvec{u}^{D\per,\fp,\bvec{k}^\per\neq\bvec{0}}(\bvec{x};\xi)
  +
  \bvec{u}^{D\per,\fp,\bvec{k}^\per=\bvec{0}}(\bvec{x};\xi),
\end{equation}
as in \eqref{eq:2p-evaluated-integrals-nonzero}--\eqref{eq:1p-evaluated-integrals-k0};
for $D=3$, cf.\ \eqref{eq:3p-ewald-decomp} and note that
$\bvec{k}^\per = \bvec{k}$ for $D=3$.
Finally, the term $\bvec{u}^{\selfterm}$ in \eqref{eq:Dp-ewald-decomp-self}
is independent of periodicity; it is zero for the rotlet and
stresslet, and given by \eqref{eq:3p-stokeslet-self-term}
for the stokeslet. Again, the expressions for $\kernel^\rp$ and
$\wh{\kernel}{}^\fp$ for the different kernels are found in
\eqref{eq:stokeslet-real}--\eqref{eq:stresslet-real} and
\eqref{eq:stokeslet-fourier}--\eqref{eq:stresslet-fourier},
respectively.

\subsection{A note on zero modes, far-field behaviour and the
stresslet integral identity}
\label{sec:ewald-summation-k0-note}

Let us note that the term $\bvec{u}^{D\per,\fp,\bvec{k}^\per=\bvec{0}}$
in \eqref{eq:Dp-ewald-Fourier-decomp} is the zero mode of the
Fourier series in the periodic directions; it is a function of
the coordinates in the free directions only. In the triply periodic
case ($D=3$), cf.\ section~\ref{sec:ewald-summation-3p}, the
requirement of zero mean flow \eqref{eq:zero-mean-flow}, was used
to fix the arbitrary constant in the zero mode
\eqref{eq:fourier-space-potential-3p-k0}. However, zero mean flow
cannot be imposed in the $D=2,1$ cases, since the computational
domain is unbounded in the free directions, and, as shown (for the
stokeslet) by \cite{Pozrikidis1996}, in order for the flow field
\eqref{eq:Dp-ewald-decomp} to be smooth it must diverge as
infinity is approached in the free directions. This far-field
behaviour in the free directions is controlled by the zero mode.

With the zero modes derived in \ref{app:2p-integrals} and
\ref{app:1p-integrals}, the far-field behaviour of
$\bvec{u}^{D\per,\fp,\bvec{k}^\per=\bvec{0}}$
in the free directions is for the stokeslet $O(\lvert z \rvert)$
in $D=2$, and $O(\log(\lvert \bvec{r} \rvert))$ in $D=1$. For
both the rotlet and stresslet, the behaviour is $O(\sgn(z))$ in $D=2$, where
$\sgn(\cdot)$ is the sign function, and $O(1/\lvert \bvec{r}
\rvert)$ in $D=1$. In $D=0$, the flow field is uniquely
determined by \eqref{eq:periodic-potential} and goes to zero at
infinity.

In boundary integral methods, one would like the stresslet
integral identity \citep[eq.~2.1.12, p.~21]{Pozrikidis1992}
\begin{equation}
  \label{eq:stresslet-identity}
  u_j^{\stresslet,D\per}(\bvec{x})
  = \int_\Gamma \sum_{\bvec{p} \in P_{D\per}}
  \stresslet_{jlm}(\bvec{x} - \bvec{y} + \bvec{p})
  q^{(0)}_l \nu_m(\bvec{y}) \, \mathrm{d}S(\bvec{y})
  =
  \begin{cases}
    0, & \bvec{x} \in \mathcal{D}_\mathrm{encl}, \\
    4\pi q^{(0)}_j, & \bvec{x} \in \Gamma, \\
    8\pi q^{(0)}_j, & \bvec{x} \in \mathcal{D}_\mathrm{ext}, \\
  \end{cases}
\end{equation}
to hold for an arbitrary constant vector $\bvec{q}^{(0)}$.
Here, $\Gamma$ is a sufficiently smooth surface
enclosing the domain $\mathcal{D}_\mathrm{encl}$, and
$\mathcal{D}_\mathrm{ext}$ is the domain outside $\Gamma$;
furthermore, $\bvec{\nu}$ is the outward-pointing unit normal of
$\Gamma$, and the set $P_{D\per}$ is as in \eqref{eq:periodic-images}.
(Note that \eqref{eq:stresslet-identity} upon discretization of
the integral becomes precisely the stresslet flow field
\eqref{eq:periodic-potential} for a constant $\bvec{q}(\bvec{y}) =
\bvec{q}^{(0)}$.)
As noted by \cite{afKlinteberg2014}, the term
$\bvec{u}^{\stresslet,3\per,\fp,\bvec{k}=\bvec{0}}$ as given by
\eqref{eq:k0-stresslet-3p} is needed for the stresslet integral
identity \eqref{eq:stresslet-identity} to hold in the $D=3$ case.
We have verified numerically that \eqref{eq:stresslet-identity}
holds for all values of $D = 3, 2, 1, 0$ when the flow field is
given by \eqref{eq:Dp-ewald-decomp}--\eqref{eq:Dp-ewald-Fourier-decomp}.

\section{The Spectral Ewald method}
\label{sec:spectral-ewald}

The goal is now to compute the Fourier-space part
\eqref{eq:fourier-potential-Dp} of the periodic potential in an
efficient way. While it is possible to truncate and directly sum
\eqref{eq:explicit-fourier-3p},
\eqref{eq:2p-evaluated-integrals-nonzero} and
\eqref{eq:1p-evaluated-integrals-nonzero}, corresponding to
$D=3,2,1$, respectively, doing so would yield a slow method that
scales at best like $O(N^{3/2})$, where $N$ is both the number of
sources and targets, and also with a large constant due to the
evaluation of special functions in the $D=2,1$ cases.
(In the $D=0$ case, \eqref{eq:periodic-potential} can be
summed directly, which would however scale like $O(N^2)$.) To get a fast
method, we instead introduce a uniform grid and compute the
interaction between sources and targets via the grid,
using the fast Fourier transform (FFT). This is the idea behind
the class of Particle--Mesh--Ewald (PME) methods, to which the
Spectral Ewald (SE) method belongs. To treat all $D=3,2,1,0$ within
the same fast framework, we will discretize the integrals that appear
in \eqref{eq:fourier-potential-Dp} for $D=2,1,0$, rather than
computing them analytically as was done in
section~\ref{sec:ewald-summation-Dp}.

The SE method, like other PME methods,
has the following steps: (i) the sources $\bmat{f}(\bvec{x}_n)$
are spread onto a uniform grid using an interpolating window function
$\window(\bvec{r})$, (ii) an FFT is applied on the grid, (iii)
the result is scaled by the kernel $\wh{\kernel}{}^\fp$, (iv) an
inverse FFT (IFFT) is applied, and finally (v) the result is
interpolated from the uniform grid to the desired target points
using the window function $\window(\bvec{r})$. The SE method
differs from other PME methods in that the support of the window
function can be varied independently of the size of the uniform
grid, which allows approximation errors from the window function
to be controlled separately from truncation errors from the grid,
as described in section~\ref{sec:estimates}. The window function
will be described in section~\ref{sec:se-window-functions};
for now it can be thought of as a generic function
$\mathbb{R}^3 \to \mathbb{R}$ with compact support.

Let us derive formulas for the steps of the SE method by
introducing the
window function into \eqref{eq:fourier-potential-Dp}.
Given a window function $\window(\bvec{r})$ with Fourier
transform $\wh{\window}(\bvec{k})$, the identity
$\wh{\window} (\wh{\window})^{-2} \wh{\window} = 1$ can be
inserted into \eqref{eq:fourier-potential-Dp}, which can then be
arranged as
\begin{equation}
  \label{eq:fourier-potential-Dp-v2}
  \bvec{u}^{D\per,\fp}(\bvec{x};\xi)
  =
  \fourier_{D\per}^{-1}\left\{
  \wh{\window}(\bvec{k})
    \frac{\wh{\kernel}{}^\fp(\bvec{k};\xi)}{[\wh{\window}(\bvec{k})]^2} \matop
  \sum_{n=1}^N
  \wh{\window}(\bvec{k})
  \bmat{f}(\bvec{x}_n)
  \E^{-\I \bvec{k} \cdot \bvec{x}_n}
  \right\} (\bvec{x})
  ,
\end{equation}
with the mixed inverse Fourier transform $\fourier_{D\per}^{-1}$
defined as in \eqref{eq:inverse-mixed-ft}. Recall that the kernel
$\kernel$ may be the stokeslet, rotlet or stresslet,
cf.~\eqref{eq:generic-kernel}--\eqref{eq:tensor-products},
and the expressions for $\wh{\kernel}{}^\fp$ are given in
\eqref{eq:stokeslet-fourier}--\eqref{eq:stresslet-fourier}.
Let us define the gridding interpolant
\begin{equation}
  \label{eq:SE-gridding}
  \Phi(\bvec{x}) :=
  \sum_{n=1}^N
  \sum_{\bvec{p} \in P_{D\per}}
  \window(\bvec{x} - \bvec{x}_n + \bvec{p})
  \bmat{f}(\bvec{x}_n),
\end{equation}
and note that the mixed Fourier transform of $\Phi$ is, by the Poisson
summation formula,
\begin{equation}
  \fourier_{D\per} \{ \Phi \} (\bvec{k}) =
  \sum_{n=1}^N
  \wh{\window}(\bvec{k})
  \bmat{f}(\bvec{x}_n)
  \E^{-\I \bvec{k} \cdot \bvec{x}_n},
\end{equation}
which appears in \eqref{eq:fourier-potential-Dp-v2}. Evaluating
\eqref{eq:SE-gridding} on the uniform grid, which is efficient
since $\window$ has compact support, corresponds to step~(i) of
the SE method, as outlined above. Step~(ii) corresponds to
computing the Fourier transform $\fourier_{D\per} \{ \Phi \} (\bvec{k})$.
Let us furthermore define
\begin{equation}
  \label{eq:SE-scaling}
  \widetilde{\Phi}(\bvec{k};\xi) :=
  \frac{\wh{\kernel}{}^\fp(\bvec{k};\xi)}{[\wh{\window}(\bvec{k})]^2} \matop
  \fourier_{D\per} \{ \Phi \} (\bvec{k}),
\end{equation}
which corresponds to step~(iii) of the SE method. Finally, by
\eqref{eq:fourier-potential-Dp-v2} and the convolution theorem,
\begin{equation}
  \label{eq:SE-gathering}
  \bvec{u}^{D\per,\fp}(\bvec{x};\xi)
  =
  \fourier_{D\per}^{-1}\{
    \wh{\window}(\bvec{k}) \widetilde{\Phi}(\bvec{k};\xi)
  \} (\bvec{x})
  =
  \int_{\boxvar_{D\per}} \fourier_{D\per}^{-1}\{\widetilde{\Phi}\}(\bvec{y})
  \sum_{\bvec{p} \in P_{D\per}}
  \window(\bvec{x} - \bvec{y} + \bvec{p}) \, \mathrm{d}\bvec{y},
\end{equation}
with $\boxvar_{D\per}$ as in
\eqref{eq:mixed-fourier-integration-domain}. Step~(iv) of the
method corresponds to computing
$\fourier_{D\per}^{-1}\{\widetilde{\Phi}\}(\bvec{y})$ on the uniform grid.
Discretizing the integral in \eqref{eq:SE-gathering} on the
uniform grid leads to step~(v) of the SE method.

Several integrals appear in the formulation above, namely in the
operator $\fourier_{D\per}$ \eqref{eq:mixed-fourier-transform} of
step~(ii), in the operator $\fourier_{D\per}^{-1}$
\eqref{eq:inverse-mixed-ft} of step~(iv) for $D=2,1,0$, and in
\eqref{eq:SE-gathering} for step~(v). These will all be
discretized using the trapezoidal rule and, whenever the
integration domain is unbounded, truncated. The discrete sums in
\eqref{eq:inverse-mixed-ft} in the periodic directions are also
truncated. The operators $\fourier_{D\per}$ and
$\fourier_{D\per}^{-1}$ can then be approximated by the FFT and
IFFT, respectively.

Special care must be taken when discretizing the integrals in
\eqref{eq:inverse-mixed-ft} for $D=2,1$, since the kernel
$\wh{\kernel}{}^\fp$ appearing in \eqref{eq:SE-scaling} is
singular for the $\bvec{k}^\per=\bvec{0}$ mode, as mentioned
in section~\ref{sec:ewald-summation-Dp}. Also for $D=0$ the
integrand of \eqref{eq:inverse-mixed-ft} is singular. These
singular cases are treated following \cite{Vico2016} by
introducing modified kernels, which are nonsingular in Fourier space and
defined as the Fourier transform of kernels that have been
truncated in real space such that they
correspond exactly to the original kernels within the primary
cell $\boxvar$ containing the sources.
Furthermore, for $D=2,1$ and
modes $\bvec{k}^\per$ which are close to zero, but not exactly
zero, the discretization of the integrals in \eqref{eq:inverse-mixed-ft}
requires upsampling, since the kernel $\wh{\kernel}{}^\fp$ varies
rapidly close to $\bvec{\kappa}=\bvec{0}$. This is handled by
using an adaptive Fourier transform (AFT) introduced by
\cite{Shamshirgar2017}, which uses a local upsampling factor.

Below, we first introduce the modified kernels in
section~\ref{sec:se-modified-kernels}; the discretized method is
then presented in section~\ref{sec:se-discrete}. The window
function is presented in section~\ref{sec:se-window-functions},
and the AFT is described in section~\ref{sec:se-aft}. In the
$D=0$ case, a precomputation scheme is used to accelerate
computations, and this is outlined in section~\ref{sec:se-precomp}.
Finally, the SE method is summarized in
section~\ref{sec:se-summary}. A large part of the method
(discretization, window function, AFT) is virtually independent
of the specific kernel, and is therefore the same as in
\cite{Shamshirgar2021}, which treats the harmonic kernel; we here
give an overview of all parts of the method, but refer to
the aforementioned paper for a more detailed discussion.

\subsection{Modified kernels of Stokes flow}
\label{sec:se-modified-kernels}

Let us now return to the formulation
\eqref{eq:explicit-fourier-3p}--\eqref{eq:explicit-fourier-0p},
before the window function was introduced. We here consider the
cases where the integrands are singular, i.e.\ the
$\bvec{k}^\per=\bvec{0}$ mode for $D=2,1$, and the integral
\eqref{eq:explicit-fourier-0p} for $D=0$. These integrals are
\begin{align}
  \label{eq:explicit-fourier-k0-2p}
  \bvec{u}^{2\per,\fp,\bvec{k}^\per=\bvec{0}}(\bvec{x};\xi)
  &= \frac{1}{L_1 L_2 2\pi} \sum_{n=1}^N
  \int_{\mathbb{R}}
  \wh{\kernel}{}^\fp(0,0,\kappa_3;\xi) \matop \bmat{f}(\bvec{x}_n)
  \E^{\I (0,0,\kappa_3) \cdot (\bvec{x}-\bvec{x}_n)}
  \, \mathrm{d}\kappa_3
  , \\
  \label{eq:explicit-fourier-k0-1p}
  \bvec{u}^{1\per,\fp,\bvec{k}^\per=\bvec{0}}(\bvec{x};\xi)
  &= \frac{1}{L_1 (2\pi)^2} \sum_{n=1}^N
  \int_{\mathbb{R}^2}
  \wh{\kernel}{}^\fp(0,\kappa_2,\kappa_3;\xi) \matop \bmat{f}(\bvec{x}_n)
  \E^{\I (0,\kappa_2,\kappa_3) \cdot (\bvec{x}-\bvec{x}_n)}
  \, \mathrm{d}\kappa_2 \, \mathrm{d}\kappa_3
  , \\
  \label{eq:explicit-fourier-k0-0p}
  \bvec{u}^{0\per,\fp}(\bvec{x};\xi)
  &= \frac{1}{(2\pi)^3} \sum_{n=1}^N
  \int_{\mathbb{R}^3}
  \wh{\kernel}{}^\fp(\kappa_1,\kappa_2,\kappa_3;\xi) \matop \bmat{f}(\bvec{x}_n)
  \E^{\I (\kappa_1,\kappa_2,\kappa_3) \cdot (\bvec{x}-\bvec{x}_n)}
  \, \mathrm{d}\kappa_1 \, \mathrm{d}\kappa_2 \, \mathrm{d}\kappa_3
  .
\end{align}
We reiterate that while the integrands are singular, these
integrals are well-defined in the distributional sense, as we
showed in section~\ref{sec:ewald-summation-Dp}. To treat them
numerically, however, we use the idea by \cite{Vico2016} to
modify the kernel $\wh{\kernel}{}^\fp$ in order to remove the singularity
at $\bvec{\kappa}=\bvec{0}$. Recall that $\wh{\kernel}{}^\fp =
\wh{\kernel} \wh{\gamma}$, where the screening function $\gamma$
is given by \eqref{eq:hasimoto-screening} for the stokeslet and
stresslet, and by \eqref{eq:ewald-screening} for the rotlet. Both
$\gamma$ and $\wh{\gamma}$ decay fast, and can essentially be
considered to have compact support; they are furthermore
nonsingular. The reason that $\wh{\kernel}$ is singular at
$\bvec{\kappa}=\bvec{0}$ is that $\kernel$ decays slowly as
$\lvert \bvec{r} \rvert \to \infty$, and the idea by Vico et al.\
is to truncate $\kernel$ outside some radius $R$, so that its
Fourier transform becomes nonsingular. The radius $R$ is
selected large enough for the truncated kernel to agree with the
original kernel within the primary cell containing all sources and
targets.

Finding the Fourier transform of the truncated kernel is easier
if the kernel is radial. For this reason, we make use of the
relation $\kernel = \Kopt A \ftpair \wh{\kernel} = \wh{\Kopt}
\wh{A}$, where the radial and scalar kernel $A$ is the biharmonic
($A=\biharmonic$) for the stokeslet and stresslet, and the
harmonic ($A=\harmonic$)
for the rotlet; cf.\ \eqref{eq:diffop-rel-1}--\eqref{eq:diffop-rel-3}
and \eqref{eq:diffop-rel-fourier-1}--\eqref{eq:diffop-rel-fourier-3}.
For each periodicity we define a truncated $A$ according to
\begin{align}
  \label{eq:truncated-A-2p}
  A_R^{2\per}(r_3) &:= \fourier^{-1}\{ \wh{A}(0,0,\cdot) \}(r_3)
  \rect\bigg(\frac{\lvert r_3 \rvert}{R}\bigg), \\
  \label{eq:truncated-A-1p}
  A_R^{1\per}(r_2,r_3) &:= \fourier^{-1}\{ \wh{A}(0,\cdot,\cdot) \}(r_2, r_3)
  \rect\bigg(\frac{\sqrt{r_2^2 + r_3^2}}{R}\bigg), \\
  \label{eq:truncated-A-0p}
  A_R^{0\per}(r_1,r_2,r_3) &:= \fourier^{-1}\{ \wh{A}(\cdot,\cdot,\cdot) \}(r_1, r_2, r_3)
  \rect\bigg(\frac{\sqrt{r_1^2 + r_2^2 + r_3^2}}{R}\bigg),
\end{align}
where the inverse Fourier transform is understood to be one-,
two- and three-dimensional in \eqref{eq:truncated-A-2p},
\eqref{eq:truncated-A-1p} and \eqref{eq:truncated-A-0p},
respectively. Here, the rectangle function is defined by
\begin{equation}
  \rect(r) = \begin{cases}
    1, & \lvert r \rvert \leq 1, \\
    0, & \lvert r \rvert > 1.
  \end{cases}
\end{equation}
Taking the Fourier transform (with appropriate dimensionality) of
$A_R^{D\per}$, we get $\wh{A}_R^{D\per}$, and we then set
\begin{align}
  \label{eq:2p-k0-relation}
  \wh{\kernel}{}_R^{2\per}(0, 0, \kappa_3) &:=
  \wh{\Kopt}(0,0,\kappa_3) \wh{A}_R^{2\per}(\kappa_3), \\
  \label{eq:1p-k0-relation}
  \wh{\kernel}{}_R^{1\per}(0, \kappa_2, \kappa_3) &:=
  \wh{\Kopt}(0,\kappa_2,\kappa_3) \wh{A}_R^{1\per}(\kappa_2,\kappa_3), \\
  \label{eq:0p-k0-relation}
  \wh{\kernel}{}_R^{0\per}(\kappa_1, \kappa_2, \kappa_3) &:=
  \wh{\Kopt}(\kappa_1,\kappa_2,\kappa_3)
  \wh{A}_R^{0\per}(\kappa_1,\kappa_2,\kappa_3).
\end{align}
Let us for clarity summarize the situation in the $D=0$ case.
Replacing $\wh{\kernel}{}^\fp = \wh{\kernel}{} \wh{\gamma}$ with
$\wh{\kernel}{}_R^{0\per} \wh{\gamma}$ in
\eqref{eq:explicit-fourier-k0-0p}, we get
\begin{align}
  \label{eq:explicit-modified-u0p-1}
  \bvec{u}_R^{0\per,\fp}(\bvec{x};\xi)
  &=
  \frac{1}{(2\pi)^3} \sum_{n=1}^N
  \bigg(
  \int_{\mathbb{R}^3}
  \wh{\Kopt}(\bvec{\kappa})
  \wh{A}_R^{0\per}(\bvec{\kappa})
  \wh{\gamma}(\bvec{\kappa};\xi)
  \E^{\I \bvec{\kappa} \cdot (\bvec{x}-\bvec{x}_n)}
  \, \mathrm{d}\bvec{\kappa}
  \bigg)
  \matop \bmat{f}(\bvec{x}_n)
  \\
  \label{eq:explicit-modified-u0p-2}
  &=
  \Kopt \matop \sum_{n=1}^N
  \bmat{f}(\bvec{x}_n)
  (A_R^{0\per} * \gamma)(\bvec{x}-\bvec{x}_n;\xi),
\end{align}
where we have used that $\wh{\Kopt} \E^{\I \bvec{\kappa} \cdot
\bvec{r}} = \Kopt \E^{\I \bvec{\kappa} \cdot \bvec{r}}$. Assuming
that the support of $\gamma$ is contained within a ball of radius
$a_\gamma$, and that $\lvert \bvec{x} - \bvec{x}_n \rvert \leq
\diam(\boxvar)$, where $\diam(\boxvar)$ is the diameter of the
primary cell $\boxvar$ containing all sources, it can be noted that
$\bvec{u}_R^{0\per,\fp}$ agrees exactly with $\bvec{u}^{0\per,\fp}$
as long as $R \geq \diam(\boxvar) + a_\gamma$. Similar conclusions hold
for $D=2,1$. Precisely how $R$ is selected in the SE method is
described further in section~\ref{sec:se-discrete}.

What remains is to derive expressions for $\wh{A}_R^{D\per}$ for
the harmonic and biharmonic kernels for $D=0,1,2$. Starting with
the case $D=0$, the kernels are three-dimensional and given by
$\harmonic(\bvec{r}) = 1/\lvert\bvec{r}\rvert$ for the harmonic,
and $\biharmonic(\bvec{r}) = \lvert\bvec{r}\rvert$ for the
biharmonic. The Fourier transforms of the truncated kernels were
derived by \cite{Vico2016} as
\begin{align}
  \label{eq:vico-harmonic-0p}
  \wh{\harmonic}_R^{0\per}(\bvec{\kappa})
  &=
  \frac{4\pi}{\kappa^2} \Big( 1 - \cos(R\kappa) \Big),
  \\
  \label{eq:vico-biharmonic-0p}
  \wh{\biharmonic}_R^{0\per}(\bvec{\kappa})
  &=
  -\frac{8\pi}{\kappa^4} \Big(
    1 - \big(1-\tfrac{1}{2}R^2\kappa^2 \big) \cos(R\kappa)
    - R \kappa \sin(R \kappa)
  \Big)
  ,
\end{align}
where $\kappa = \lvert \bvec{\kappa} \rvert = \sqrt{\kappa_1^2 +
\kappa_2^2 + \kappa_3^2}$. Both expressions have finite limits as
$\bvec{\kappa} \to \bvec{0}$, namely
$\wh{\harmonic}_R^{0\per}(\bvec{0}) = 2\pi R^2$ and
$\wh{\biharmonic}_R^{0\per}(\bvec{0}) = \pi R^4$. As noted by
\cite{afKlinteberg2017}, the truncated biharmonic \eqref{eq:vico-biharmonic-0p}
has a slower decay in Fourier space than the exact biharmonic
$\wh{B}(\bvec{\kappa}) = -8\pi/\kappa^4$; for large values of
$\kappa$, we have
\begin{equation}
  \label{eq:0p-vico-Btrunc-decay}
  \frac{\wh{\biharmonic}_R^{0\per}(\bvec{\kappa})}{\wh{\biharmonic}(\bvec{\kappa})}
  \sim R^2\kappa^2.
\end{equation}
This slower decay would affect the convergence of the SE method
for kernels based on the biharmonic (i.e.\ the stokeslet and
stresslet). (The rotlet, which is based on the harmonic, is
unaffected since $\wh{\harmonic}_R^{0\per}(\bvec{\kappa}) /
\wh{\harmonic}(\bvec{\kappa}) \sim 1$ for large $\kappa$.)
One way to understand the slower decay of \eqref{eq:vico-biharmonic-0p}
is that the truncation of $\lvert \bvec{r} \rvert$ at $R$ makes the
kernel discontinuous, which introduces terms proportional to
$1/\kappa^2$ in the Fourier transform. Fortunately, this can be
solved in the following simple way. Let us redefine the
biharmonic kernel as
\begin{equation}
  \label{eq:biharmonic-with-constants}
  \biharmonic(\bvec{r}) = \lvert\bvec{r}\rvert + \aB + \bB \lvert \bvec{r} \rvert^2,
\end{equation}
where $\aB$ and $\bB$ are arbitrary real constants. These extra
terms vanish when the biharmonic operator $\nabla^4$ is applied,
reflecting the gauge freedom of the biharmonic equation; thus,
\eqref{eq:biharmonic-with-constants} represents a family of
fundamental solutions to the biharmonic equation, i.e.\
$-\nabla^4 \biharmonic(\bvec{r}) = 8\pi \delta(\bvec{r})$.
Repeating the derivation of Vico et al.\ with this more general
biharmonic kernel, we get
\begin{multline}
  \label{eq:general-biharmonic-0p}
  \wh{\biharmonic}_R^{0\per}(\bvec{\kappa})
  =
  -\frac{8\pi}{\kappa^4} \Big[
    1
    -\Big(1-\tfrac{1}{2}(\aB+R+\bB R^2)R\kappa^2 + 3\bB R \Big) \cos(R\kappa)
    \\
    - \Big( \tfrac{1}{2}(\aB+2R+3\bB R^2) \kappa^2 - 3\bB \Big) \frac{\sin(R \kappa)}{\kappa}
  \Big]
  .
\end{multline}
We have full freedom in choosing $\aB$ and $\bB$, and may select
them to get optimal decay as $\lvert \bvec{\kappa} \rvert \to \infty$.
From \eqref{eq:general-biharmonic-0p}, we see that this
corresponds to $\aB + R + \bB R^2 = 0$ and $\aB + 2R + 3\bB R^2 =
0$, which has the solution $\aB = -\tfrac{1}{2} R$ and $\bB =
-\tfrac{1}{2}R^{-1}$.
(Note that these values of $\aB$ and $\bB$ are precisely
the ones that make \eqref{eq:biharmonic-with-constants}
continuously differentiable everywhere outside the origin when
truncated at $\lvert \bvec{r} \rvert = R$.)
With this selection, the modified biharmonic becomes
\begin{equation}
  \label{eq:optimal-biharmonic-0p}
  \wh{\biharmonic}_R^{0\per}(\bvec{\kappa})
  =
  -\frac{8\pi}{\kappa^4} \left(
    1 + \frac{1}{2} \cos(R\kappa)
    - \frac{3}{2} \frac{\sin(R \kappa)}{R \kappa}
  \right)
  .
\end{equation}
The finite limit as $\bvec{\kappa} \to \bvec{0}$ is
\begin{equation}
  \label{eq:optimal-biharmonic-0p-limit}
  \wh{\biharmonic}_R^{0\per}(\bvec{0}) =
  -\frac{1}{15} \pi R^4.
\end{equation}
With \eqref{eq:optimal-biharmonic-0p}, we have
$\wh{\biharmonic}_R^{0\per}(\bvec{\kappa})/\wh{\biharmonic}(\bvec{\kappa})
\sim 1$ for large $\kappa$.

We now get the modified stokeslet and stresslet by applying \eqref{eq:0p-k0-relation}
to $\wh{A}_R^{0\per}=\wh{\biharmonic}_R^{0\per}$ given by
\eqref{eq:optimal-biharmonic-0p}--\eqref{eq:optimal-biharmonic-0p-limit},
and the modified rotlet by applying \eqref{eq:0p-k0-relation}
to $\wh{A}_R^{0\per}=\wh{\harmonic}_R^{0\per}$ given by
\eqref{eq:vico-harmonic-0p}. The arbitrary constants $\aB$ and
$\bB$ that were added to the biharmonic have no effect on the
stresslet flow field, but one can show, by applying
$\Kop^{\stokeslet}_{jl} = \delta_{jl} \nabla^2 - \nabla_j \nabla_l$
to the extra terms of \eqref{eq:biharmonic-with-constants}, that
the stokeslet $\stokeslet_{jl}$ gains an extra term $4 \bB
\delta_{jl}$, and the stokeslet flow field \eqref{eq:explicit-fourier-k0-0p}
gains an extra contribution
\begin{equation}
  \label{eq:0p-stokeslet-extra-contribution}
  \bvec{u}^{\stokeslet,0\per,\fp,\mathrm{extra}}
  =
  4 \bB \sum_{n=1}^N \bvec{f}(\bvec{x}_n).
\end{equation}
This is just a constant, so it can easily be adjusted for
afterwards, without affecting the time complexity of the
algorithm. In numerical experiments, we will subtract
\eqref{eq:0p-stokeslet-extra-contribution} afterwards such that
the final stokeslet flow field goes to zero at infinity, as
expected from the definition of the stokeslet \eqref{eq:stokes-kernels}.

We move on to $D=1$, where the harmonic and biharmonic kernels
are two-dimensional and given by \citep{Tornberg2016}
\citep[eq.~2.6.16, p.~60]{Pozrikidis1992} \citep[eq.~39]{Paalsson2020}
\begin{align}
  \label{eq:2d-harmonic-kernel}
  \harmonic^{1\per}(r_2,r_3) &= -2 \log(\rho/\ell_\harmonic), \\
  \label{eq:2d-biharmonic-kernel}
  \biharmonic^{1\per}(r_2,r_3) &= - \rho^2 \log(\rho/\ell_\biharmonic) + \cB,
\end{align}
with $\rho = \sqrt{r_2^2 + r_3^2}$. Here, we have included some
gauge constants, similar to the $D=0$ case: $\ell_\harmonic$ and
$\ell_\biharmonic$ are arbitrary positive constants, while $\cB$
is an arbitrary real constant.
The Fourier transforms of the truncated versions of
\eqref{eq:2d-harmonic-kernel}--\eqref{eq:2d-biharmonic-kernel}
were derived by \cite{Vico2016} for $\ell_\harmonic=1$ and
$\ell_\biharmonic=\mathrm{e}$, $\cB=0$. Repeating the derivations
with arbitrary constants, the Fourier transforms become
\begin{align}
  \label{eq:general-1p-harmonic}
  \wh{\harmonic}_R^{1\per}(\bvec{\kappa})
  &=
  \frac{4\pi}{\kappa^2} \Big(
  1 - J_0(R\kappa)
  - R\kappa \log(R/\ell_\harmonic) J_1(R \kappa)
  \Big)
  ,
  \\
  \notag
  \wh{\biharmonic}_R^{1\per}(\bvec{\kappa})
  &=
  -\frac{8\pi}{\kappa^4} \Big(
    1 - J_0(R\kappa)
    - R\kappa(1+\log(R/\ell_\biharmonic)) J_1(R\kappa) \\
  \label{eq:general-1p-biharmonic}
    &\qquad\qquad
    + \tfrac{1}{4} R^2\kappa^2 (1+2\log(R/\ell_\biharmonic)) J_0(R\kappa)
    - \tfrac{1}{4} R\kappa^3 (\cB-R^2\log(R/\ell_\biharmonic)) J_1(R\kappa)
  \Big)
  ,
\end{align}
where $\kappa = \lvert \bvec{\kappa} \rvert = \sqrt{\kappa_2^2 +
\kappa_3^2}$, and $J_\nu(\cdot)$ is the Bessel function of the
first kind and order~$\nu$. We may now select $\ell_\harmonic$,
$\ell_\biharmonic$ and $\cB$ to optimize the decay of
\eqref{eq:general-1p-harmonic}--\eqref{eq:general-1p-biharmonic}
as $\lvert \bvec{\kappa} \rvert \to \infty$. Noting that
\begin{equation}
  \label{eq:bessel-approximation}
  J_\nu(t) = \sqrt{\frac{2}{\pi t}} \Big(
  \cos(t - \tfrac{\pi}{2}\nu - \tfrac{\pi}{4})
  + O(t^{-1})
  \Big),
  \qquad t>0
\end{equation}
holds as $t \to \infty$ \citep[p.~364, 9.2.1]{Abramowitz1972}, we
see that $J_\nu(R \kappa) \sim (R \kappa)^{-1/2}$ for large
$\kappa$. To get optimal decay, we should have
$\log(R/\ell_\harmonic) = 0$, $1+2\log(R/\ell_\biharmonic)=0$,
and $\cB-R^2\log(R/\ell_\biharmonic)=0$, which leads to the
choices $\ell_\harmonic=R$, $\ell_\biharmonic=R \sqrt{\E}$ and
$\cB=-\tfrac{1}{2}R^2$. With these choices, the modified kernels
become
\begin{align}
  \label{eq:optimal-1p-harmonic}
  \wh{\harmonic}_R^{1\per}(\bvec{\kappa})
  &=
  \frac{4\pi}{\kappa^2} \Big( 1 - J_0(R\kappa) \Big)
  ,
  \\
  \label{eq:optimal-1p-biharmonic}
  \wh{\biharmonic}_R^{1\per}(\bvec{\kappa})
  &=
  -\frac{8\pi}{\kappa^4} \Big(
    1 - J_0(R\kappa)
    - \tfrac{1}{2} R\kappa J_1(R\kappa)
  \Big)
  ,
\end{align}
with finite limits
\begin{align}
  \label{eq:optimal-1p-harmonic-limit}
  \wh{\harmonic}_R^{1\per}(\bvec{0})
  &= \pi R^2,
  \\
  \label{eq:optimal-1p-biharmonic-limit}
  \wh{\biharmonic}_R^{1\per}(\bvec{0})
  &= -\frac{\pi R^4}{8}.
\end{align}
The asymptotic behaviours are
$\wh{\harmonic}_R^{1\per}(\bvec{\kappa}) / \wh{\harmonic}(\bvec{\kappa}) \sim 1$
and $\wh{\biharmonic}_R^{1\per}(\bvec{\kappa}) /
\wh{\biharmonic}(\bvec{\kappa}) \sim (R \kappa)^{1/2}$ for large
$\kappa$, which turns out to be sufficient.
For the harmonic kernel, it turns out that the value of
$\ell_\harmonic$ is not important, since the error in the SE
method is dominated by the $\bvec{k}^\per \neq \bvec{0}$ modes in
this case (nevertheless, we select $\ell_\harmonic=R$). For the
biharmonic, the $\bvec{k}^\per = \bvec{0}$ mode dominates the error,
and selecting \eqref{eq:optimal-1p-biharmonic} leads to
noticeable better convergence compared to e.g.\ the choices
$\ell_\biharmonic=\E$, $\cB=0$ from the original derivation by
Vico et al.

Before applying \eqref{eq:1p-k0-relation} to get the modified
stokeslet, stresslet and rotlet for $D=1$, let us note that the
biharmonic kernel can be avoided for the first component (the one
along the periodic direction) of the stokeslet flow field. The
reason is that, cf.\ \eqref{eq:1p-k0-relation},
\begin{equation}
  \wh{\bmat{\stokeslet}}(0,\kappa_2,\kappa_3)
  =
  \wh{\Kopt}^\stokeslet(0,\kappa_2,\kappa_3)
  \wh{\biharmonic}(0,\kappa_2,\kappa_3)
  =
  \begin{bmatrix}
    -\kappa^2 & 0 & 0 \\
    0 & -\kappa_3^2 & \kappa_2 \kappa_3 \\
    0 & \kappa_2 \kappa_3 & -\kappa_2^2 \\
  \end{bmatrix}
  \wh{\biharmonic}(0,\kappa_2,\kappa_3)
  ,
\end{equation}
and since $-\kappa^2 \wh{\biharmonic} = 2 \wh{\harmonic}$, we get
that $\wh{\stokeslet}_{11}(0,\kappa_2,\kappa_3) = 2
\wh{\harmonic}(0,\kappa_2,\kappa_3)$. Thus, the first component
of the stokeslet can be based on the truncated harmonic kernel
$\wh{\harmonic}^{1\per}_R$, which marginally reduces the number
of floating-point operations needed to compute it. A similar
result holds for the stresslet, namely
\begin{equation}
  \wh{\stresslet}_{1lm}(0,\kappa_2,\kappa_3) =
  2 \wh{\harmonic}(0,\kappa_2,\kappa_3)
  \begin{bmatrix}
    0 & \I\kappa_2 & \I\kappa_3 \\
    \I\kappa_2 & 0 & 0 \\
    \I\kappa_3 & 0 & 0 \\
  \end{bmatrix}_{lm}
  ,
\end{equation}
which means that the first component of the stresslet
flow field can also be related to the harmonic kernel. Note that
the relation is a differential relation, as it must be for the
equivalent of \eqref{eq:explicit-modified-u0p-1}--\eqref{eq:explicit-modified-u0p-2}
to hold.

We can now write down the modified kernels in the $D=1$ case.
The modified rotlet, which is completely based on the harmonic
kernel, is given by
\begin{equation}
  \label{eq:mod-rotlet-1p}
  \wh{\bvec{\rotlet}}{}_R^{1\per}(k_1,\kappa_2,\kappa_3) := \begin{cases}
    \wh{\Kopt}^\rotlet(k_1,\kappa_2,\kappa_3) \wh{\harmonic}(k_1,\kappa_2,\kappa_3), &
    k_1 \neq 0, \\
    \wh{\Kopt}^\rotlet(0,\kappa_2,\kappa_3)
    \wh{\harmonic}_R^{1\per}(\kappa_2, \kappa_3), &
    k_1 = 0,
  \end{cases}
\end{equation}
with $\wh{\Kopt}^\rotlet$ as in \eqref{eq:diffop-rel-fourier-2},
$\wh{\harmonic}(\bvec{k}) = 4\pi/\lvert\bvec{k}\rvert^2$,
and $\wh{\harmonic}_R^{1\per}$ as in \eqref{eq:optimal-1p-harmonic}
with the limit \eqref{eq:optimal-1p-harmonic-limit} for
$(\kappa_2,\kappa_3)=(0,0)$.
The modified stokeslet is given by
\begin{equation}
  \label{eq:mod-stokeslet-1p}
  \wh{\bmat{\stokeslet}}{}_R^{1\per}(k_1,\kappa_2,\kappa_3) := \begin{cases}
    \wh{\Kopt}^\stokeslet(k_1,\kappa_2,\kappa_3) \wh{\biharmonic}(k_1,\kappa_2,\kappa_3), &
    k_1 \neq 0, \\[5pt]
    \begin{bmatrix}
      2 & 0 & 0 \\
      0 & 0 & 0 \\
      0 & 0 & 0 \\
    \end{bmatrix}
    \wh{\harmonic}_R^{1\per}(\kappa_2, \kappa_3)
    +
    \begin{bmatrix}
      0 & 0 & 0 \\
      0 & -\kappa_3^2 & \kappa_2 \kappa_3 \\
      0 & \kappa_2 \kappa_3 & -\kappa_2^2 \\
    \end{bmatrix}
    \wh{\biharmonic}_R^{1\per}(\kappa_2, \kappa_3)
    , &
    k_1 = 0,
  \end{cases}
\end{equation}
with $\wh{\Kopt}^\stokeslet$ as in \eqref{eq:diffop-rel-fourier-1},
$\wh{\biharmonic}(\bvec{k})=-8\pi/\lvert\bvec{k}\rvert^4$, and
$\wh{\biharmonic}_R^{1\per}(\kappa_2, \kappa_3)$ as in
\eqref{eq:optimal-1p-biharmonic} with the limit
\eqref{eq:optimal-1p-biharmonic-limit} for
$(\kappa_2,\kappa_3)=(0,0)$.
Finally, the modified stresslet is given by
\begin{equation}
  \label{eq:mod-stresslet-1p}
  \wh{\bmat{\stresslet}}{}_R^{1\per}(k_1,\kappa_2,\kappa_3) := \begin{cases}
    \wh{\Kopt}^\stresslet(k_1,\kappa_2,\kappa_3) \wh{\biharmonic}(k_1,\kappa_2,\kappa_3), &
    k_1 \neq 0, \\
    2\I \bmat{C}^{\stresslet\harmonic}(\kappa_2,\kappa_3)
    \wh{\harmonic}_R^{1\per}(\kappa_2, \kappa_3)
    -
    \I \bmat{C}^{\stresslet\biharmonic}(\kappa_2,\kappa_3)
    \wh{\biharmonic}_R^{1\per}(\kappa_2, \kappa_3)
    , &
    k_1 = 0,
  \end{cases}
\end{equation}
with $\wh{\Kopt}^\stresslet$ as in \eqref{eq:diffop-rel-fourier-3},
and where $\bmat{C}^{\stresslet\harmonic}$ and
$\bmat{C}^{\stresslet\biharmonic}$ are symmetric tensors with
entries given by
\begin{equation}
  C^{\stresslet\harmonic}_{1lm}(\kappa_2,\kappa_3)
  =
  \begin{bmatrix}
    0 & \kappa_2 & \kappa_3 \\
    \kappa_2 & 0 & 0 \\
    \kappa_3 & 0 & 0 \\
  \end{bmatrix}_{lm}
  ,
  \quad
  C^{\stresslet\harmonic}_{2lm}(\kappa_2,\kappa_3)
  =
  \begin{bmatrix}
    \kappa_2 & 0 & 0 \\
    0 & 0 & 0 \\
    0 & 0 & 0 \\
  \end{bmatrix}_{lm}
  ,
  \quad
  C^{\stresslet\harmonic}_{3lm}(\kappa_2,\kappa_3)
  =
  \begin{bmatrix}
    \kappa_3 & 0 & 0 \\
    0 & 0 & 0 \\
    0 & 0 & 0 \\
  \end{bmatrix}_{lm}
  ,
\end{equation}
and
\begin{equation}
  \begin{array}{c}
  C^{\stresslet\biharmonic}_{1lm}(\kappa_2,\kappa_3)
  =
  \begin{bmatrix}
    0 & 0 & 0 \\
    0 & 0 & 0 \\
    0 & 0 & 0 \\
  \end{bmatrix}_{lm}
  ,
  \quad
  C^{\stresslet\biharmonic}_{2lm}(\kappa_2,\kappa_3)
  =
  \begin{bmatrix}
    0 & 0 & 0 \\
    0 & \kappa_2(3\kappa^2 - 2\kappa_2^2) & \kappa_3(\kappa^2 - 2\kappa_2^2) \\
    0 & \kappa_3(\kappa^2 - 2\kappa_2^2) & \kappa_2(\kappa^2 - 2\kappa_3^2) \\
  \end{bmatrix}_{lm}
  ,
  \\[20pt]
  C^{\stresslet\biharmonic}_{3lm}(\kappa_2,\kappa_3)
  =
  \begin{bmatrix}
    0 & 0 & 0 \\
    0 & \kappa_3(\kappa^2 - 2\kappa_2^2) & \kappa_2(\kappa^2 - 2\kappa_3^2) \\
    0 & \kappa_2(\kappa^2 - 2\kappa_3^2) & \kappa_3(3\kappa^2 - 2\kappa_3^2) \\
  \end{bmatrix}_{lm}
  .
  \end{array}
\end{equation}
The values of the gauge constants $\ell_\harmonic$,
$\ell_\biharmonic$ and $\cB$ have no
effect on the rotlet or stresslet flow fields. For the stokeslet,
on the other hand, one can show, similar to how
\eqref{eq:constant-stokeslet-field-contribution-1p} was derived
(see \ref{app:1p-integrals-k0}),
that the extra contribution to the flow field
\eqref{eq:explicit-fourier-k0-1p} is
\begin{equation}
  \label{eq:1p-stokeslet-extra-contribution}
  \bvec{u}^{\stokeslet,1\per,\fp,\bvec{k}^\per=\bvec{0},\mathrm{extra}}
  =
  \left(
    \log(\ell_\harmonic \E) \begin{bmatrix}
      4 & 0 & 0 \\
      0 & 0 & 0 \\
      0 & 0 & 0 \\
    \end{bmatrix}
    +
    \log(\ell_\biharmonic) \begin{bmatrix}
      0 & 0 & 0 \\
      0 & 2 & 0 \\
      0 & 0 & 2 \\
    \end{bmatrix}
  \right)
  \frac{1}{L_1} \sum_{n=1}^N \bvec{f}(\bvec{x}_n).
\end{equation}
As in the $D=0$ case, this is just a constant, so it can easily
be adjusted for afterwards if other values of $\ell_\harmonic$
and $\ell_\biharmonic$ are wanted in the computation of the
actual flow field. In numerical experiments, we use
$\ell_\harmonic = R$ and $\ell_\biharmonic = R\sqrt{\E}$ in the
modified kernels (to get the optimal decay), but (somewhat
arbitrarily) adjust the final flow field such that
\eqref{eq:1p-stokeslet-extra-contribution} vanishes, by setting
$\ell_\biharmonic = 1$ and $\ell_\harmonic = \E^{-1}$. (It can
here be noted that the relation $\ell_\biharmonic /
\ell_\harmonic = \E$ is needed for the relation $\nabla^2
\biharmonic^{1\per} = 2 \harmonic^{1\per}$ to hold.)

For $D=2$, the kernels are one-dimensional, and the relation
\eqref{eq:2p-k0-relation} can in fact be simplified to the point
where the biharmonic kernel can be avoided altogether.
Introducing the kernels
\begin{equation}
  \label{eq:2p-vico-Z-kernel}
  Z^{2\per}(r_3) = 2\pi \sgn(r_3)
  \quad \ftpair \quad
  \wh{Z}^{2\per}(\kappa_3) = -\frac{4\pi\I}{\kappa_3},
\end{equation}
and
\begin{equation}
  \label{eq:2p-vico-H-kernel}
  \harmonic^{2\per}(r_3) = -2\pi \lvert r_3 \rvert
  \quad \ftpair \quad
  \wh{\harmonic}^{2\per}(\kappa_3) = \frac{4\pi}{\kappa_3^2},
\end{equation}
we can write down the relations
\begin{equation}
  \wh{\bmat{\stokeslet}}(0,0,\kappa_3) =
  \wh{\harmonic}^{2\per}(\kappa_3)
  \bmat{D}^\stokeslet,
  \qquad
  \wh{\bvec{\rotlet}}(0,0,\kappa_3) =
  \wh{Z}^{2\per}(\kappa_3)
  \bmat{D}^\rotlet,
  \qquad
  \wh{\bmat{\stresslet}}(0,0,\kappa_3) =
  \wh{Z}^{2\per}(\kappa_3)
  \bmat{D}^\stresslet,
\end{equation}
where the constant tensors $\bmat{D}^\stokeslet$ and
$\bmat{D}^\rotlet$ are given by
\begin{align}
  \bmat{D}^\stokeslet &:= \begin{bmatrix}
    2 & 0 & 0 \\
    0 & 2 & 0 \\
    0 & 0 & 0 \\
  \end{bmatrix}, \\[3pt]
  \bmat{D}^\rotlet &:= \begin{bmatrix}
    0 & 1 & 0 \\
    -1 & 0 & 0 \\
    0 & 0 & 0 \\
  \end{bmatrix},
\end{align}
and the constant symmetric tensor $\bmat{D}^\stresslet$ has
entries given by
\begin{equation}
  D^{\stresslet}_{1lm}
  =
  \begin{bmatrix}
    0 & 0 & -2 \\
    0 & 0 & 0 \\
    -2 & 0 & 0 \\
  \end{bmatrix}_{lm}
  ,
  \quad
  D^{\stresslet}_{2lm}
  =
  \begin{bmatrix}
    0 & 0 & 0 \\
    0 & 0 & -2 \\
    0 & -2 & 0 \\
  \end{bmatrix}_{lm}
  ,
  \quad
  D^{\stresslet}_{3lm}
  =
  \begin{bmatrix}
    -2 & 0 & 0 \\
    0 & -2 & 0 \\
    0 & 0 & -2 \\
  \end{bmatrix}_{lm}
  .
\end{equation}
Truncating the kernels \eqref{eq:2p-vico-Z-kernel} and
\eqref{eq:2p-vico-H-kernel} at $\lvert r_3 \rvert = R$ and taking
the Fourier transform, we get
\begin{align}
  \wh{Z}_R^{2\per}(\kappa_3)
  &=
  -\frac{4\pi \I}{\kappa_3} \Big( 1 - \cos(R\kappa_3) \Big),
  \\
  \wh{\harmonic}_R^{2\per}(\kappa_3)
  &=
  \frac{4\pi}{\kappa_3^2} \Big(
    1 - \cos(R \kappa_3) - R \kappa_3 \sin(R \kappa_3)
  \Big).
\end{align}
The finite limits as $\kappa_3 \to 0$ are
\begin{align}
  \wh{Z}_R^{2\per}(0)
  &=
  0,
  \\
  \wh{\harmonic}_R^{2\per}(0)
  &=
  -2 \pi R^2.
\end{align}
(The asymptotic behaviours are
$\wh{Z}_R^{2\per}(\kappa_3) / \wh{Z}^{2\per}(\kappa_3) \sim 1$
and
$\wh{\harmonic}_R^{2\per}(\kappa_3) / \wh{\harmonic}^{2\per}(\kappa_3) \sim R\kappa_3$
for large $\kappa_3$, which turns out to be sufficient.)
The modified stokeslet, rotlet and stresslet are given by
\begin{align}
  \label{eq:mod-stokeslet-2p}
  \wh{\bmat{\stokeslet}}{}_R^{2\per}(k_1,k_2,\kappa_3) &:= \begin{cases}
    \wh{\Kopt}^\stokeslet(k_1,k_2,\kappa_3) \wh{\biharmonic}(k_1,k_2,\kappa_3), &
    (k_1,k_2) \neq (0,0), \\
    \wh{\harmonic}_R^{2\per}(\kappa_3)
    \bmat{D}^\stokeslet, &
    (k_1,k_2) = (0,0),
  \end{cases}
  \\[2pt]
  \label{eq:mod-rotlet-2p}
  \wh{\bvec{\rotlet}}{}_R^{2\per}(k_1,k_2,\kappa_3) &:= \begin{cases}
    \wh{\Kopt}^\rotlet(k_1,k_2,\kappa_3) \wh{\harmonic}(k_1,k_2,\kappa_3), &
    (k_1,k_2) \neq (0,0), \\
    \wh{Z}_R^{2\per}(\kappa_3)
    \bmat{D}^\rotlet, &
    (k_1,k_2) = (0,0),
  \end{cases}
  \\[2pt]
  \label{eq:mod-stresslet-2p}
  \wh{\bmat{\stresslet}}{}_R^{2\per}(k_1,k_2,\kappa_3) &:= \begin{cases}
    \wh{\Kopt}^\stresslet(k_1,k_2,\kappa_3) \wh{\biharmonic}(k_1,k_2,\kappa_3), &
    (k_1,k_2) \neq (0,0), \\
    \wh{Z}_R^{2\per}(\kappa_3)
    \bmat{D}^\stresslet, &
    (k_1,k_2) = (0,0),
  \end{cases}
\end{align}
respectively. Here, $\wh{\Kopt}^\stokeslet$, $\wh{\Kopt}^\rotlet$
and $\wh{\Kopt}^\stresslet$ are given by
\eqref{eq:diffop-rel-fourier-1}, \eqref{eq:diffop-rel-fourier-2}
and \eqref{eq:diffop-rel-fourier-3}, respectively;
$\wh{\biharmonic}(\bvec{k}) = -8\pi/\lvert \bvec{k} \rvert^4$,
$\wh{\harmonic}(\bvec{k}) = 4\pi/\lvert \bvec{k} \rvert^2$,
and other variables are as above.

For the sake of completeness, we also define modified kernels in
the $D=3$ case, namely
\begin{equation}
  \label{eq:mod-kernel-3p}
  \wh{\kernel}{}_R^{3\per}(k_1,k_2,k_3) = \begin{cases}
    \wh{\Kopt}(k_1,k_2,k_3) \wh{A}(k_1,k_2,k_3), &
    (k_1,k_2,k_3) \neq (0,0,0),
    \\
    \bmatu{0}, &
    (k_1,k_2,k_3) = (0,0,0),
  \end{cases}
\end{equation}
where $A=\biharmonic$ for the stokeslet and stresslet, and
$A=\harmonic$ for the rotlet. The stresslet has a nonzero zero
mode given by \eqref{eq:k0-stresslet-3p}, but that is added
separately since it is in fact a nonperiodic term that depends
on the target point. Since $\bvec{x}$ can be moved out of the sum
in \eqref{eq:k0-stresslet-3p}, the zero mode can clearly be
computed in $O(N)$ operations for $N$ sources and targets.

In summary, the modified kernels $\wh{\kernel}{}_R^{D\per}$ are
given by \eqref{eq:0p-k0-relation} for $D=0$ (with
$\wh{A}_R^{0\per}$ given by \eqref{eq:optimal-biharmonic-0p}
for the stokeslet and stresslet, and by \eqref{eq:vico-harmonic-0p}
for the rotlet), by
\eqref{eq:mod-rotlet-1p}--\eqref{eq:mod-stresslet-1p} for $D=1$,
by \eqref{eq:mod-stokeslet-2p}--\eqref{eq:mod-stresslet-2p} for
$D=2$, and by \eqref{eq:mod-kernel-3p} for $D=3$. These will
replace the kernel $\wh{\kernel}$ in the expression
$\wh{\kernel}{}^\fp = \wh{\kernel} \wh{\gamma}$ in
\eqref{eq:SE-scaling}, i.e.\ step (iii) of the SE method, which
will thus become
\begin{equation}
  \label{eq:SE-scaling-2}
  \widetilde{\Phi}_R(\bvec{k};\xi) :=
  \frac{\wh{\kernel}{}_R^{D\per}(\bvec{k})
  \wh{\gamma}(\bvec{k};\xi)}{[\wh{\window}(\bvec{k})]^2} \matop
  \fourier_{D\per} \{ \Phi \} (\bvec{k}).
\end{equation}

\subsection{Discrete formulation}
\label{sec:se-discrete}

We are now ready to discretize the integrals appearing in
\eqref{eq:SE-gridding}--\eqref{eq:SE-gathering}. As a first step,
we truncate the unbounded domain $\boxvar_{D\per}$
\eqref{eq:mixed-fourier-integration-domain} that appear in the
integrals. The truncated domain
\begin{equation}
  \label{eq:extended-box}
  \extbox_{D\per} := \begin{cases}
    [0,L_1) \times [0,L_2) \times [0,L_3), & \text{if $D=3$}, \\
    [0,L_1) \times [0,L_2) \times
    [-\frac{1}{2}\delta L_3, L_3+\frac{1}{2}\delta L_3), & \text{if $D=2$}, \\
    [0,L_1) \times [-\frac{1}{2}\delta L_2, L_2+\frac{1}{2}\delta L_2)
    \times
    [-\frac{1}{2}\delta L_3, L_3+\frac{1}{2}\delta L_3), & \text{if $D=1$}, \\
    [-\frac{1}{2}\delta L_1, L_1+\frac{1}{2}\delta L_1) \times [-\frac{1}{2}\delta L_2, L_2+\frac{1}{2}\delta L_2)
    \times [-\frac{1}{2}\delta L_3, L_3+\frac{1}{2}\delta L_3), & \text{if $D=0$}.
  \end{cases}
\end{equation}
will also be called the extended box, since it extends the box
$\boxvar = [0,L_1) \times [0,L_2) \times [0,L_3)$,
which contains all sources, by length $\delta L_i$ in each free direction.
The padding of the box $\boxvar$ by $\delta L_i$ in the free
directions is necessary to ensure that the window function is
fully contained in the extended box $\extbox_{D\per}$ (even when
a source point is at the boundary of the box $\boxvar$);
furthermore, $\delta L_i$ must be selected such that
also the screening function, which is introduced through
$\wh{\kernel}{}^\fp$ in \eqref{eq:SE-scaling}, has decayed
sufficiently at the boundary of $\extbox_{D\per}$.
The minimal box pad length $\delta L_i$ is thus determined by the
window function, but depends also on the screening function;
we determine its appropriate value through numerical experiments
as stated in section~\ref{sec:estimates-aft}.

\begin{figure}[htp]
  \centering
  \includegraphics{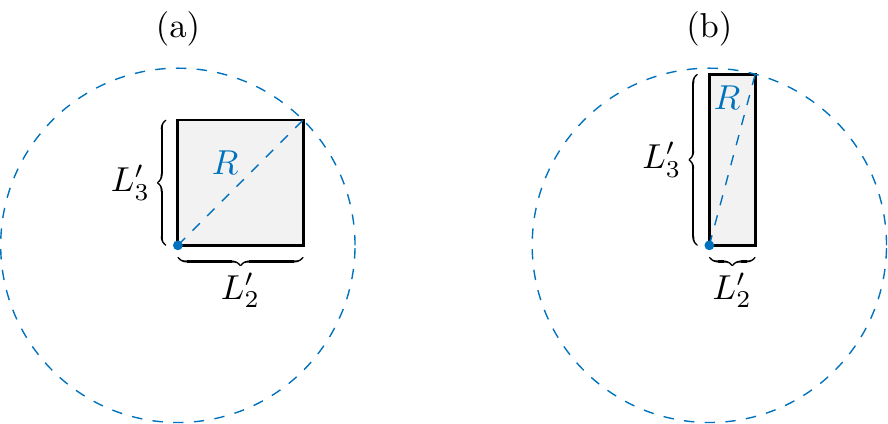}
  \caption{Illustration of how the truncation radius $R$ is set,
  two-dimensional case (corresponding to $D=1$). The rectangle is
  the extended box $\extbox_{D\per}$ projected on the $x_2x_3$-plane,
  and $R$ is its diagonal, given by \eqref{eq:truncation-radius}.
  In (a), the box is a square, while in (b), it has a higher aspect ratio.
  Since the modified kernels are based on radial kernels, $R$
  must be the same in all directions (it must define a circle,
  not an ellipse). This means that the upsampling factor $s_0$,
  given by \eqref{eq:upsampling-factor-s0}, will tend to be
  larger for boxes with high aspect ratio. For this reason,
  the SE method will be most efficient when applied to a box with
  low aspect ratio.}
  \label{fig:truncation-radius}
\end{figure}

The modified kernels defined in
section~\ref{sec:se-modified-kernels} are required to
agree with the original kernels in the extended
box~$\extbox_{D\per}$. The truncation radius $R$ of the modified
kernels should be set as small as possible, since large values of
$R$ will make the kernels more oscillatory, cf.\ e.g.\
\eqref{eq:vico-harmonic-0p} and \eqref{eq:optimal-biharmonic-0p}.
Thus, defining the extended side length
\begin{equation}
  \extL_i = L_i + \delta L_i, \qquad i=1,2,3,
\end{equation}
we set the truncation radius $R$ of the modified kernels to
\begin{equation}
  \label{eq:truncation-radius}
  R = \begin{cases}
    \extL_3, & D=2, \\
    \sqrt{(\extL_2)^2 + (\extL_3)^2}, & D=1, \\
    \sqrt{(\extL_1)^2 + (\extL_2)^2 + (\extL_3)^2}, & D=0,
  \end{cases}
\end{equation}
which are the smallest possible values given that $R$ must be
direction-independent, as illustrated in
Figure~\ref{fig:truncation-radius}.

A uniform Cartesian grid with grid spacing $\gridh$ is introduced
on $\extbox_{D\per}$; the grid has $M_i = L_i/h$ subintervals in
each periodic direction and $\extM_i = \extL_i/h$ subintervals in
each free direction. It is assumed that $L_i$ and $h$ are
selected such that $M_i$ becomes an even integer in each periodic
direction. In each free direction, we adjust $\delta L_i$ to make
sure that $\extM_i$ also becomes an even integer, as described in
section~\ref{sec:estimates-aft}.

The uniform grid also sets the resolution and bandwidth in Fourier
space. In each periodic direction, the resolution is given by
\begin{equation}
  \Delta k_i = \frac{2\pi}{L_i},
\end{equation}
and the $M_i$ discrete wavenumbers are given by
\begin{equation}
  \label{eq:periodic-disc-wavenumbers}
  k_i(\isc{k}_i) = \isc{k}_i \Delta k_i, \qquad
  \isc{k}_i \in \{-\isc{k}_i^\infty, -\isc{k}_i^\infty+1, \ldots,
  \isc{k}_i^\infty-2, \isc{k}_i^\infty-1 \},
\end{equation}
where $\isc{k}_i^\infty := M_i/2$. Note that $-\pi/h
\leq k_i < \pi/h$. In each free direction, the resolution
is given by
\begin{equation}
  \label{eq:free-disc-wavenumber-resolution}
  \Delta \kappa_i = \frac{2\pi}{\extL_i s(\bvec{k}^\per)},
\end{equation}
where $s(\bvec{k}^\per)$ is an adaptive upsampling factor
depending on the periodic wavenumber vector $\bvec{k}^\per$,
described further in section~\ref{sec:se-aft} for $D=1,2$.
For $D=0$, there are no periodic directions, and the upsampling
factor $s = s_0$ in \eqref{eq:free-disc-wavenumber-resolution} is
uniform and fixed to
\begin{equation}
  \label{eq:upsampling-factor-s0}
  s_0 \approx 1 + \frac{R}{\min_{i} \extL_i},
\end{equation}
which is needed to resolve the oscillations of the modified kernel
$\wh{\kernel}{}_R^{0\per}$ \citep{afKlinteberg2017}.
For a cubic box $\boxvar = [0,L)^3$,
\eqref{eq:upsampling-factor-s0} simplifies to $s_0 \approx 1 + \sqrt{3}
\approx 2.8$.
(For $D=3$, there are no free directions and thus no upsampling.)
For $D=0,1,2$, the discrete wavenumbers in each free direction
are given by
\begin{equation}
  \label{eq:free-disc-wavenumbers}
  \kappa_i(\isc{\kappa}_i) = \isc{\kappa}_i \Delta \kappa_i, \qquad
  \isc{\kappa}_i \in \{-\isc{\kappa}_i^\infty, -\isc{\kappa}_i^\infty+1, \ldots,
  \isc{\kappa}_i^\infty-2, \isc{\kappa}_i^\infty-1 \},
\end{equation}
where $\isc{\kappa}_i^\infty := s(\bvec{k}^\per) \extM_i/2$. Note
that again $-\pi/h \leq \kappa_i < \pi/h$.
The upsampling factor $s(\bvec{k}^\per)$ is in
practice adjusted upwards such that $s(\bvec{k}^\per) \extM_i$
becomes an even integer.

The Fourier operators $\fourier_{D\per}$ and $\fourier_{D\per}^{-1}$
are now approximated by the FFT and IFFT on the uniform grid, denoted by
$\fourier_{h,D\per}$ and $\fourier_{h,D\per}^{-1}$, respectively.
For $D=1,2$, these are adaptive FFTs with adaptive upsampling
factor $s(\bvec{k}^\per)$, described further in section~\ref{sec:se-aft}.
Finally, the integral in \eqref{eq:SE-gathering} is approximated
by the trapezoidal rule. The steps of the SE method can now be
written down in discrete from:
\begin{enumerate}
  \renewcommand*{\labelenumi}{(\roman{enumi})}
  \item
    \textit{Gridding:} The gridding interpolant
    \eqref{eq:SE-gridding} is evaluated at the grid points
    $\bvec{x}_j$ of the uniform grid:
    \begin{equation}
      \label{eq:SE-gridding-disc}
      \Phi(\bvec{x}_j) =
      \sum_{n=1}^N
      \sum_{\bvec{p} \in P_{D\per}}
      \window(\bvec{x}_j - \bvec{x}_n + \bvec{p})
      \bmat{f}(\bvec{x}_n).
    \end{equation}
  \item
    \textit{FFT:} An FFT is applied, resulting in
    \begin{equation}
      \label{eq:SE-FFT-disc}
      F_h(\bvec{k}_l) :=
      \fourier_{h,D\per}\{ \Phi(\bvec{x}_j) \}(\bvec{k}_l),
    \end{equation}
    where $\bvec{k}_l$ is the vector of discrete wavenumbers
    given by \eqref{eq:periodic-disc-wavenumbers} and
    \eqref{eq:free-disc-wavenumbers}.
  \item
    \textit{Scaling:} Evaluate the equivalent of
    \eqref{eq:SE-scaling-2}, i.e.\
    \begin{equation}
      \label{eq:SE-scaling-disc}
      \widetilde{\Phi}_{h}(\bvec{k}_l;\xi) :=
      \frac{\wh{\kernel}{}_R^{D\per}(\bvec{k}_l)
      \wh{\gamma}(\bvec{k}_l;\xi)}{[\wh{\window}(\bvec{k}_l)]^2} \matop
      F_h(\bvec{k}_l).
    \end{equation}
  \item
    \textit{IFFT:} Apply an IFFT, resulting in
    \begin{equation}
      \label{eq:SE-IFFT-disc}
      \widetilde{F}_{h}(\bvec{x}_j;\xi) :=
      \fourier_{h,D\per}^{-1} \{ \widetilde{\Phi}_{h}(\bvec{k}_l;\xi) \}(\bvec{x}_j),
    \end{equation}
    where $\bvec{x}_j$ are the grid points of the uniform grid.
  \item
    \textit{Gathering:} Evaluate the trapezoidal rule
    approximation of \eqref{eq:SE-gathering} at the target points
    $\bvec{x}_m$, i.e.\
    \begin{equation}
      \label{eq:SE-gathering-disc}
      \bvec{u}_{h}^{D\per,\fp}(\bvec{x}_m;\xi)
      :=
      h^3 \sum_j \widetilde{F}_h(\bvec{x}_j;\xi)
      \sum_{\bvec{p} \in P_{D\per}}
      \window(\bvec{x}_m - \bvec{x}_j + \bvec{p}),
    \end{equation}
    where the sum $\sum_j$ is over all grid points of the uniform
    grid.
\end{enumerate}

\subsection{The window function}
\label{sec:se-window-functions}

\begin{figure}[h]
  \centering
  \includegraphics{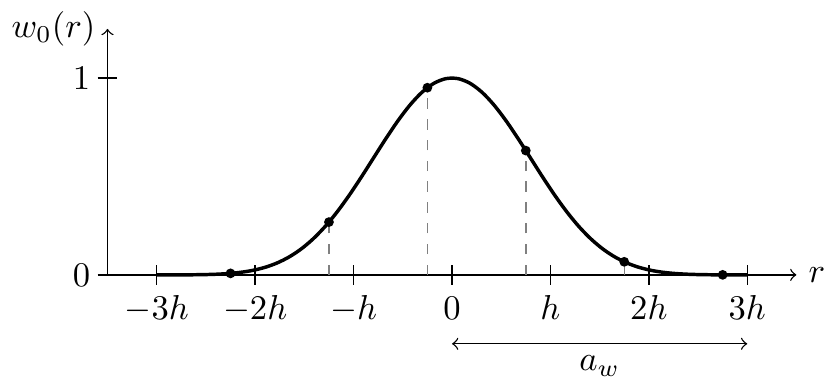}
  \caption{Plot of window function $\window_0$ with support $[-3h, 3h]$
  and evaluation points as in \eqref{eq:1d-window} indicated by
  points. The offset of $-x_{n,i}+p_i$ with respect to the
  uniform grid determines the offset of the evaluation points in
  this figure. The number of evaluation points (window size) is here $P=2
  a_\window/h=6$. The window function shown here is
  \eqref{eq:KB-window} with $\beta=15$.}
  \label{fig:window}
\end{figure}

The time has come to describe the window function
$\window(\bvec{r})$ appearing in the method in more detail. As
already mentioned, the window function should be a function from
$\mathbb{R}^3$ to $\mathbb{R}$ with compact support; we will
furthermore assume that it is given by a tensor product
$\window(\bvec{r}) = \window_0(r_1) \window_0(r_2) \window_0(r_3)$,
where the one-dimensional window $\window_0 : \mathbb{R} \to
\mathbb{R}$ has compact support $[-a_\window, a_\window]$, with
$a_\window > 0$. When evaluating the window function
$\window(\bvec{x}_j - \bvec{x}_n + \bvec{p})$ in
\eqref{eq:SE-gridding-disc}, the point $\bvec{x}_j$ lies on the
uniform grid, and the one-dimensional window~$\window_0$ should
thus be evaluated in
\begin{equation}
  \label{eq:1d-window}
  \window_0(l_i h - x_{n,i} + p_i), \qquad i=1,2,3,
\end{equation}
where $l_i$ is an integer such that $\lvert l_i h - x_{n,i} + p_i
\rvert \leq a_\window$. For an example, see
Figure~\ref{fig:window}, where $a_\window=3h$.
(The situation in \eqref{eq:SE-gathering-disc} is analogous. In
the special case where target points and source points are the
same, it is enough to evaluate the window function for
\eqref{eq:SE-gridding-disc}, since the points in \eqref{eq:SE-gathering-disc}
have the opposite sign, and we assume the window function to be
even, i.e.\ $\window_0(-r) = \window_0(r)$.)
In general, we assume $a_\window$ to be a
multiple of the grid spacing $h$, which means that $\window_0$ is
to be evaluated in $P = 2 a_\window / h$ grid points; for
brevity, we will call the number of evaluation points $P$ the
``window size'' in the following. Clearly, the window size $P$
will be an even integer.

In the original SE method, the window function $\window_0$ was a
truncated Gaussian. Saffar Shamshirgar et al.\ \cite{Shamshirgar2021} compared the truncated
Gaussian window against the Kaiser--Bessel (KB) window function
\citep{Kaiser1980}
\begin{equation}
  \label{eq:KB-window}
  \window_{0,\text{KB}}(r) =
  \begin{cases}
    \displaystyle
    \frac{I_0(\beta\sqrt{1 - (r/a_\window)^2})}{I_0(\beta)}, & \lvert r \rvert \leq a_\window, \\
    0, & \lvert r \rvert > a_\window,
  \end{cases}
  \quad \ftpair \quad
  \wh{\window}_{0,\text{KB}}(k) = \frac{2 a_\window \sinh(\sqrt{\beta^2 - k^2 a_\window^2})}{
    I_0(\beta) \sqrt{\beta^2 - k^2 a_\window^2}}
  ,
\end{equation}
where $I_0(\cdot)$ is the modified Bessel function of the first
kind and order 0, and $\beta > 0$ is a shape parameter. (Note
that the Fourier transform $\wh{\window}_{0,\text{KB}}$ is used in
\eqref{eq:SE-scaling-disc}.) It was found that the window size
$P$ required to achieve a given error tolerance is significantly smaller
(about 40\% smaller) for the KB window \eqref{eq:KB-window} than
the truncated Gaussian window.
While evaluating the KB window $\window_{0,\text{KB}}$ directly is
expensive, it can be approximated by a piecewise polynomial with
an adaptively selected polynomial degree~$\nu$, inspired by the
FINUFFT library \citep{Barnett2019,FINUFFT},
without affecting the overall error of the SE method \citep{Shamshirgar2021}.
The polynomial approximation is constructed by interpolating the
exact KB window in $\nu+1$ Chebyshev points in each of the $P$
subintervals $[l h, (l+1)h]$ shown in Figure~\ref{fig:window}
($l=-P/2, -P/2+1, \ldots, P/2-1$);
the piecewise polynomial is allowed to be discontinuous where
two subintervals meet. The resulting approximation is referred to
as the polynomial Kaiser--Bessel (PKB) window; for more details,
we refer to \cite{Shamshirgar2021}.

Both the shape parameter $\beta$ and the polynomial degree $\nu$
can be tied to the window size $P$ by setting
\begin{equation}
  \label{eq:PKB-parameters}
  \beta = 2.5 P, \qquad
  \nu = \min(\tfrac{1}{2}P+2, 10),
\end{equation}
which means that the PKB window is uniquely determined by $P$.
Using these parameter choices, it was shown by
\cite{Shamshirgar2021} that the PKB window is superior to the
truncated Gaussian window traditionally used in the SE method,
in the sense that the window size $P$ and hence the computational time needed to achieve a
given error tolerance is smaller for the PKB window. The above
mentioned paper treated only the harmonic kernel, but there is
good reason to believe that the conclusion holds for any kernel,
and in section~\ref{sec:results-window} we show that the PKB
window is indeed superior to the truncated Gaussian also for the
stokeslet, stresslet and rotlet kernels. For this reason, we will
focus exclusively on the PKB window throughout
sections~\ref{sec:spectral-ewald} and \ref{sec:estimates}.

\subsection{Adaptive Fourier transform and upsampling}
\label{sec:se-aft}

We now return to the adaptive upsampling factor
$s(\bvec{k}^\per)$ introduced in
\eqref{eq:free-disc-wavenumber-resolution}. For $D=0$, the
upsampling factor is simply given by
\eqref{eq:upsampling-factor-s0}, while for $D=3$ no upsampling is
needed; what remains to explain here are the cases $D=1,2$.
Upsampling is needed to increase the resolution in two cases:
firstly, for the oscillatory modified kernels $\wh{\kernel}{}^{D\per}_R$ used for the
$\bvec{k}^\per=\bvec{0}$ mode (see
section~\ref{sec:se-modified-kernels}), as in the $D=0$ case; secondly, for modes
$\bvec{k}^\per \neq \bvec{0}$ which are close to zero, where the
regular kernels $\wh{\kernel}$ have rapid variations.

\begin{figure}[bp]
  \centering
  \includegraphics{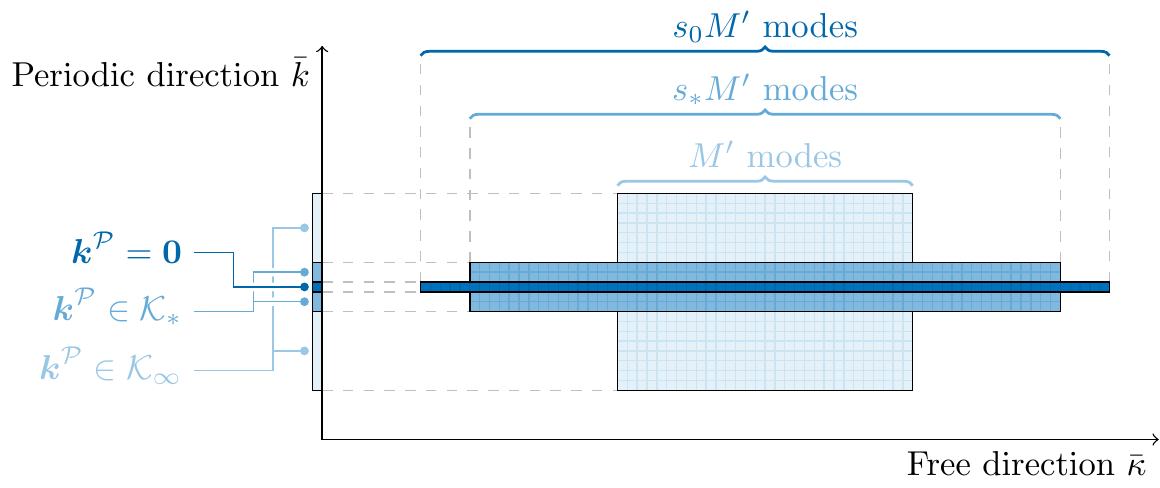}
  \caption{Schematic illustration of adaptive Fourier transform
  (AFT), two-dimensional example. The upsampling factor $s$ in the free
  direction is different for different wavenumbers
  $\bvec{k}^\per$ in the periodic direction, as given by
  \eqref{eq:adaptive-upsampling-factor}. The areas of the shaded
  rectangles represent the number of discrete modes that are
  stored; for example, $\lvert \wavenumset_* \rvert s_*
  \extM$ modes are stored for the set $\wavenumset_*$.}
  \label{fig:aft}
\end{figure}

Thus, in the adaptive Fourier transform (AFT) framework, the
upsampling factor $s(\bvec{k}^\per)$ in the free directions
depends on the wavenumber $\bvec{k}^\per$ in the periodic
directions, as shown in Figure~\ref{fig:aft}.
In practice, this means that modes $\bvec{k}^\per$
that require upsampling in the free directions are stored
separately from modes $\bvec{k}^\per$ that do not. For
simplicity, the set of periodic wavenumbers $\bvec{k}^\per$ is
partitioned into three classes, namely (i)
$\bvec{k}^\per=\bvec{0}$, (ii) the set
\begin{equation}
  \label{eq:wavenumset-upsampling}
  \wavenumset_* := \{ \bvec{k}^\per \in \wavenumset^D
  \setminus \bvec{0} : \lvert k_i \rvert \leq \frac{2\pi}{L_i}
  \isc{k}_i^*, \: i=1,\ldots,D \}
\end{equation}
of nonzero modes to upsample, and (iii) the set
\begin{equation}
  \label{eq:wavenumset-no-upsampling}
  \wavenumset_\infty := \{ \bvec{k}^\per \in \wavenumset^D
  \setminus \wavenumset_* : -\frac{2\pi}{L_i}
  \isc{k}_i^\infty \leq k_i \leq \frac{2\pi}{L_i}
  (\isc{k}_i^\infty - 1), \: i=1,\ldots,D \}
\end{equation}
of modes not to upsample. Here, $\wavenumset^D$ is given in
\eqref{eq:wavenumbers-Dp}, and $\isc{k}_i^\infty = M_i/2$ as in
section~\ref{sec:se-discrete}; the threshold $\isc{k}_i^*$ in
\eqref{eq:wavenumset-upsampling} is a parameter to be selected.
The adaptive upsampling factor is
\begin{equation}
  \label{eq:adaptive-upsampling-factor}
  s(\bvec{k}^\per) = \begin{cases}
    s_0, & \bvec{k}^\per = \bvec{0}, \\
    s_*, & \bvec{k}^\per \in \wavenumset_*, \\
    1, & \bvec{k}^\per \in \wavenumset_\infty,
  \end{cases}
\end{equation}
where $s_0$ is given by \eqref{eq:upsampling-factor-s0} with the
minimization over $i$ restricted to the free directions, and $s_*$
is a parameter to be selected. Note that for a cubic box $\boxvar
= [0,L)^3$, \eqref{eq:upsampling-factor-s0} simplifies to $s_0
\approx 1 + \sqrt{3-D}$, i.e.\ around 2.5 for $D=1$ and 2 for
$D=2$. The values of $s_*$ and $\isc{k}_i^*$ are given in
section~\ref{sec:estimates-aft}.

The AFT is computed by first applying an FFT in the periodic
directions, and then separate FFTs in the free directions for
each of the three cases in \eqref{eq:adaptive-upsampling-factor},
with the appropriate upsampling achieved by zero-padding the free
directions in real space before the second round of FFTs. The
inverse transform (adaptive inverse Fourier transform, AIFT)
first applies three separate IFFTs in the free directions,
truncates and merges the results, and finally applies an IFFT in the
periodic directions. For details, we refer to
\cite{Shamshirgar2021}.

\subsection{Precomputation in the $D=0$ case}
\label{sec:se-precomp}

As mentioned in section~\ref{sec:se-discrete}, the $D=0$ case
requires an upsampling factor $s_0$ given by
\eqref{eq:upsampling-factor-s0} ($s_0 \approx 2.8$ for a cubic box)
in order to resolve the modified kernel
$\wh{\kernel}{}^{0\per}_R$ in the scaling step. This is due to
the oscillatory nature of the modified harmonic \eqref{eq:vico-harmonic-0p}
and biharmonic \eqref{eq:optimal-biharmonic-0p}. However, by
precomputing effective kernels corresponding to
$\wh{\harmonic}^{0\per}_R$ and $\wh{\biharmonic}^{0\per}_R$, the
upsampling factor in the SE method itself can be reduced to 2,
which is the minimum factor needed to compute an aperiodic
convolution by FFTs.

The precomputation is done as follows, where $A$ is used to denote
either the harmonic~$\harmonic$ or biharmonic~$\biharmonic$
kernel: (i) The modified kernel $\wh{A}^{0\per}_R$, i.e.\
\eqref{eq:vico-harmonic-0p} or \eqref{eq:optimal-biharmonic-0p}, is evaluated
on the upsampled uniform grid of size $(s_0 \extM_1) \times (s_0
\extM_2) \times (s_0 \extM_3)$, with $s_0$ as in \eqref{eq:upsampling-factor-s0}.
(ii) A three-dimensional IFFT is applied to get
$A^{0\per}_R$ in real space on the grid. (iii) The result is
truncated to the $(2\extM_1) \times (2\extM_2) \times (2\extM_3)$
points around the origin, resulting in $A^{0\per}_{R,\text{tr}}$.
For the biharmonic only, a ``mollification'' step is applied,
described below.
(iv) A three-dimensional FFT is applied to get
$\wh{A}^{0\per}_{R,\text{tr}}$ in Fourier space on the grid
with upsampling factor 2. For more details, we refer to
\cite{afKlinteberg2017}.

The result $\wh{A}^{0\per}_{R,\text{tr}}$ from the precomputation
is stored, and the tensorial kernels (stokeslet, stresslet,
rotlet) are computed by applying \eqref{eq:0p-k0-relation} to it.
The precomputation is beneficial when computations are to be done
for several different source configurations in the same
computational box $\boxvar$, for example in a time-dependent
simulation, since the precomputation is done only once at the
beginning (or if $\boxvar$ must be resized), and the lower
upsampling factor of 2 can then be used in the SE method for each
source configuration.

For the biharmonic, a special ``mollification'' step is applied
to $A^{0\per}_{R,\text{tr}}$ in step (iii) above. This is needed
since the modified biharmonic kernel \eqref{eq:optimal-biharmonic-0p}
decays as $1/\kappa^4$ in Fourier space, and the truncation in
step (iii) introduces sharp corners which decay as $1/\kappa^2$.
Thus, without mollifying the result, the truncation error would
increase compared to using the grid with upsampling factor $s_0$
throughout the SE method (i.e.\ without precomputation). For the
harmonic \eqref{eq:vico-harmonic-0p}, and indeed for the more
slowly decaying biharmonic \eqref{eq:vico-biharmonic-0p} used by
\cite{afKlinteberg2017}, there is no difference, since those
kernels already decay as $1/\kappa^2$ before truncating. The
mollification for the biharmonic \eqref{eq:optimal-biharmonic-0p}
is done in a tensor product fashion by introducing a
$(4+4)$-point transition band furthest away from the
origin in each spatial direction, and multiplying
$A^{0\per}_{R,\text{tr}}$ by a smooth mollifier in the transition
band. More specifically, we compute
\begin{equation}
  (A^{0\per}_{R,\text{tr},\text{mollified}})_{jlm}
  =
  (A^{0\per}_{R,\text{tr}})_{jlm}
  \mu^{(1)}_j \mu^{(2)}_l \mu^{(3)}_m,
\end{equation}
where indices range over $j = 1, \ldots, 2\extM_1$; $l = 1, \ldots, 2\extM_2$;
$m = 1, \ldots, 2\extM_3$, and $\bvec{\mu}^{(\nu)}$ are vectors
given by
\begin{equation}
  \bvec{\mu}^{(\nu)}
  = (\underbrace{f_\mu(3), f_\mu(2), f_\mu(1), f_\mu(0),}_{\text{4 points}}
  \underbrace{1, 1, \ldots, 1, 1,}_{\text{$2\extM_\nu-8$ points}}
  \underbrace{f_\mu(0), f_\mu(1), f_\mu(2), f_\mu(3)}_{\text{4 points}})
  \in \mathbb{R}^{2 \extM_\nu},
  \qquad \nu=1,2,3.
\end{equation}
This assumes that the origin is at index $\extM_\nu+1$ in each
direction. The mollifier function $f_\mu$ is selected as a sum of
Gaussians
\begin{equation}
  f_\mu(t) = \E^{-a^2 t^2} + \E^{-a^2 (t-7)^2}, \qquad
  \text{with } a = \frac{2}{7} \sqrt{\log(10^2)},
\end{equation}
based on numerical experiments; this function improves the
Fourier-space decay of $A^{0\per}_{R,\text{tr},\text{mollified}}$
enough to not introduce an increase in the truncation error (cf.\
Figure~\ref{fig:estimates-truncation-fourier} in
section~\ref{sec:estimates-truncation}). While it may seem like
selecting the transition band from the original $2 \extM_\nu$
points, rather than extending the grid, would increase the risk
of polluting the result, this does not appear to be a problem in
practice, as shown by the success of the parameter selection
procedure presented in section~\ref{sec:estimates}.

\subsection{Summary of the SE method}
\label{sec:se-summary}

The Spectral Ewald algorithm for computing the Fourier-space part
\eqref{eq:fourier-potential-Dp} of the periodic potential
\eqref{eq:periodic-potential} at arbitrary target locations
$\bvec{x}_m \in \boxvar$ (which may or may not coincide with
source locations) can now be summarized as in
Algorithm~\ref{alg:se}. As usual,
the arbitrary kernel $\kernel$
denotes one of the stokeslet~$\bmat{\stokeslet}$,
rotlet~$\bvec{\rotlet}$ and stresslet~$\bmat{\stresslet}$,
cf.\ \eqref{eq:tensor-products},
while $A$ denotes the harmonic $\harmonic$ (for the rotlet) or biharmonic
$\biharmonic$ (for the stokeslet or stresslet).
For $D=0$, the precomputation step described in
section~\ref{sec:se-precomp} is done separately, in advance.

\begin{breakablealgorithm}
  \caption{Spectral Ewald method (Fourier-space part)}
\small
\label{alg:se}
\begin{algorithmic}[1]\setcounter{ALG@line}{-1}
  \INPUT Source locations $\bvec{x}_n \in \boxvar$ and strengths
  $\bmat{f}(\bvec{x}_n)$ for $n = 1, \ldots, N$; target locations
  $\bvec{x}_m$ for $m=1,\ldots,N_\mathrm{t}$; primary cell~$\boxvar = [0,L_1) \times [0,L_2) \times [0,L_3)$,
  periodicity~$D$, decomposition parameter $\xi$, uniform grid
  spacing $h$, box paddings $\delta L_i$, window size $P$,
  upsampling parameters $s_0$, $s_*$, $\isc{k}_i^*$ ($i=1,2,3$).
  \State Define the extended box $\extbox_{D\per}$ according to
  \eqref{eq:extended-box}.
  If $D \neq 3$, set extended side lengths $\extL_i = L_i +
  \delta L_i$ and truncation radius $R$ according to
  \eqref{eq:truncation-radius}.
  \State \textit{(Gridding)} Introduce a uniform grid
  $\gridvar_{D\per}$ on $\extbox_{D\per}$ with grid spacing $h$
  in each direction. Evaluate $\Phi(\bvec{x}_j)$ for $\bvec{x}_j
  \in \gridvar_{D\per}$ as in \eqref{eq:SE-gridding-disc}.
  \State \textit{(FFT)} If $D=3$, apply a three-dimensional FFT
  to $\Phi(\bvec{x}_j)$ on the grid $\gridvar_{D\per}$, to
  compute $F_h(\bvec{k}_l)$ as in \eqref{eq:SE-FFT-disc}.
  If $D=2,1$, apply an AFT with parameters $s_0$, $s_*$ and
  $\isc{k}_i^*$, as in section~\ref{sec:se-aft}, to compute
  $F_h(\bvec{k}_l)$. If $D=0$, apply a three-dimensional FFT with
  upsampling factor 2 to compute $F_h(\bvec{k}_l)$.
  \State \textit{(Scaling)} Use \eqref{eq:SE-scaling-disc} to
  compute $\widetilde{\Phi}_{h}(\bvec{k}_l;\xi)$. If $D=0$, use
  the precomputed kernel $\wh{\kernel}{}^{0\per}_{R,\text{tr}}
  := \wh{\Kopt} \wh{A}^{0\per}_{R,\text{tr}}$. Otherwise,
  evaluate the modified kernel $\wh{\kernel}{}^{D\per}_{R}$
  as in section~\ref{sec:se-modified-kernels}.
  \State \textit{(IFFT)} Apply a three-dimensional IFFT (if
  $D=3,0$) or AIFT (if $D=2,1$) to compute
  $\widetilde{F}_{h}(\bvec{x}_j;\xi)$ for $\bvec{x}_j \in
  \gridvar_{D\per}$, as in \eqref{eq:SE-IFFT-disc}.
  \State \textit{(Gathering)} Evaluate \eqref{eq:SE-gathering-disc}
  at the target locations $\bvec{x}_m$.
  \OUTPUT Approximation $\bvec{u}_{h}^{D\per,\fp}(\bvec{x}_m;\xi)$
  to the Fourier-space part of the periodic potential,
  $m = 1,\ldots,N_\mathrm{t}$.
\end{algorithmic}
\end{breakablealgorithm}

In particular, for the stokeslet and rotlet kernels, the
source strengths $\bmat{f}(\bvec{x}_n)$ are vectors in
$\mathbb{R}^3$, and the gridding step is done componentwise such
that $\Phi(\bvec{x}_j)$ is also in $\mathbb{R}^3$. Since the
output $\bvec{u}_{h}^{D\per,\fp}(\bvec{x}_m;\xi)$ is also a
vector in $\mathbb{R}^3$, a total of 3 three-dimensional FFTs and
IFFTs are needed (one 3D FFT+IFFT for each vector component).

For the stresslet kernel, the source strengths
$\bmat{f}(\bvec{x}_n)$ are $3 \times 3$-tensors, and gridding is
again done componentwise; in the scaling step, a $3 \times
3$-tensor is transformed into a vector in $\mathbb{R}^3$, such
that the output $\bvec{u}_{h}^{D\per,\fp}(\bvec{x}_m;\xi)$ is a
vector. Thus, the stresslet requires 9 three-dimensional FFTs
and 3 three-dimensional IFFTs.

\section{Error estimates and parameter selection}
\label{sec:estimates}

The Spectral Ewald (SE) method as stated in Algorithm~\ref{alg:se} has
a multitude of parameters, such as the decomposition
parameter~$\xi$, uniform grid spacing~$h$, box paddings~$\delta
L_i$, window size~$P$, and several upsampling parameters. There
is an inherent freedom in the choice of $\xi$, which serves to
shift computational effort between the Fourier-space and
real-space parts. The optimal value of $\xi$, which minimizes the
total computation time, depends on the implementation and
machine, and must be determined by experiment.

This section serves to establish an automated procedure to, given
$\xi$, select all other parameters, such that the error of the
method is as close as possible to a given error tolerance (to be
made more precise below). The
error will be measured in the root mean squared (rms) sense; we
define the absolute rms error as
\begin{equation}
  \label{eq:abs-rms-error}
  \Erms := \sqrt{\frac{1}{N_\mathrm{t}} \sum_{m=1}^{N_\mathrm{t}}
  \lvert \bvec{u}_h(\bvec{x}_m) - \bvec{u}_\mathrm{ref}(\bvec{x}_m)\rvert^2},
\end{equation}
where $N_\mathrm{t}$ is the number of evaluation points
$\bvec{x}_m$, and $\bvec{u}_h$ is the approximation of the
potential given by the SE method, while $\bvec{u}_\mathrm{ref}$
is a reference potential with negligible error. (The nature of
the reference potential will be stated in each case below.)
We also define the relative rms error as
\begin{equation}
  \label{eq:rel-rms-error}
  \Erel := \Erms \Bigg\slash
  \sqrt{\frac{1}{N_\mathrm{t}} \sum_{m=1}^{N_\mathrm{t}}
  \lvert \bvec{u}_\mathrm{ref}(\bvec{x}_m)\rvert^2}
  .
\end{equation}
In \eqref{eq:abs-rms-error} and \eqref{eq:rel-rms-error}, the
potential $\bvec{u}$ may be either the full periodic potential
$\bvec{u}^{D\per}$ given by \eqref{eq:periodic-potential}, the
Fourier-space part $\bvec{u}^{D\per,\fp}$ given by
\eqref{eq:fourier-potential-Dp}, or the real-space part
$\bvec{u}^{D\per,\rp}$ given by \eqref{eq:real-space-potential-Dp};
to clarify what is meant in each case we will write $\Erms^\fp$
for the Fourier-space part error and $\Erms^\rp$ for the
real-space part error (and similarly $\Erel^\fp$ and $\Erel^\rp$
for the relative errors).

The automated parameter selection procedure to be described in
this section will be given an absolute error tolerance $\abstol$,
and is to select parameters such that the actual error is within
one order of magnitude of the given tolerance, i.e.
\begin{equation}
  \frac{\abstol}{10} \leq \Erms \leq 10 \abstol .
\end{equation}
To simplify formulas, we will in sections~\ref{sec:estimates} and
\ref{sec:results} restrict ourselves to the case where the
primary cell $\boxvar$ is a cube of side length $L$, i.e.\
$L_1 = L_2 = L_3 = L$. Note that this is not a restriction of the
method itself, nor of the implementation, which works for any
rectangular cuboid $\boxvar$.

Errors in the SE method come from different sources, such as
truncation errors caused by truncating the Ewald sum at some maximum
wavenumber $k_\infty = \pi/h$ related to the grid spacing~$h$
(cf.\ section~\ref{sec:se-discrete}), and approximation errors
caused by approximating integrals by discrete sums (both in
\eqref{eq:SE-gathering-disc} and in the Fourier transforms in
free directions) and discretizing the window function (cf.\
section~\ref{sec:se-window-functions}). Due to e.g.\ aliasing,
these errors are interdependent, and the total error $\Erms$ is in
general not equal to the sum of the individual error sources
considered in isolation from each other. We will refer to this
interdependence as ``error pollution'', and it needs to be taken
into account when selecting parameters.

We will start by motivating the different aspects of the
parameter selection procedure: estimates for truncation errors,
the potential rms value, and approximation errors, are given in
sections~\ref{sec:estimates-truncation},
\ref{sec:estimates-potential}, and
\ref{sec:estimates-approximation}, respectively; error
pollution and its implications are discussed in
section~\ref{sec:estimates-pollution}; the selection of the box
padding and upsampling parameters in the free directions is described in
section~\ref{sec:estimates-aft}. Finally, the parameter selection
procedure is summarized and demonstrated in
section~\ref{sec:estimates-summary}.

\subsection{Truncation error estimates}
\label{sec:estimates-truncation}

\subsubsection{Fourier-space part truncation error estimates}
\label{sec:estimates-truncation-fourier}

The nonzero grid spacing $h$ used for the uniform grid in the SE
method corresponds to truncation of modes with absolute wavenumbers
greater than $k_\infty := \pi/h$ in the Fourier-space Ewald
sum \eqref{eq:fourier-potential-Dp}. For the harmonic kernel,
there exists an excellent estimate by \cite{Kolafa1992} for the
error caused by omitting the modes above $k_\infty$ (assuming
other error sources are negligible, such as when using direct
summation of the Ewald sum). For the
rotlet, an equally excellent estimate was derived by
\cite{afKlinteberg2016}, namely
\begin{equation}
  \label{eq:fs-trunc-est-rotlet}
  \text{(Rotlet)} \qquad \Etrunc^\fp \approx
  \sqrt{\frac{8 \xi^2 Q}{3\pi L^3 k_\infty}}
  \E^{-k_\infty^2/(2\xi)^2},
\end{equation}
with
\begin{equation}
  \label{eq:vector-source-quantity}
  Q := \sum_{n=1}^{N} \lvert \bvec{f}(\bvec{x}_n) \rvert^2,
\end{equation}
where $\bvec{f}(\bvec{x}_n)$ are the rotlet source strengths.
Here, $\Etrunc^\fp$ is understood to mean the absolute rms error
in the Fourier-space part potential caused by truncation (i.e.\
when other error sources are negligible).
The estimate \eqref{eq:fs-trunc-est-rotlet} works well for any
periodicity~$D = 0,1,2,3$, just like the harmonic estimate does.

\begin{figure}[ht]
  \centering
  \includegraphics[width=0.5\textwidth]{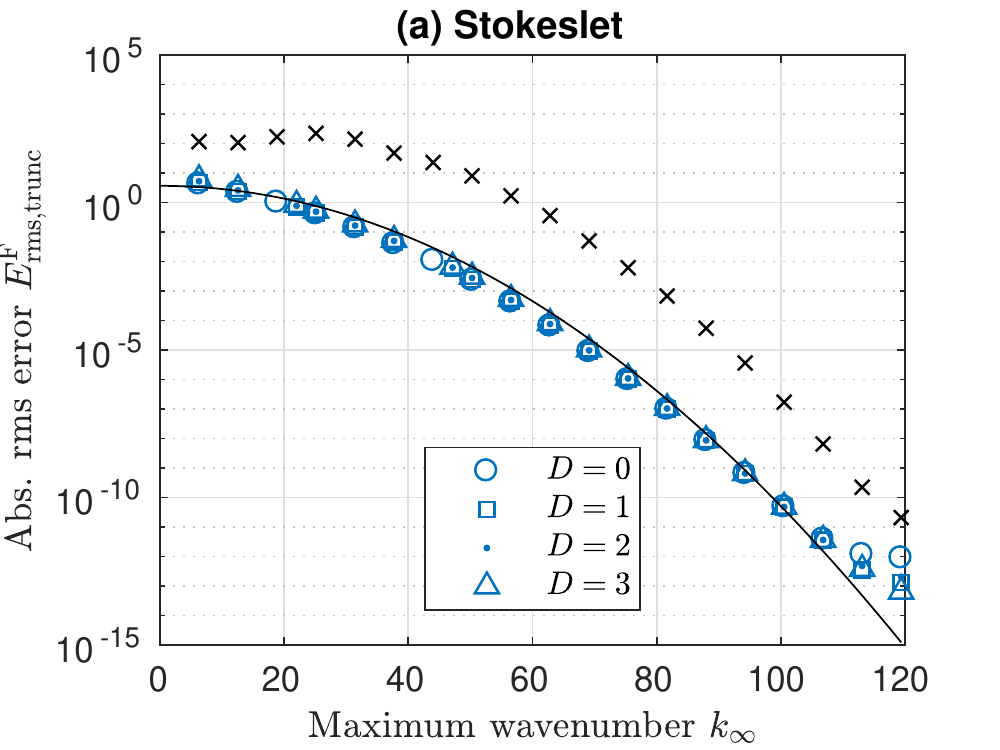}%
  \includegraphics[width=0.5\textwidth]{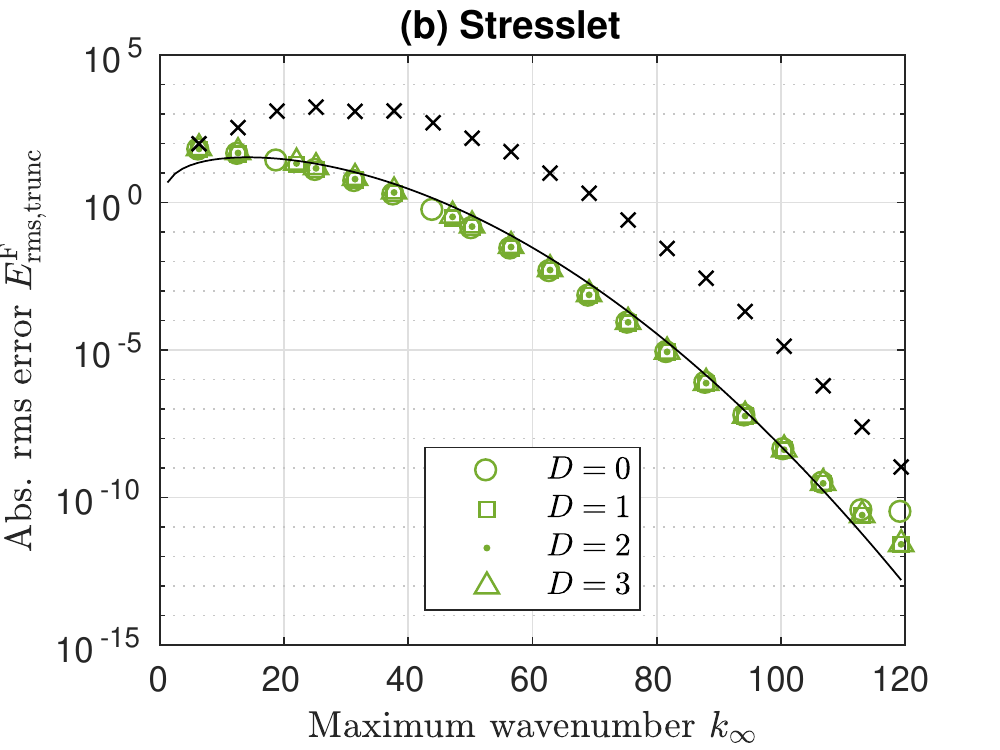}\\
  \includegraphics[width=0.5\textwidth]{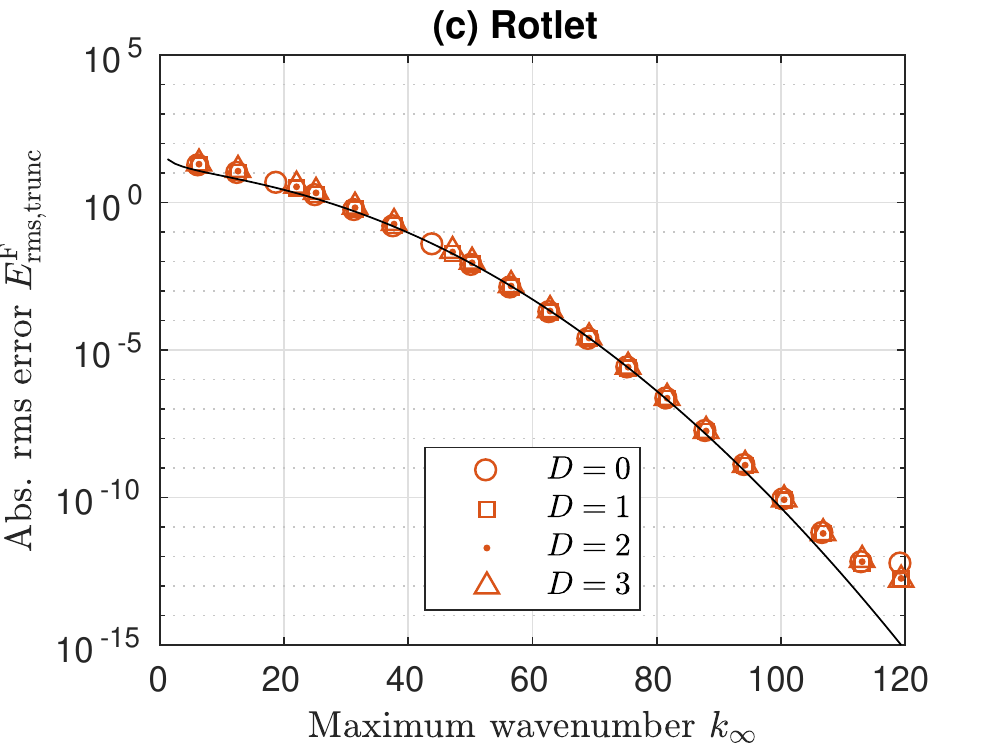}%
  \includegraphics[width=0.5\textwidth]{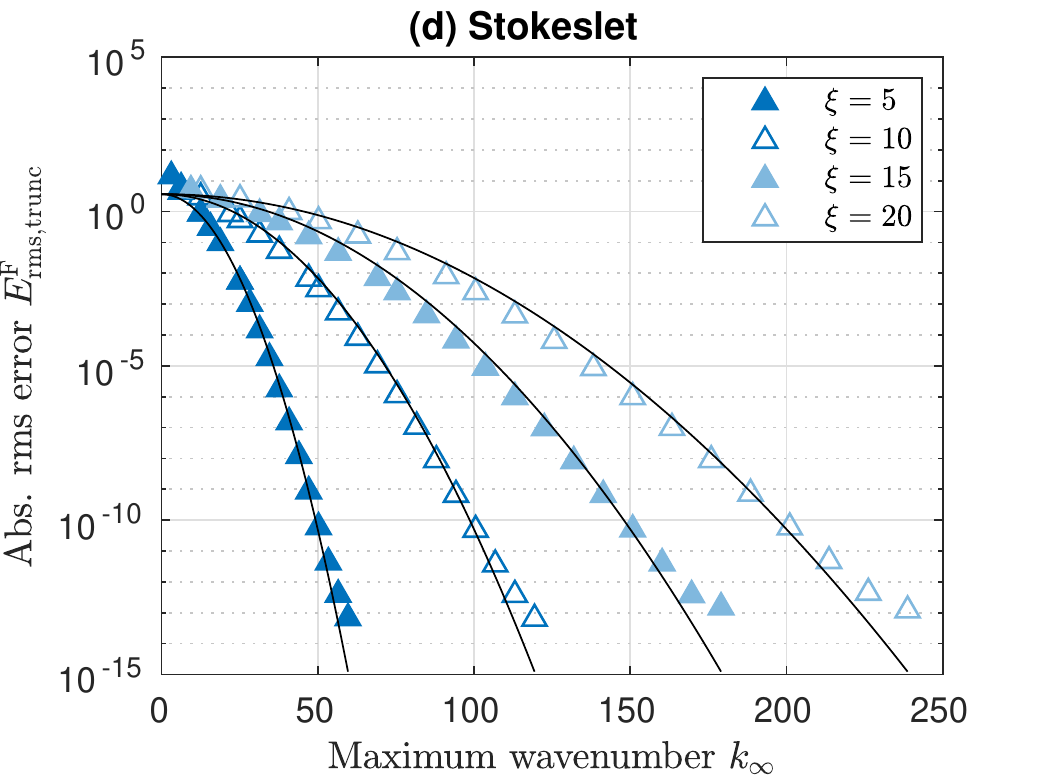}\\
  \caption{Fourier-space part truncation errors and estimates.
  Actual errors are shown as colored symbols,
  estimates \eqref{eq:fs-trunc-est-rotlet}--\eqref{eq:fs-trunc-est-stresslet}
  are shown as solid black curves.
  For comparison, black crosses ($\times$) mark stokeslet and
  stresslet errors for $D=0$ when using the modified biharmonic
  \eqref{eq:general-biharmonic-0p} with $\aB=\bB=0$, which was the
  case in \cite{afKlinteberg2017}.
  Plots (a)--(c) show the three kernels for $\xi = 10$, $D=0,1,2,3$;
  (d) shows the stokeslet for $\xi = 5, 10, 15, 20$ (left to right), $D=3$.
  In all plots, there are $N=10^4$ random (see explanation in
  text) sources which also serve as evaluation points, and $L=2$,
  $Q=100$. Parameters other than $\xi$ and $h=\pi/k_\infty$ are
  selected such that errors other than the truncation error are
  negligible. The reference potential is given by the SE method
  with $k_\infty = 4\pi \xi$.}
  \label{fig:estimates-truncation-fourier}
\end{figure}

For the stokeslet, the situation is a bit more complicated. Since
the modified stokeslet (used for $D=0$ and for the
$\bvec{k}^\per=\bvec{0}$ mode for $D=1,2$) is based on the
truncated biharmonic kernel, special care is needed to ensure
that the decay is fast enough as $\kappa \to \infty$, as
described in section~\ref{sec:se-modified-kernels}. If this is
not done, the truncation error becomes much larger in the $D=0$
case, as illustrated by the black crosses in
Figure~\ref{fig:estimates-truncation-fourier}; a truncation error
estimate specifically for this case was derived by \cite{afKlinteberg2017}.
However, using our
optimized modified kernel from section~\ref{sec:se-modified-kernels},
the truncation error becomes independent of $D$, as it is for the
rotlet. Truncation error estimates for $D=3,2$ have previously
been derived by \cite{Lindbo2010} and \cite{Lindbo2011b}. In
\ref{app:truncation-estimates}, we derive an improved, sharper
estimate using the technique by \cite{afKlinteberg2017} adapted
to $D=3$; since the truncation error is now independent of $D$,
this estimate is valid for any periodicity. Thus,
\begin{align}
  \label{eq:fs-trunc-est-stokeslet}
  \text{(Stokeslet)} \qquad \Etrunc^\fp &\approx
  \frac{4}{\pi L} \sqrt{\frac{Q}{3}}
  \E^{-k_\infty^2/(2\xi)^2},
\end{align}
with $Q$ given by \eqref{eq:vector-source-quantity} where
$\bvec{f}(\bvec{x}_n)$ are the stokeslet source strengths.

The situation for the stresslet, which is also based on the
biharmonic, is very similar to the stokeslet, and in fact the
stresslet estimate can be related to the stokeslet estimate.
Previously, an estimate for $D=3$ has been constructed by
\cite{afKlinteberg2014} using curve fitting, and an estimate for
$D=0$ was derived by \cite{afKlinteberg2017}. We derive an
improved estimate in \ref{app:truncation-estimates}, which is
valid for any periodicity. Our truncation error estimate is
\begin{align}
  \label{eq:fs-trunc-est-stresslet}
  \text{(Stresslet)} \qquad \Etrunc^\fp &\approx
  \frac{4 k_\infty}{3 \pi L} \sqrt{\frac{7Q}{2}}
  \E^{-k_\infty^2/(2\xi)^2},
\end{align}
with
\begin{equation}
  \label{eq:tensor-source-quantity}
  Q := \sum_{n=1}^N \sum_{l,m=1}^3 (f_{lm}(\bvec{x}_n))^2,
\end{equation}
where $\bmat{f}(\bvec{x}_n)$ are the stresslet source strengths.

The estimates are illustrated by an example with random sources in
Figure~\ref{fig:estimates-truncation-fourier}.
Whenever we refer to random sources in this paper, we will,
unless otherwise stated, mean that the locations of the sources
$\bvec{x}_n$ are uniformly distributed within the primary cell
$\boxvar$, and that each component of the source strengths
$\bmat{f}(\bvec{x}_n)$ is uniformly distributed in the interval
$[-a,a]$, where $a$ is adjusted a posteriori to get the value of
$Q$ stated in the example.

Note that the estimates
\eqref{eq:fs-trunc-est-rotlet}--\eqref{eq:fs-trunc-est-stresslet}
can all be solved for $k_\infty = \pi/h$, in most cases using the
Lambert $W$ function, defined as the solution to $W(t) \E^{W(t)} =
t$; thus, $k_\infty$, and hence the grid spacing $h$, can be
computed given $\xi$ and an absolute error tolerance.
Once the desired grid spacing $h=\pi/k_\infty$ has been computed
from the estimates, the grid size $M$, i.e.\ the number of
subintervals of the uniform grid in each periodic direction, can
be computed as $M=L/h$. However, $M$ must at the very least be an
integer, and we will require it to be a multiple of
$\gridfactor$, where $\gridfactor$ is some positive integer.
Thus, the grid spacing is adjusted according to
\begin{equation}
  \label{eq:grid-adjustment}
  h_\text{actual} = \frac{L}{M_\text{actual}}, \qquad
  M_\text{actual} := \gridfactor \left\lceil
  \dfrac{L/h_\text{target}}{\gridfactor}
  \right\rceil,
\end{equation}
where $h_\text{actual}$ is the grid spacing used in the SE
method, $h_\text{target}$ is the one given by the
truncation error estimate, and $\lceil\cdot\rceil$ is the ceiling
function. Typically, $\gridfactor$ is set to a small
power of two to increase the efficiency of the FFTs. Throughout
section~\ref{sec:estimates}, we use $\gridfactor=2$ in the
numerical examples.

\subsubsection{Real-space part truncation error estimates}
\label{sec:estimates-truncation-real}

To efficiently compute the real-space part potential
\eqref{eq:real-space-potential-Dp}, a cut-off radius $\rc > 0$ is
introduced, and only terms for which
\begin{equation}
  \lvert \bvec{x} - \bvec{x}_n + \bvec{p} \rvert < \rc
\end{equation}
are included in the sum. The primary cell $\boxvar$ is
divided into a uniform Cartesian mesh of rectangular cuboid
subcells, such that the side lengths of each subcell are no less than
$\rc$. A list of the points in each subcell is constructed, so that
the computation of \eqref{eq:real-space-potential-Dp} can be done
efficiently, only considering the 27 neighbouring subcells of each
evaluation point.

\begin{figure}[ht]
  \centering
  \includegraphics[width=0.33\textwidth]{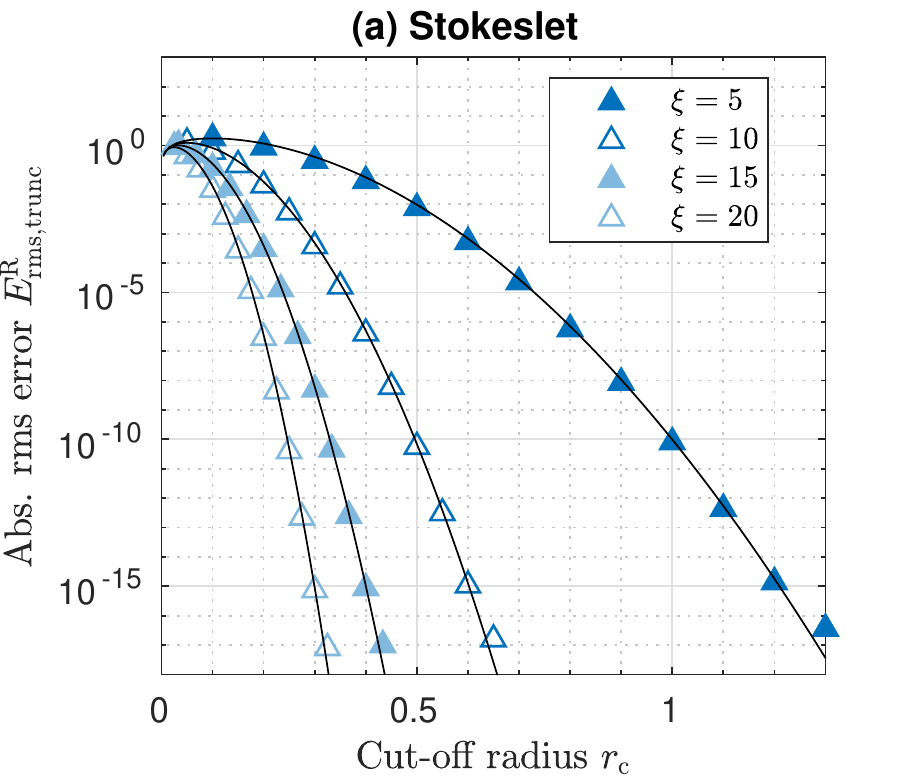}%
  \includegraphics[width=0.33\textwidth]{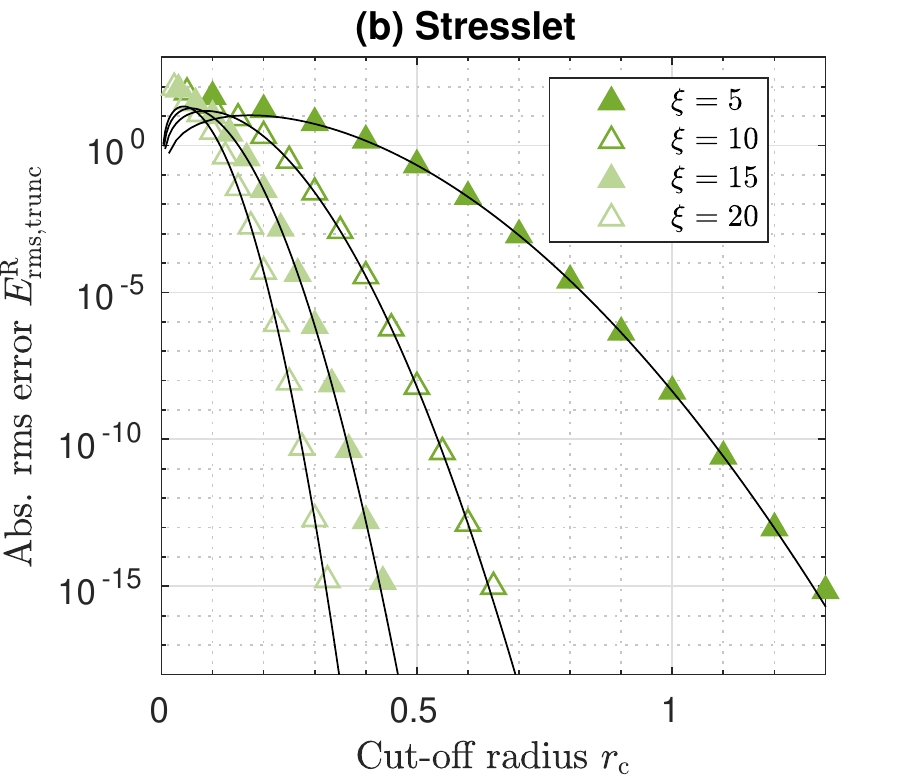}%
  \includegraphics[width=0.33\textwidth]{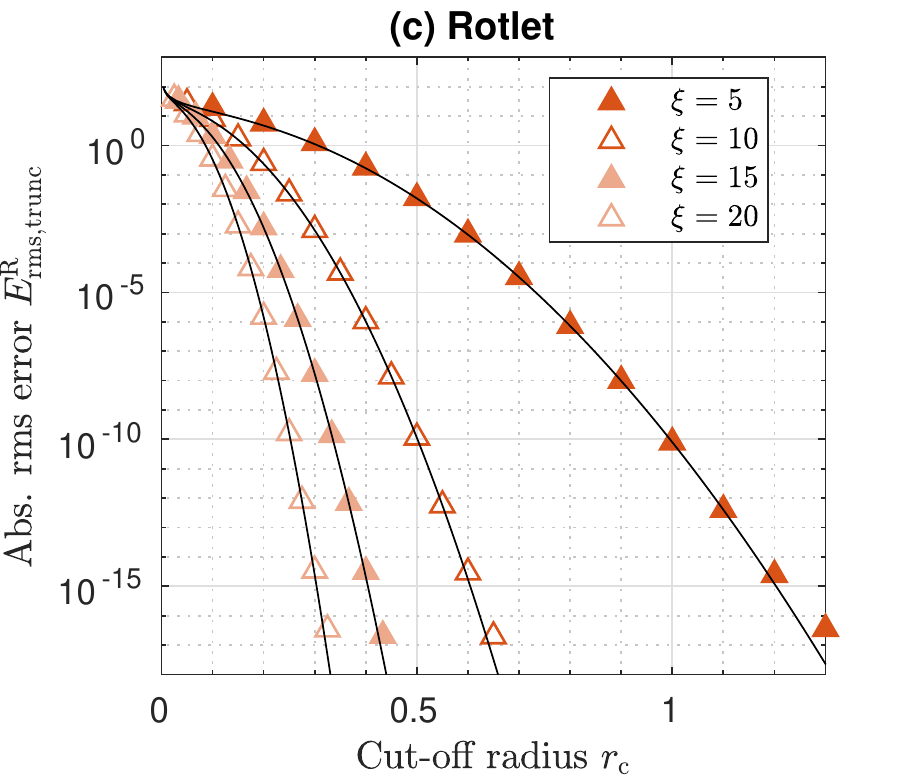}%
  \caption{Real-space part truncation errors and estimates.
  (a)--(c) show the three kernels for $\xi = 5, 10, 15, 20$
  (right to left),
  $D=3$ (errors for other periodicities are almost identical and
  not shown). In
  all plots, there are $N=10^4$ random sources which also serve as
  evaluation points, and $L=2$, $Q=100$. The reference potential
  is given by the same method with $\rc = L$. Actual
  errors are shown as symbols. Estimates
  \eqref{eq:rs-trunc-est-stokeslet}--\eqref{eq:rs-trunc-est-rotlet}
  are shown as solid curves.}
  \label{fig:estimates-truncation-real}
\end{figure}

Truncating the real-space Ewald sum at $\rc$ introduces an error,
which is well described by existing estimates for the stokeslet
\citep{Lindbo2011b}, stresslet \citep{afKlinteberg2017}, and
rotlet \citep{afKlinteberg2016}:
\begin{align}
  \label{eq:rs-trunc-est-stokeslet}
  \text{(Stokeslet)} \qquad \Etrunc^\rp &\approx
  \sqrt{\frac{4Q\rc}{L^3}}
  \E^{-\xi^2 \rc^2},
  \\
  \label{eq:rs-trunc-est-stresslet}
  \text{(Stresslet)} \qquad \Etrunc^\rp &\approx
  \sqrt{\frac{112Q\xi^4\rc^3}{9L^3}}
  \E^{-\xi^2 \rc^2},
  \\
  \label{eq:rs-trunc-est-rotlet}
  \text{(Rotlet)} \qquad \Etrunc^\rp &\approx
  \sqrt{\frac{8Q}{3 L^3 \rc}}
  \E^{-\xi^2 \rc^2},
\end{align}
where $Q$ is given by \eqref{eq:vector-source-quantity} for the
stokeslet and rotlet, and by \eqref{eq:tensor-source-quantity}
for the stresslet. Estimates \eqref{eq:rs-trunc-est-stokeslet}--\eqref{eq:rs-trunc-est-rotlet}
work well for any periodicity $D=0,1,2,3$, for all kernels.
The estimates are illustrated in
Figure~\ref{fig:estimates-truncation-real}. They can all be
solved for $\rc$, again using the Lambert $W$ function, so that
they allow us to determine $\rc$ given $\xi$ and an absolute
error tolerance.

\subsection{Potential rms value estimates}
\label{sec:estimates-potential}

The truncation error estimates in
section~\ref{sec:estimates-truncation} are for the absolute error
\eqref{eq:abs-rms-error}, and this is the reason that the
parameter selection procedure takes an absolute error tolerance,
rather than a relative tolerance. Of course, once the potential
has been computed, one can compute the relative error
\eqref{eq:rel-rms-error} a posteriori (dividing the absolute
error by the rms value of the potential), and relate absolute and
relative errors. If one can find an estimate of the rms value
of the potential, i.e.\ the denominator of \eqref{eq:rel-rms-error},
it would be possible to make an a priori estimate of the relative
error, such that a relative tolerance can be directly given.

For the Fourier-space part of the potential, it turns out that
one can estimate the rms value using only the box side length
$L$, source strength quantity $Q$ and decomposition parameter
$\xi$, as we outline below. Not only can this be used to estimate
the relative error of the Fourier-space part, but it will also
turn out to play a crucial role in the error estimates to follow
in sections~\ref{sec:estimates-approximation} and
\ref{sec:estimates-aft}, and is thus a critical piece of the
parameter selection procedure. In contrast, the full potential or
real-space part cannot have such a simple rms value estimate,
since they depend heavily on how close to the source points the
potential is evaluated, and thus have large variations in space.
While it might be possible to derive more complicated estimates
for the full or real-space part potentials, for example taking
into account the distance to the source points, having such
estimates is not critical for the parameter selection procedure,
and we will therefore not do it here.

Shifting our attention to the Fourier-space part, we want to
find, for each kernel, a quantity $\potmag(L, Q, \xi)$ such that
\begin{equation}
  \potmag(L, Q, \xi) \approx
  \sqrt{\frac{1}{N_\mathrm{t}} \sum_{m=1}^{N_\mathrm{t}}
  \lvert \bvec{u}^{D\per,\fp}(\bvec{x}_m)\rvert^2}
  .
\end{equation}
Using numerical experiments and curve fitting, we have determined
the estimates
\begin{align}
  \label{eq:pot-mag-est-stokeslet}
  \text{(Stokeslet)} \qquad \potmag^\stokeslet(L,Q,\xi) &\approx
  C_\potmag^\stokeslet \sqrt{Q} f_\potmag^\stokeslet(\xi L) / L,
  \\
  \label{eq:pot-mag-est-stresslet}
  \text{(Stresslet)} \qquad \potmag^\stresslet(L,Q,\xi) &\approx
  C_\potmag^\stresslet \sqrt{Q} f_\potmag^\stresslet(\xi L) / L^2,
  \\
  \label{eq:pot-mag-est-rotlet}
  \text{(Rotlet)} \qquad \potmag^\rotlet(L,Q,\xi) &\approx
  C_\potmag^\rotlet \sqrt{Q} f_\potmag^\rotlet(\xi L) / L^2,
\end{align}
with $Q$ as in section~\ref{sec:estimates-truncation}. Here, the
functions $f_\potmag$ for each kernel are given by (again from
curve fitting)
\begin{align}
  \text{(Stokeslet)} \qquad
  f_\potmag^\stokeslet(t) &= \Big(1 + (1.323 \times 10^{-2})t +
  (2.469 \times 10^{-4})t^2 \Big) \E^{-5.205/t^2},
  \\
  \text{(Stresslet)} \qquad
  f_\potmag^\stresslet(t) &= \sqrt{t},
  \\
  \text{(Rotlet)} \qquad
  f_\potmag^\rotlet(t) &= \sqrt{t} \E^{-11.60/t^2}.
\end{align}
The constants $C_\potmag$ in
\eqref{eq:pot-mag-est-stokeslet}--\eqref{eq:pot-mag-est-rotlet}
depend on the specific configuration of source points and
evaluation points. Assuming that the points are drawn from a uniformly
random distribution (as described in the paragraph following
\eqref{eq:tensor-source-quantity}),
the constants $C_\potmag$ have approximately
symmetrical distributions, with mean and standard deviation given
in Table~\ref{tab:estimates-potential}. The data in this table is
generated using 160 different random point systems per kernel,
and using $\xi L = 5, 10, 15, 20, 25, 30$ for each system;
furthermore, $L=1$ and $Q=1$ since $C_\potmag$ does not depend on
these parameters. This was done in the fully periodic ($D=3$)
case, but the order of magnitude of the Fourier-space part
potential does not depend strongly on the periodicity, so
\eqref{eq:pot-mag-est-stokeslet}--\eqref{eq:pot-mag-est-rotlet}
can be used in any periodicity. When using
\eqref{eq:pot-mag-est-stokeslet}--\eqref{eq:pot-mag-est-rotlet}
in the parameter selection procedure, the mean value of
$C_\potmag$ from Table~\ref{tab:estimates-potential} is used.

\begin{table}[h]
  \centering
  \caption{Mean and standard deviation of the constant $C_\potmag$ in
  \eqref{eq:pot-mag-est-stokeslet}--\eqref{eq:pot-mag-est-rotlet},
  assuming uniform point distribution.}
  \label{tab:estimates-potential}
  \begin{tabular}{ccc}
    \toprule
    Kernel & Mean & Standard deviation \\
    \midrule
    Rotlet & 2.4 & 0.2 \\
    Stokeslet & 1.8 & 0.3 \\
    Stresslet & 7.2 & 0.7 \\
    \bottomrule
  \end{tabular}
\end{table}

These potential rms value estimates allow us to approximate the
relative error of the Fourier-space part as $\Erel^\fp \approx \Erms^\fp /
\potmag(L,Q,\xi)$. We do not account for further tests of these
estimates here, but they will be implicitly tested as part of the
parameter selection procedure in
section~\ref{sec:estimates-summary}.

\subsection{Approximation error estimates}
\label{sec:estimates-approximation}

\begin{figure}[ht]
  \centering
  \includegraphics[width=0.5\textwidth]{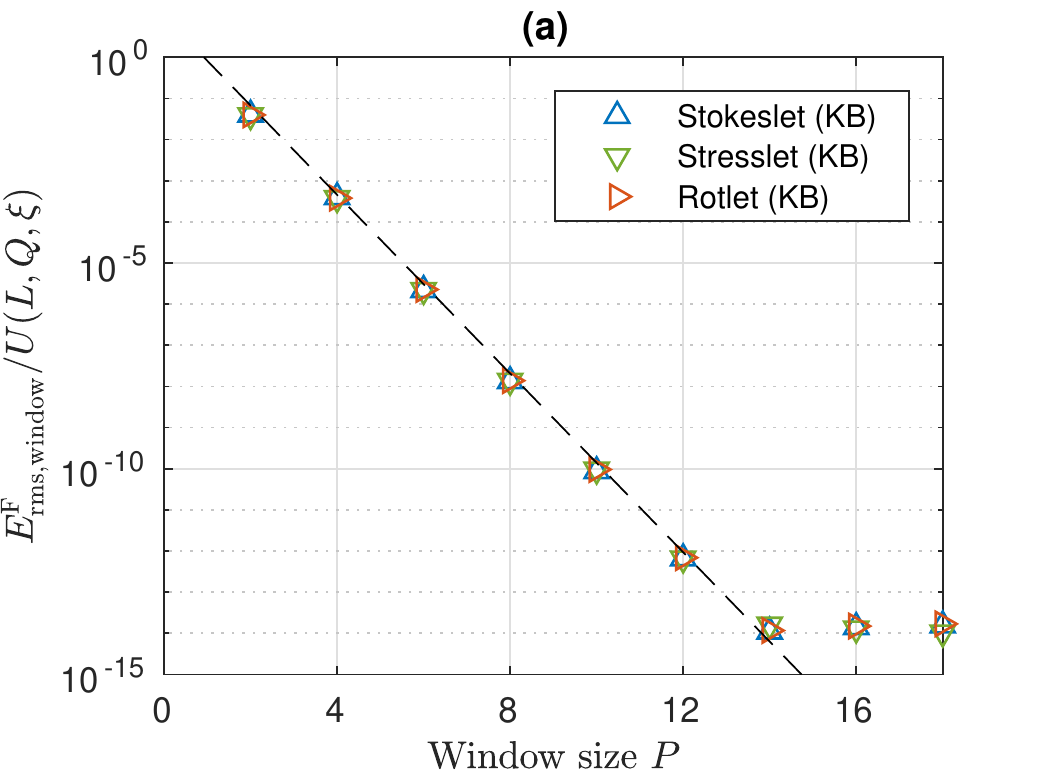}%
  \includegraphics[width=0.5\textwidth]{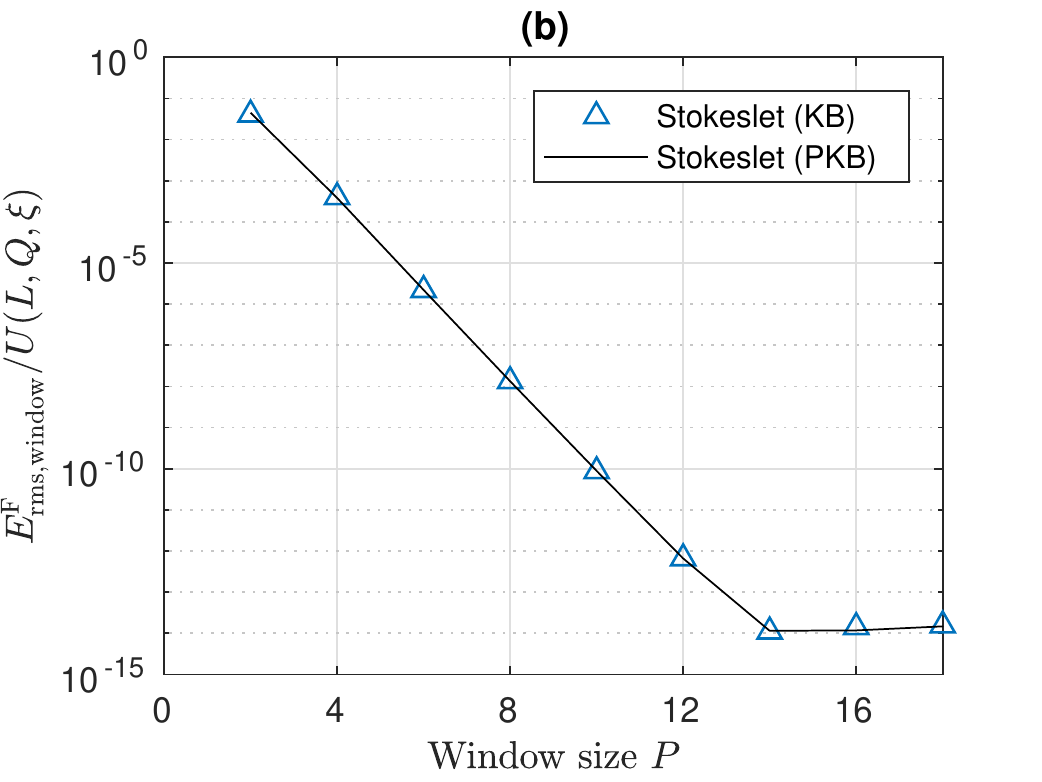}%
  \caption{Window approximation errors. (a) Actual error for the
  three kernels when the exact KB window \eqref{eq:KB-window} is
  used (symbols), and the estimate \eqref{eq:approximation-estimate}, i.e.\
  $\Ewindow^\fp/U(L,Q,\xi) = 10 \E^{-2.5 P}$ (dashed black line);
  the data points are almost on top of each other.
  (b) Actual error for the stokeslet when the exact KB window
  (blue triangles), and the PKB window (solid black line) with
  degree~$\nu$ selected as in \eqref{eq:PKB-parameters},
  are used; results for other kernels are almost identical and
  not shown. In both (a) and (b), $D=3$ (other periodicities are the same
  assuming that other error sources are negligible), there are
  $N=10^3$ sources, $Q=10^3$, the box side length is $L=1$, the
  decomposition parameter is $\xi=10$, and the grid size is
  $M=46$ (stokeslet, rotlet) or $M=48$ (stresslet), corresponding
  to a relative truncation error below $10^{-23}$ (negligible).
  Note that the quantity on the $y$-axis is an approximation of
  the relative window approximation error $\Ewindowrel^\fp$.
  The reference potential has the same parameters but is computed
  using the exact KB window with $P=20$.}
  \label{fig:estimates-approximation}
\end{figure}

We will now characterize approximation errors introduced by the
window function, and will therefore for the moment assume that
the truncation errors described in
section~\ref{sec:estimates-truncation}, and any errors introduced
by additional approximations in the free directions, are negligible.
In \cite{Shamshirgar2021}, the (absolute) window approximation error
estimate
\begin{equation}
  \label{eq:full-approximation-estimate}
  \Ewindow^\fp \approx
  5 U(L,Q,\xi) \left(
  \E^{-2\pi P^2/\beta} + \erfc(\sqrt{\beta})
  \right)
\end{equation}
was constructed for the exact KB window \eqref{eq:KB-window}.
Here, $P$ is the window size and $\beta$ is the shape parameter;
furthermore, $U(L,Q,\xi)$ is an estimate of the rms value of the
Fourier-space part potential, cf.\
section~\ref{sec:estimates-potential} and in particular
\eqref{eq:pot-mag-est-stokeslet}--\eqref{eq:pot-mag-est-rotlet}.
In \eqref{eq:full-approximation-estimate}, the exponential term
approximates the error caused by discretizing the window function
and applying the trapezoidal rule to it in
\eqref{eq:SE-gathering-disc}; the $\erfc$ term approximates the
error caused by the discontinuity of the KB window
\eqref{eq:KB-window} at $\lvert r \rvert = a_w = h P/2$. It was
shown that the two terms are approximately balanced if $\beta
\approx \sqrt{2\pi} P \approx 2.5 P$, at which point
\eqref{eq:full-approximation-estimate} simplifies to
\begin{equation}
  \label{eq:approximation-estimate}
  \Ewindow^\fp \approx
  10 U(L,Q,\xi) \E^{-2.5 P}.
\end{equation}
This window approximation error is independent of both the kernel
and the periodicity, assuming that other error sources are
negligible. The estimate \eqref{eq:approximation-estimate} is
illustrated in Figure~\ref{fig:estimates-approximation}~(a), and
is in excellent agreement with the actual error. In our parameter
selection procedure, we can use \eqref{eq:approximation-estimate}
to compute the window size $P$ given $\xi$ and an error tolerance.

In the SE method, we do not use the exact KB window, but a
piecewise polynomial approximation of it, the PKB window (see
section~\ref{sec:se-window-functions}). In
Figure~\ref{fig:estimates-approximation}~(b), we show that when
the polynomial degree $\nu$ is selected according to
\eqref{eq:PKB-parameters}, the PKB window introduces no further
error compared to the exact KB window.

\subsection{Error pollution}
\label{sec:estimates-pollution}

\begin{figure}[htb]
  \centering
  \includegraphics[width=0.5\textwidth]{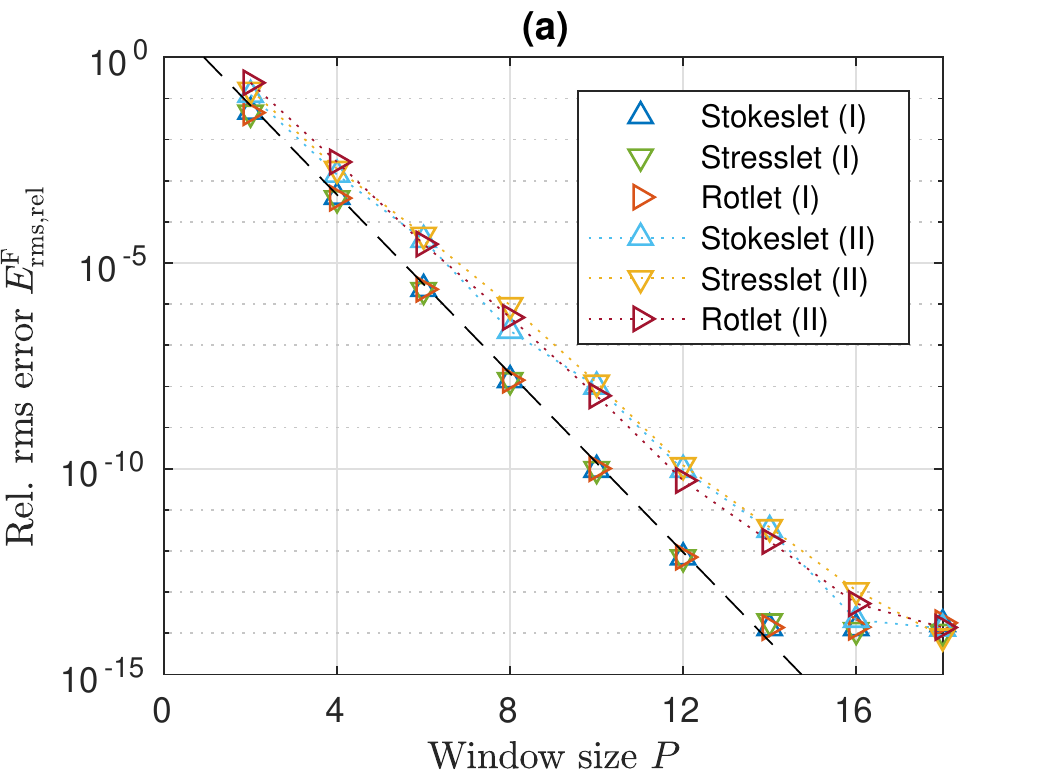}%
  \includegraphics[width=0.5\textwidth]{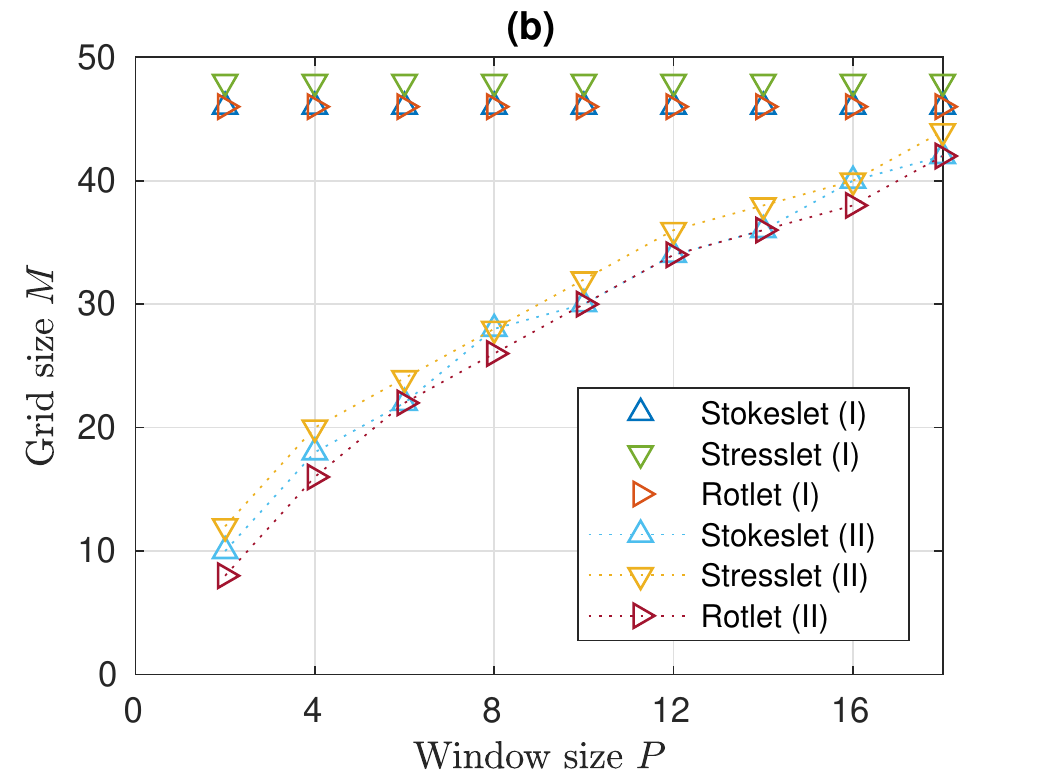}%
  \caption{Error pollution illustration.
  (a) Actual error for the three kernels when the PKB window
  is used; for (I), the grid size~$M$ is selected as in
  Figure~\ref{fig:estimates-approximation} such that the
  truncation error is negligible; for (II), the grid size~$M$
  is selected from estimates to make the truncation error
  approximately equal to the window approximation error for
  each $P$; the dashed black line is the estimate
  \eqref{eq:approximation-estimate}, i.e.\
  $E^\fp_\mathrm{rms,rel} = 10\E^{-2.5 P}$.
  (b) The grid sizes used to generate the data in (a). Other
  parameters are $D=3$, $N=10^3$, $Q=10^3$, $L=1$ and
  $\xi=10$, as in Figure~\ref{fig:estimates-approximation}.
  The reference potential is the same as (I) but computed
  with $P=20$.}
  \label{fig:estimates-pollution}
\end{figure}

So far, we have considered errors, such as the truncation error and
window approximation error, in isolation. It turns out that this
is not sufficient, as shown by Figure~\ref{fig:estimates-pollution}.
This figure has two sets of error curves: (I), for which the grid
size $M$ is selected large enough to make the truncation error
negligible (so that the window approximation error dominates),
and (II), for which $M$ is selected to make the truncation error
approximately of the same size as the window approximation error
for each given value of $P$
(according to the error estimates from
sections~\ref{sec:estimates-truncation}--\ref{sec:estimates-approximation}).
It is seen for (II) that the actual error is significantly larger
than one would expect from the error estimates (by simply summing
the individual error sources). From
Figure~\ref{fig:estimates-pollution}~(a), it is clear that the
convergence is still spectral in $P$ for (II), but with a
slower rate than expected from the window approximation error
estimate \eqref{eq:approximation-estimate}. Thus, if both $M$ and
$P$ are set from error estimates, the error tolerance will not be
met. We call this interdependence of the truncation and
approximation errors ``error pollution''.

We will not attempt an exhaustive analysis of the mechanics
behind error pollution here. However, we note that it is related
to aliasing errors that appear when evaluating the discrete
Fourier transform coefficients of the window function. Aliasing
causes the actual approximation error to depend not only on $P$,
but also on $M$, an effect that is responsible for the ``error pollution''
and is missing from the approximation error estimates found in
section~\ref{sec:estimates-approximation}.
A more detailed analysis of the approximation error would include
the effect of approximating the continuous Fourier transform of
the window function by the discrete Fourier transform, which
would take aliasing, and thus the ``error pollution'', into
account; this is left as future work. A related analysis of
aliasing errors in the nonuniform fast Fourier transform (NUFFT),
although not directly applicable to our situation, is found in
\cite{Potts2021,Barnett2021}.

Error pollution can be reduced either by increasing
$M$ (which reduces aliasing errors) or $P$ (which reduces the
approximation error, albeit slower than expected from the error
estimates).
In practice, the best strategy with respect to minimizing the run
time of the method is probably to increase both $M$ and $P$
moderately. We have found heuristically that if one either
increases $M$ to $1.05 M$ and $P$ to $P+4$, or $M$ to $1.1 M$ and
$P$ to $P+2$, compared to the values given by the error
estimates, the total error decreases to the level one would
expect from the estimates. Which of these rules is the fastest
will depend on the periodicity~$D$, since, as will be seen in
section~\ref{sec:estimates-aft}, increasing $P$ also increases
the extended grid size $\extM$ in the free directions, thus
increasing the cost of FFTs and the scaling step. It is thus more
expensive to increase $P$ the smaller $D$ is. We have found that
a suitable adjustment rule to mitigate error pollution is
\begin{equation}
  \label{eq:error-pollution-fix}
  (h_\mathrm{target}, P_\mathrm{target}) = \begin{cases}
    (h_\mathrm{estimate}/1.1, P_\mathrm{estimate}+2),
    & \text{for $D=0$}, \\
    (h_\mathrm{estimate}/1.05, P_\mathrm{estimate}+4),
    & \text{for $D=1,2,3$},
  \end{cases}
\end{equation}
where, as usual, $h=L/M$. To reiterate, the dependence on $D$ is
only to try to minimize the overall run time for each $D$.
The rule \eqref{eq:error-pollution-fix} is used whenever both $h$
and $P$ are selected from the error estimates, unless otherwise
stated. It is verified together with the rest of the parameter
selection procedure in section~\ref{sec:estimates-summary}.

\subsection{Box padding and upsampling}
\label{sec:estimates-aft}

The estimates covered in
sections~\ref{sec:estimates-truncation}--\ref{sec:estimates-pollution}
are enough to select all parameters that appear in the fully
periodic ($D=3$) case. When some, or all, spatial directions are
free, i.e.\ in the $D=2,1,0$ cases, additional parameters appear
in the method: the box paddings $\delta L_i$, upsampling factors
$s_0$ and $s_*$, and upsampling threshold $\isc{k}_i^*$. Since we
assume here that the primary cell is a cube, the parameters
will be the same in all directions, i.e.\ $\delta L_i = \delta L$
and $\isc{k}_i^* = \isc{k}_*$, $i=1,2,3$. The parameters that we
discuss here are to some extent interdependent, since, as seen by
\eqref{eq:free-disc-wavenumber-resolution}, both the upsampling
factor $s$ (i.e.\ $s_0$ or $s_*$) and the box padding $\delta L$
(which appears in the side length $\extL = L + \delta L$
of the extended box) will influence the wavenumber resolution in
the free directions. Thus, if $\delta L$ is selected larger, $s$
can be selected smaller, and vice versa. Below, we first describe
how $\delta L$ is selected, and then how the upsampling
parameters $s_0$, $s_*$ and $\isc{k}_*$ are selected, which is
also the order that the parameters are selected in the final parameter
selection algorithm.

The box padding $\delta L$ relates the side lengths of the
original box $L$ and the extended box $\extL = L + \delta
L$, in the free directions.
Recall from section~\ref{sec:se-discrete} that the reason that
the original box must be padded is that (i) the window
function must be fully contained in the extended box, and (ii)
the screening function must have decayed sufficiently at the
boundary of the extended box. Point (i) is guaranteed by
requiring that $\delta L \geq 2 a_w = hP$, where $a_w$ is the
halfwidth of the window function (see Figure~\ref{fig:window}).
Point (ii) is less straightforward, since the screening function
does not truly have compact support, so it must be truncated at
some level where it has decayed ``sufficiently''. We have found
it useful to let $\delta L \approx \lambda hP$, where the factor
$\lambda \geq 1$ depends on the kernel (and thus on the screening
function); this naturally connects points (i) and (ii), and
furthermore links $\delta L$ to the window size $P$, and thus to
the overall precision of the method (so that $\delta L$ becomes larger when more
precision is required). Since the Hasimoto screening function
\eqref{eq:hasimoto-screening} decays slower than the Ewald
screening function \eqref{eq:ewald-screening}, $\lambda$ is
expected to be larger for the stokeslet and stresslet, than for
the rotlet.

In practice, we want the extended grid size $\extM =
\extL/h$ to be a multiple of $\gridfactor$, just like $M$
is (cf.\ \eqref{eq:grid-adjustment}). Therefore, we first set
$\extM$, and then compute $\extL = h \extM$
(and hence $\delta L = \extL-L$), rather than the other
way around. The extended grid size is set to
\begin{equation}
  \label{eq:extended-grid-size}
  \extM = \gridfactor \left\lceil
  \frac{L/h + P + (\lambda-1)\max(P,\theta)}{\gridfactor}
  \right\rceil,
\end{equation}
where $\gridfactor$ is some positive integer and $\lceil \cdot
\rceil$ is the ceiling function; here, $\theta$ is a safety
threshold and $\theta=8$ for the stokeslet and stresslet, while
$\theta=0$ for the rotlet. Note that
\eqref{eq:extended-grid-size} is approximately the same as
setting $\delta L/h \approx \lambda P$. It remains to select
$\lambda$, but we will postpone that for a few paragraphs until
the upsampling parameters have been discussed.

Recall from section~\ref{sec:se-aft} that upsampling is needed
both to resolve the modified kernels used for the
$\bvec{k}^\per=\bvec{0}$ mode, which is upsampled by the factor
$s_0$, and to resolve the regular kernels for modes
$\bvec{k}^\per \in \wavenumset_*$, given by \eqref{eq:wavenumset-upsampling},
which are upsampled by the factor $s_*$. The parameter $\isc{k}_*$
controls the size of the set $\wavenumset_*$. It has been
established by \cite{afKlinteberg2017} that $s_0$ should be
selected according to \eqref{eq:upsampling-factor-s0}, where $i$
in the minimization is restricted to the free directions ($D+1
\leq i \leq 3$). In practice, we round $s_0$ upwards to one
decimal, i.e.\ we use $\lceil 10 s_0 \rceil / 10$.

Ref.~\cite{Shamshirgar2021} determined a rule for selecting $s_*$ and
$\isc{k}_*$ for the harmonic kernel, namely
\begin{align}
  \label{eq:upsampling-factor-sstar}
  s_* &=
  \frac{M}{\extM} \left( 1 + \frac{1}{2\pi}
  \log\left(
  \frac{U(L,Q,\xi)}{2\Erms^\fp}
  \right) \right),
  \\
  \label{eq:estimate-kbarstar}
  \isc{k}_* &= \left\lceil
  \frac{M}{(\extM - M)} \frac{1}{2\pi}
  \log\left(\frac{U(L,Q,\xi)}{2\Erms^\fp}\right) - 1
  \right\rceil,
\end{align}
where $\Erms^\fp$ is the absolute rms error and $U(L,Q,\xi)$ is
an estimate of the rms value of the Fourier-space part potential.
It turns out that \eqref{eq:upsampling-factor-sstar}--\eqref{eq:estimate-kbarstar}
work well also for the stokeslet, stresslet and rotlet, using the
appropriate $U(L,Q,\xi)$ for each kernel according to
\eqref{eq:pot-mag-est-stokeslet}--\eqref{eq:pot-mag-est-rotlet}.
Finally, we want the upsampled grid sizes $s_0 \extM$ and
$s_* \extM$ to be multiples of $\gridfactor$, and will
therefore adjust both $s_0$ and $s_*$ upwards such that this
holds.

Note that the formulas \eqref{eq:upsampling-factor-s0} and
\eqref{eq:upsampling-factor-sstar}--\eqref{eq:estimate-kbarstar}
to select $s_0$, $s_*$ and $\isc{k}_*$ depend on the extended box
size $\extL$ and grid size $\extM$, which makes the
upsampling parameters depend on the box padding $\delta L$. The
only missing piece of our parameter selection algorithm is the
value of $\lambda$ which is needed to compute $\extM$, and
thus $\delta L$, from \eqref{eq:extended-grid-size}. We determine
$\lambda$ manually by numerical experiments, by considering a
plot such as Figure~\ref{fig:estimates-approximation}~(a),
initializing $\lambda$ to 1 (which makes the errors large) and
then gradually increasing $\lambda$ until the errors follow the
estimates as in Figure~\ref{fig:estimates-approximation}. At
every point of this process, the upsampling parameters are
selected according to \eqref{eq:upsampling-factor-s0} and
\eqref{eq:upsampling-factor-sstar}--\eqref{eq:estimate-kbarstar}.

The resulting values of $\lambda$ are given in
Table~\ref{tab:estimates-lambda}.
For $D=0$, the smallest value of $\lambda$ given by the process
just described is used; as expected, $\lambda$ is larger for the
stokeslet and stresslet than for the rotlet, due to the different
screening functions used for these kernels. (The stokeslet and
stresslet also have slightly different $\lambda$ since the
relation between $\xi$, $h$ and $P$ is different for the
different kernels.) For $D=1,2$, we select a larger $\lambda$
than strictly necessary, namely $\lambda=2.4$ for all kernels.
This allows us to use slightly smaller values for $s_*$ and
$\isc{k}_*$ (since they are inversely related to
$\extM$ and thus to $\lambda$, cf.\
\eqref{eq:upsampling-factor-sstar}--\eqref{eq:estimate-kbarstar}
and \eqref{eq:extended-grid-size}), which reduces the overall
computation time of the method slightly.

\begin{table}[h]
  \centering
  \caption{Value of $\lambda$ for different kernels and
  periodicities.}
  \label{tab:estimates-lambda}
  \begin{tabular}{ccc}
    \toprule
    & $D=0$ & $D=1,2$ \\
    \midrule
    Stokeslet & $\lambda=2.2$ & $\lambda=2.4$ \\
    Stresslet & $\lambda=2.4$ & $\lambda=2.4$ \\
    Rotlet    & $\lambda=1.5$ & $\lambda=2.4$ \\
    \bottomrule
  \end{tabular}
\end{table}

\subsection{Summary of the parameter selection procedure}
\label{sec:estimates-summary}

For a given kernel (stokeslet, stresslet or rotlet) and
periodicity~$D$ (0, 1, 2 or 3), the automated parameter selection
procedure can be summarized as follows. The input consists of the
primary cell $\boxvar =
[0,L)^3$, the source strength quantity $Q$ as given by
\eqref{eq:vector-source-quantity} or
\eqref{eq:tensor-source-quantity}, the decomposition parameter
$\xi$, and an absolute error tolerance $\abstol$.
The general structure of the procedure is the same as in
\cite{Shamshirgar2021}, with the following steps:
\begin{enumerate}
  \item
    Given $L$, $Q$ and $\xi$, compute the real-space part cut-off
    radius $\rc$ from estimate \eqref{eq:rs-trunc-est-stokeslet},
    \eqref{eq:rs-trunc-est-stresslet} or
    \eqref{eq:rs-trunc-est-rotlet}, depending on the kernel,
    with $\Etrunc^\rp = \abstol$. (This is the only step needed
    for the real-space part.)
  \item
    Compute a preliminary grid spacing $h_\mathrm{estimate}=\pi/k_\infty$ by
    computing $k_\infty$ from one of the estimates
    \eqref{eq:fs-trunc-est-rotlet},
    \eqref{eq:fs-trunc-est-stokeslet} or \eqref{eq:fs-trunc-est-stresslet},
    depending on the kernel, with $\Etrunc^\fp = \abstol$.
  \item
    Compute a preliminary window size $P_\mathrm{estimate}$ from the estimate
    \eqref{eq:approximation-estimate}, with $\Ewindow^\fp =
    \abstol$, and with $U(L,Q,\xi)$ for the kernel in question
    given by
    \eqref{eq:pot-mag-est-stokeslet}--\eqref{eq:pot-mag-est-rotlet};
    the constant $C_U$ is given by its mean value from
    Table~\ref{tab:estimates-potential}.
  \item
    Apply the error pollution adjustment rule
    \eqref{eq:error-pollution-fix} to adjust the preliminary
    values $h_\mathrm{estimate}$ and $P_\mathrm{estimate}$,
    yielding $h_\mathrm{target}$ and $P_\mathrm{target}$. Then
    apply rounding according to \eqref{eq:grid-adjustment}, with
    $\gridfactor$ set to a small power of two, to
    compute the actual grid spacing $h_\mathrm{actual}$. This
    also yields the grid size $M_\mathrm{actual} =
    L/h_\mathrm{actual}$ in each periodic direction.
  \item
    Compute $P_\mathrm{actual} = 2 \lceil P_\mathrm{target}/2 \rceil$,
    then compute the shape parameter $\beta$ and polynomial
    degree $\nu$ from \eqref{eq:PKB-parameters}. If $D=3$, stop
    here; the parameter selection is complete.
  \item
    Compute the extended grid size $\extM$ in the free
    directions from \eqref{eq:extended-grid-size}, with $\lambda$
    given by Table~\ref{tab:estimates-lambda}, and using
    $\gridfactor$ from step~4. Compute
    $\extL = h \extM$. This also sets $\delta L =
    \extL - L$ and the truncation radius $R$ through
    \eqref{eq:truncation-radius}.
  \item
    Compute the upsampling parameters $s_0$ from
    \eqref{eq:upsampling-factor-s0} with $D+1 \leq i \leq 3$,
    $s_*$ from \eqref{eq:upsampling-factor-sstar}, and
    $\isc{k}_*$ from \eqref{eq:estimate-kbarstar}. Here,
    $\Erms^\fp = \abstol$, and $U(L,Q,\xi)$ is as in step~3 above.
    Round $s_0$ to $\lceil 10 s_0 \rceil/10$, then adjust both upsampling
    factors such that $s_0 \extM$ and $s_* \extM$
    are multiples of $\gridfactor$.
\end{enumerate}
The grid sizes in steps 4, 6 and 7 are consistently rounded
upwards to multiples of $\gridfactor$, with $\gridfactor$ set to
a small power of two in order to increase the efficiency of FFTs.
Typically, a small value for $\gridfactor$ (such as 2) will put
the error of the SE method close to the target error tolerance,
while a slightly larger value (such as 4 or 8) may lead to more
consistent run times but give a smaller error than expected.
Here, in section~\ref{sec:estimates}, we set $\gridfactor=2$ in
the numerical examples, while in section~\ref{sec:results} we
will set $\gridfactor=4$. (The default value in the SE package is
$\gridfactor=4$, but the optimal value may depend on the FFT
implementation.)

\begin{figure}[htp]
  \centering
  \includegraphics[width=0.33\textwidth]{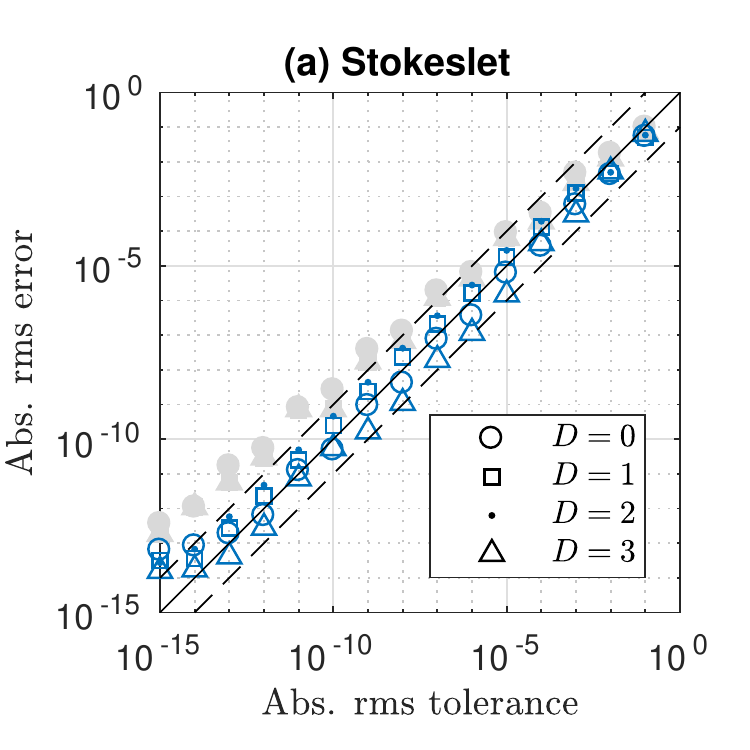}%
  \includegraphics[width=0.33\textwidth]{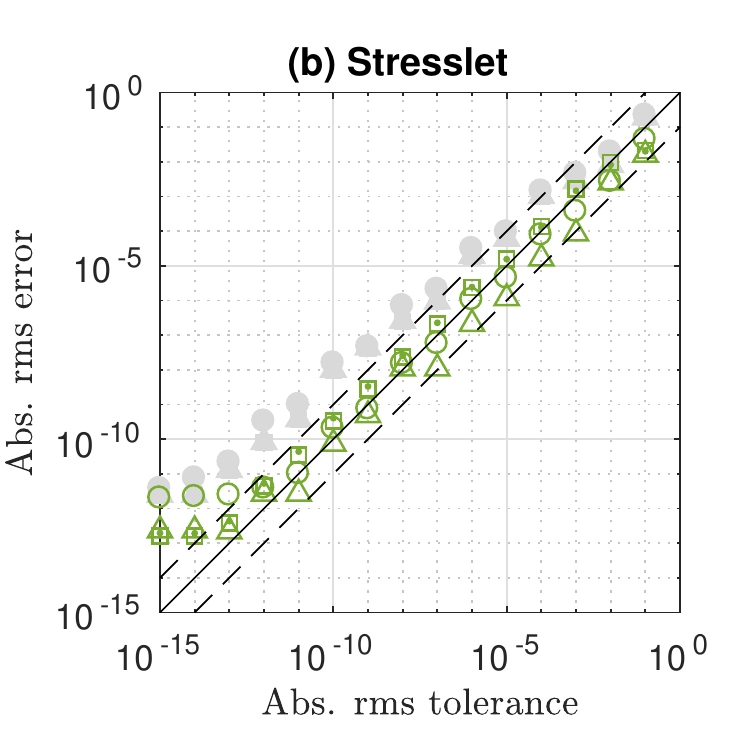}%
  \includegraphics[width=0.33\textwidth]{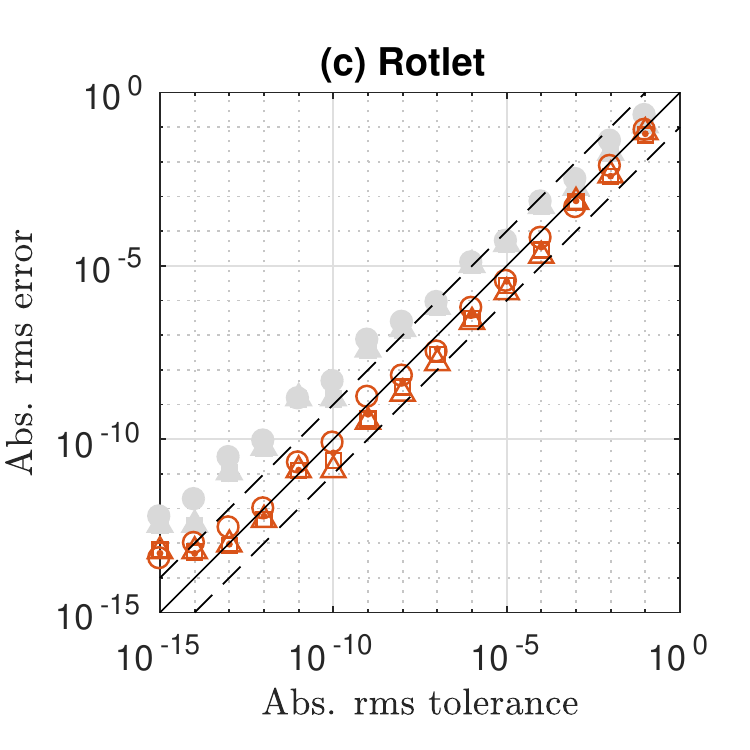}\\
  \includegraphics[width=0.33\textwidth]{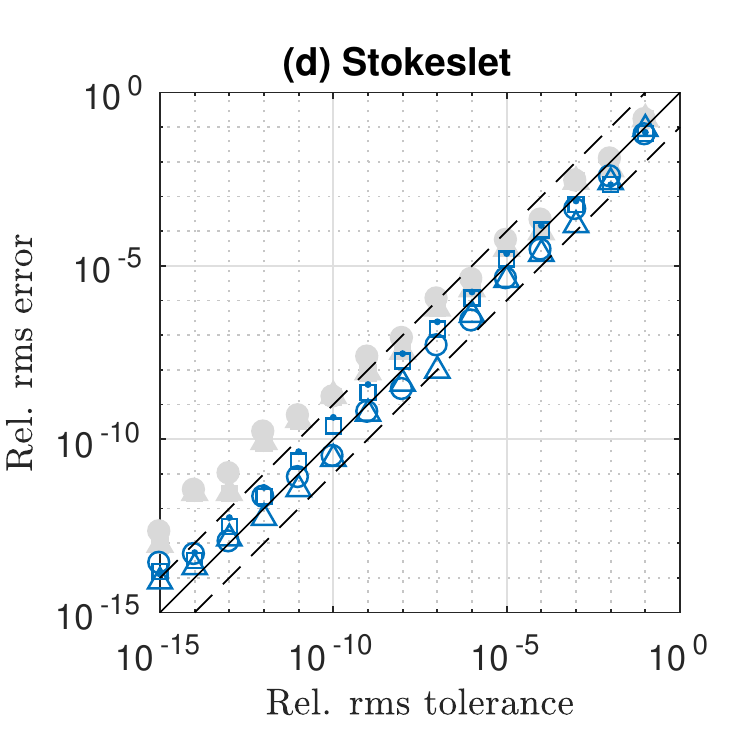}%
  \includegraphics[width=0.33\textwidth]{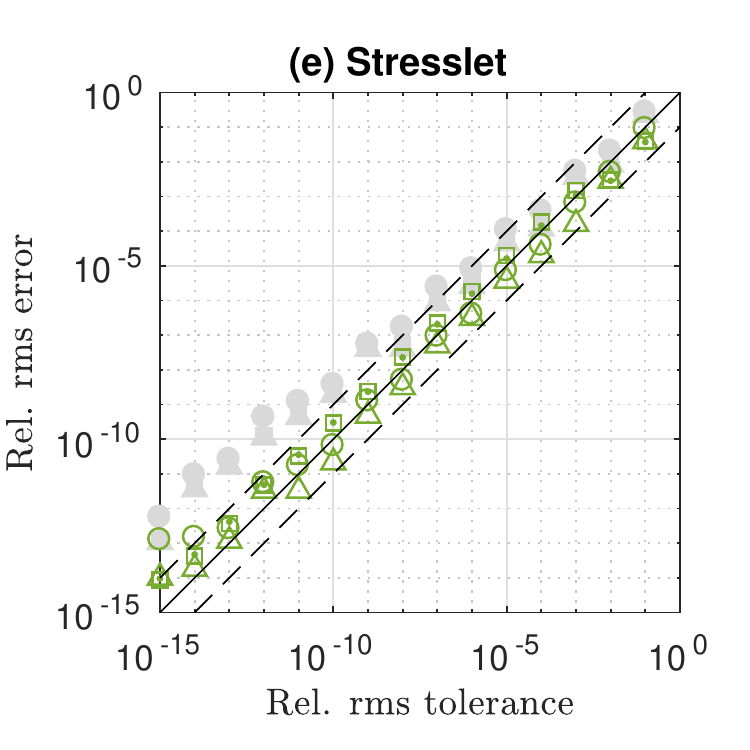}%
  \includegraphics[width=0.33\textwidth]{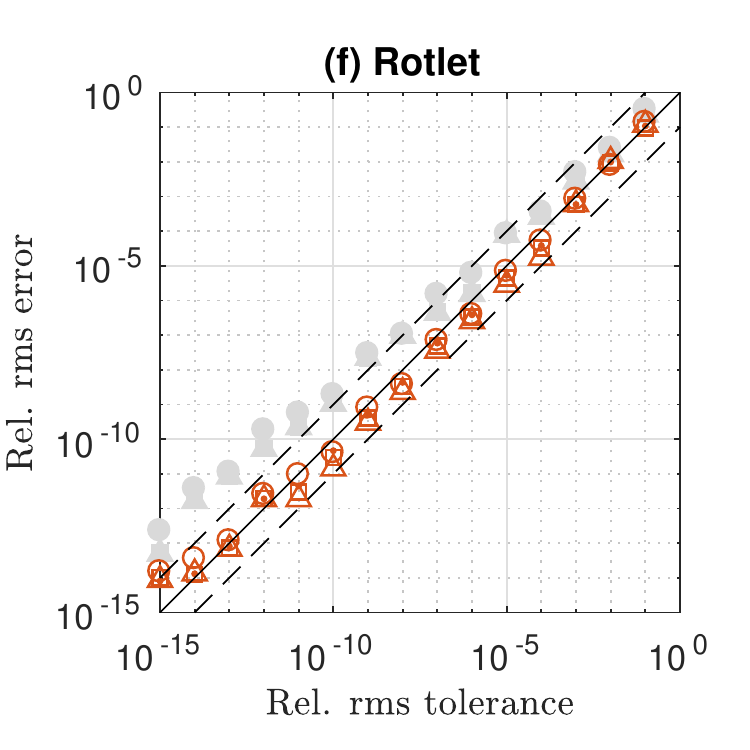}%
  \caption{Parameter selection test, part 1. Actual rms error is
  plotted versus the error tolerance; the top row shows absolute
  errors while the bottom row shows relative errors; each column
  shows one kernel. Colored unfilled symbols (blue, green, red)
  are with the error pollution adjustment \eqref{eq:error-pollution-fix},
  light gray filled symbols are without. The shape of the symbols shows
  the periodicity~$D$, according to the legend. The solid black
  line shows error = tolerance, while the dashed lines are offset
  by a factor of 10 in each direction. Parameters are $N=1000$, $L=1$, $Q=1$,
  $\xi=10$; all other parameters are set by the parameter
  selection procedure. The reference potential is computed using
  the same method with a relative error tolerance of $10^{-17}$.
  For this case, the rms value of the potential is around 2 for
  the stokeslet, 20 for the stresslet, and 7 for the rotlet,
  also captured well by
  \eqref{eq:pot-mag-est-stokeslet}--\eqref{eq:pot-mag-est-rotlet}.}
  \label{fig:validation-part1}
\end{figure}

To verify that the parameter selection procedure works, we will
consider some examples. Here, we consider only the error in the
Fourier-space part of the potential; the full potential will be
considered in section~\ref{sec:results}. First, we consider a single
particle system per kernel, consisting of $N=1000$ sources in a
box of side length $L=1$; the sources are also the evaluation
points. The sources are random and selected
such that $Q=1$. We pick $\xi=10$ and set the absolute tolerance
$\abstol$ to different values between $10^{-1}$ and $10^{-15}$.
For each kernel and periodicity, the actual Fourier-space part
error is plotted versus the tolerance in
Figure~\ref{fig:validation-part1}~(a)--(c).
We show that the error pollution adjustment
\eqref{eq:error-pollution-fix} is needed by plotting data both
when it is used (colored unfilled symbols) and when it is not
used (light gray filled symbols); the error in the latter case is
too large. With the error pollution adjustment, the error is within one order of
magnitude of the tolerance, as desired, until it flattens out.

Using the potential rms value estimate $U(L,Q,\xi)$, it is also
possible to prescribe a relative error tolerance $\reltol$ for
the Fourier-space part, using the relation $\abstol = U(L,Q,\xi)
\reltol$. In Figure~\ref{fig:validation-part1}~(d)--(f), we plot
the relative error versus the relative tolerance. Again, the
error is within one order of magnitude of the tolerance when the
error pollution adjustment is used, which shows that $U(L,Q,\xi)$
is a good estimate of the potential rms value.

\begin{figure}[ht]
  \centering
  \includegraphics[width=0.33\textwidth,trim=65mm 117mm 72mm 116mm,clip]{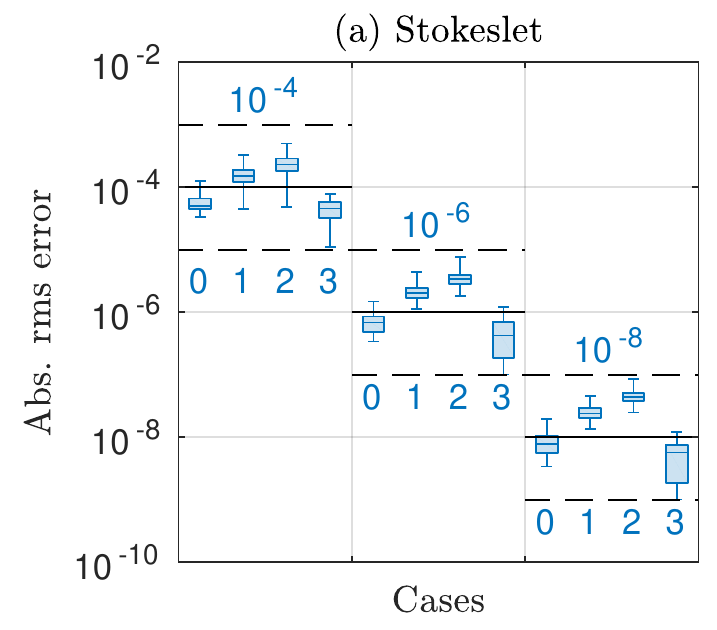}%
  \includegraphics[width=0.33\textwidth,trim=65mm 117mm 72mm 116mm,clip]{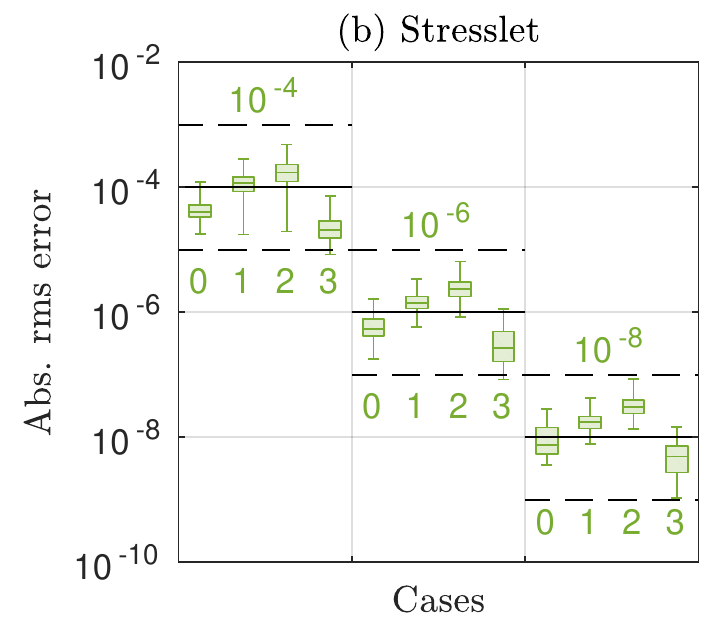}%
  \includegraphics[width=0.33\textwidth,trim=65mm 117mm 72mm 116mm,clip]{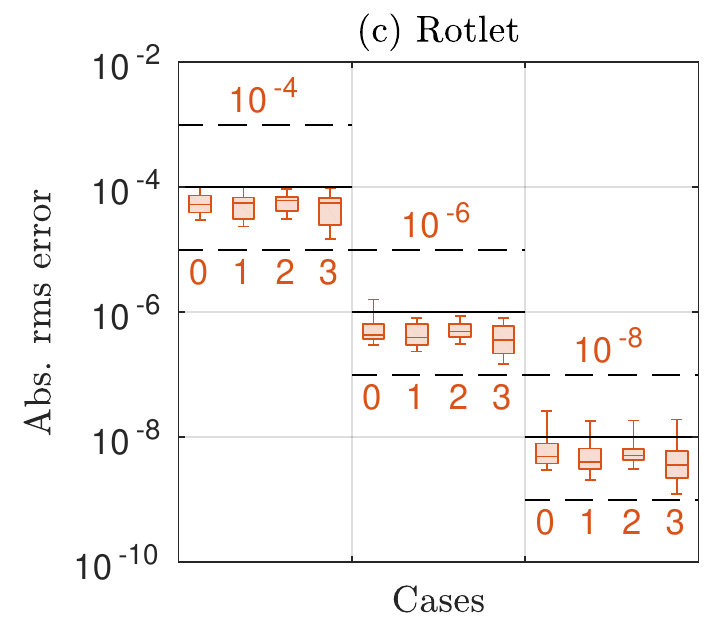}\\[\baselineskip]
  \includegraphics[width=0.33\textwidth,trim=65mm 117mm 72mm 116mm,clip]{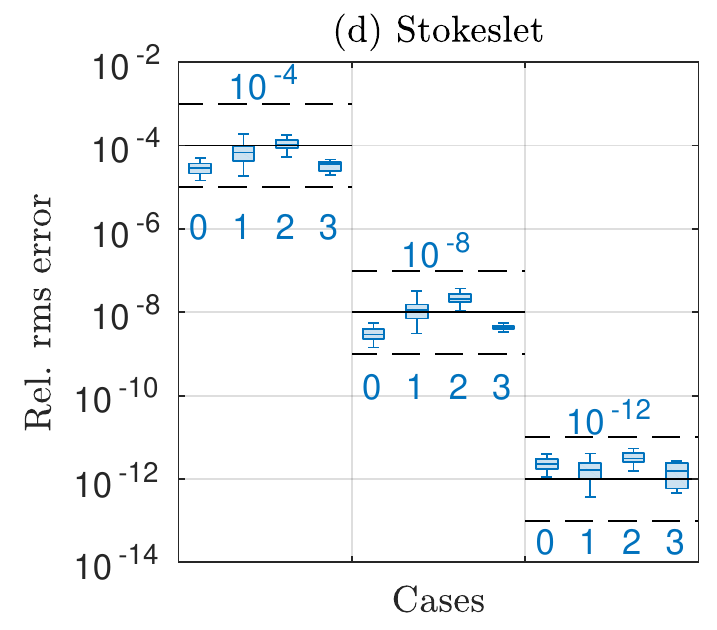}%
  \includegraphics[width=0.33\textwidth,trim=65mm 117mm 72mm 116mm,clip]{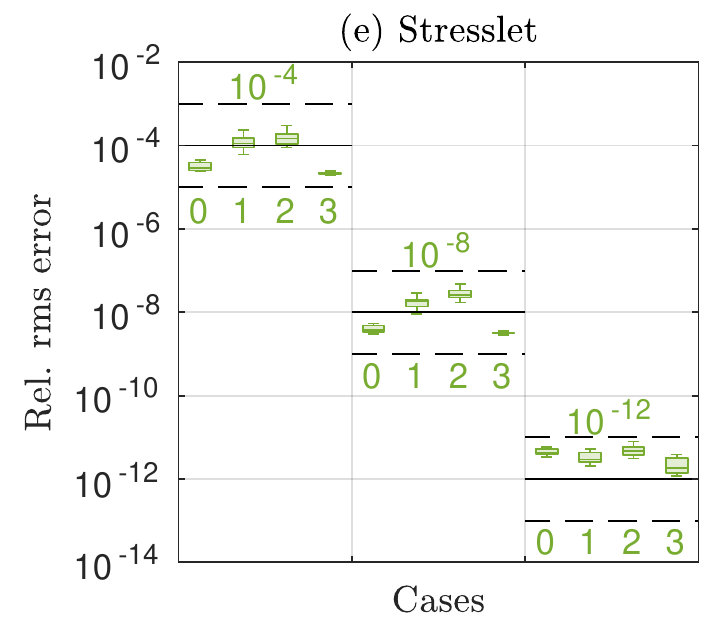}%
  \includegraphics[width=0.33\textwidth,trim=65mm 117mm 72mm 116mm,clip]{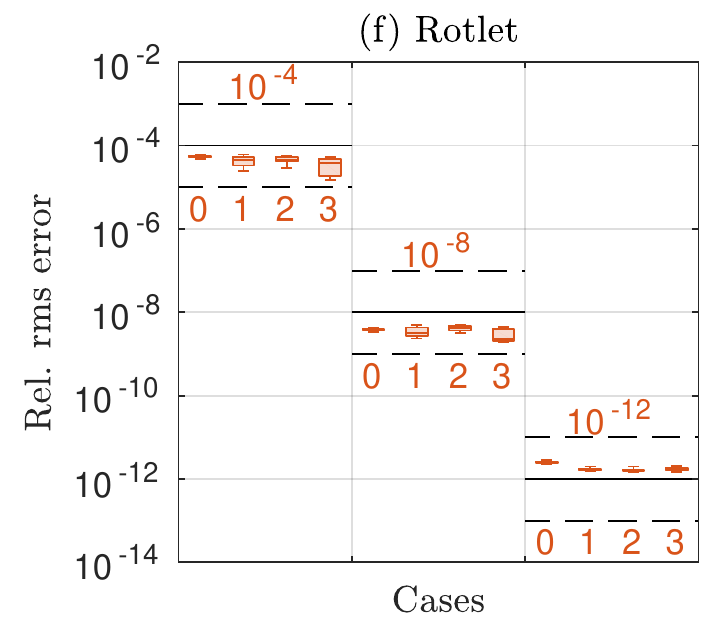}%
  \caption{Parameter selection test, part 2. Box plots of actual
  rms error for different tolerances and periodicities
  $D=0,1,2,3$ (indicated inside the plots). Each box plot is
  constructed using 270 samples (generated by every combination
  of $L=0.1, 1, 10$; $Q=0.01, 1, 100$; 10 different random seeds;
  and $\xi L = 10,20,30$). Boxes indicate upper/lower quartiles
  and the median, while whiskers indicate the maximum and minimum
  of the samples. The top row shows absolute errors for absolute error
  tolerances $10^{-4}$, $10^{-6}$ and $10^{-8}$; the bottom row
  shows relative errors for relative error tolerances $10^{-4}$,
  $10^{-8}$ and $10^{-12}$; the tolerance is marked with a solid
  black line. The error pollution adjustment
  \eqref{eq:error-pollution-fix} is always used. Throughout, the number of
  source points is $N=1000$. Other parameters than those already
  mentioned are set by the parameter selection procedure. The
  reference potential is computed using the same method with a
  relative error tolerance of $10^{-17}$. The rms value of the
  potential varies with $L$, $Q$ and $\xi$, and ranges from 0.02
  to 1000 for the stokeslet, from 0.006 to 10\,000 for the
  rotlet, and from 0.02 to 40\,000 for the stresslet.}
  \label{fig:validation-part2}
\end{figure}

To further test the robustness of the procedure, we consider 90
different particle systems per kernel, given by all combinations
of $L=0.1, 1, 10$; $Q=0.01, 1, 100$; and 10 different random seeds.
For each system, we set $\xi$ such that  $\xi L = 10, 20, 30$ and
compute the actual Fourier-space part error for a few different tolerances (again
evaluating the potential at the source points). The result is shown in
Figure~\ref{fig:validation-part2}. In the vast majority of cases, the error is
within one order of magnitude of the tolerance, the only exception
being that the error in the $D=3$ case sometimes becomes ``too
small'', i.e.\ smaller than one tenth of the tolerance. This is
of course not a problem, so we conclude that the parameter
selection procedure works well.

\section{Numerical results}
\label{sec:results}

Having established an automated procedure for selecting the
parameters of the SE method, we move on to other numerical
results. We start by considering the pointwise error in
section~\ref{sec:results-pointwise}, then measure the
computational time of the method in
section~\ref{sec:results-timing}, and finally in
section~\ref{sec:results-window} compare the PKB
window function with the truncated Gaussian window that appeared
in previous iterations of the SE method.

All computations are done on 6 cores on a machine with an Intel
Core i7-8700 processor, running at 4.6~GHz with 32~GB of memory.
The SE method is implemented in \Matlab, with critical routines
written in C and called through \Matlab's MEX interface. Fourier
transforms are computed using \Matlab's \texttt{fft} function, which is
based on the FFTW library \citep{Frigo2005}.
Throughout section~\ref{sec:results}, we let the grid sizes be
multiples of $\gridfactor=4$, since we have observed that this
leads to more consistent timings than $\gridfactor=2$. The reason
is that the FFT can be substantially faster when the grid size is
a multiple of 4.

\subsection{Pointwise error}
\label{sec:results-pointwise}

\begin{figure}[htbp]
  \centering
  \includegraphics[width=0.33\textwidth]{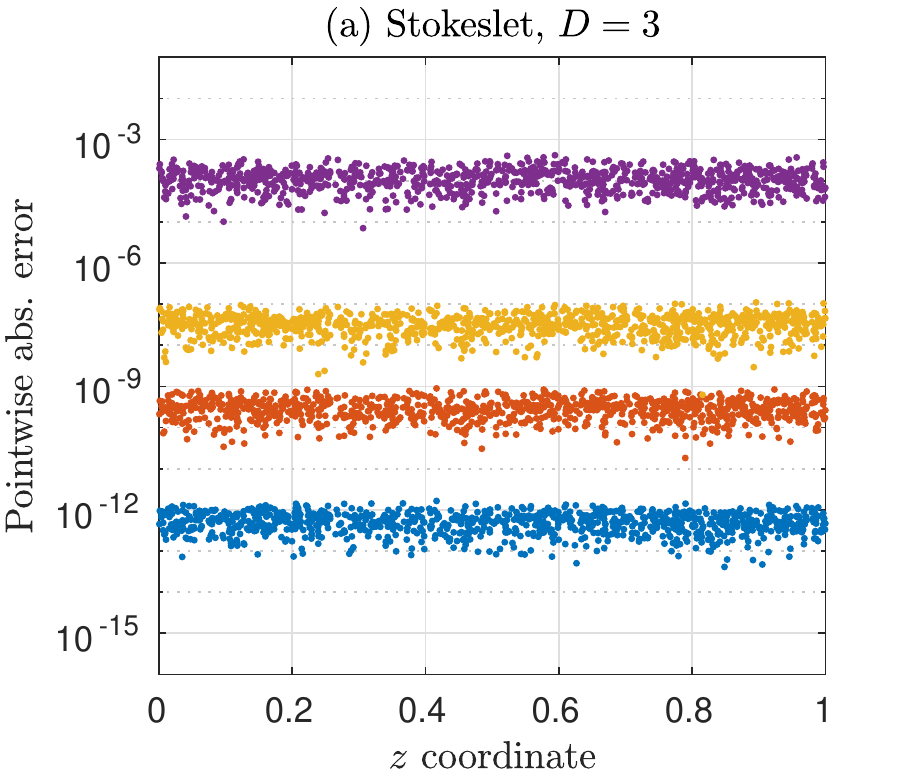}%
  \includegraphics[width=0.33\textwidth]{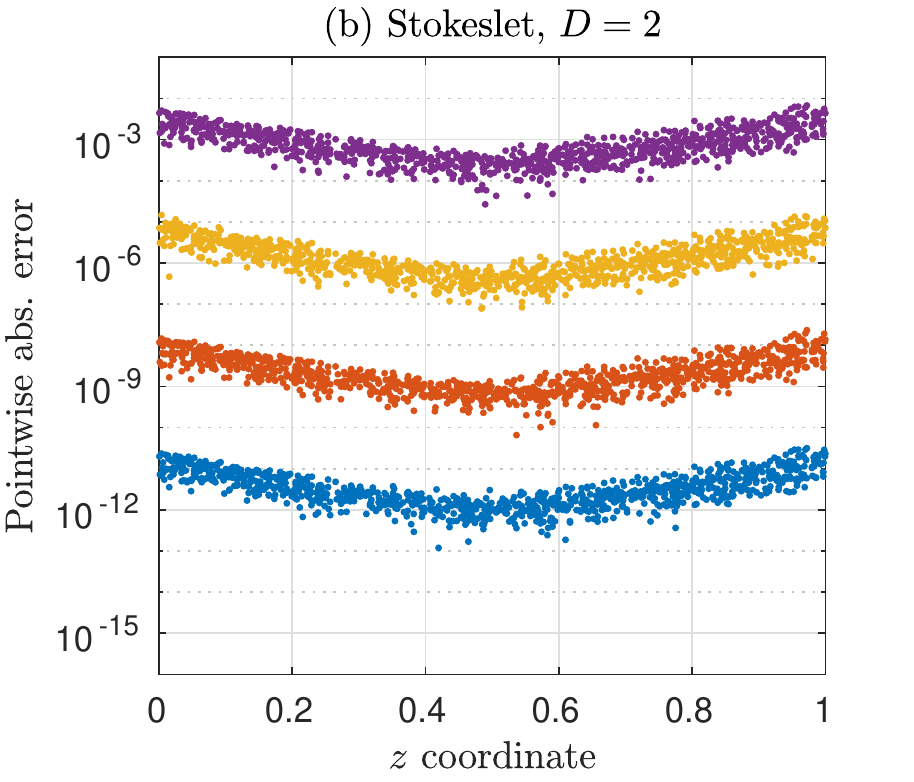}%
  \includegraphics[width=0.33\textwidth]{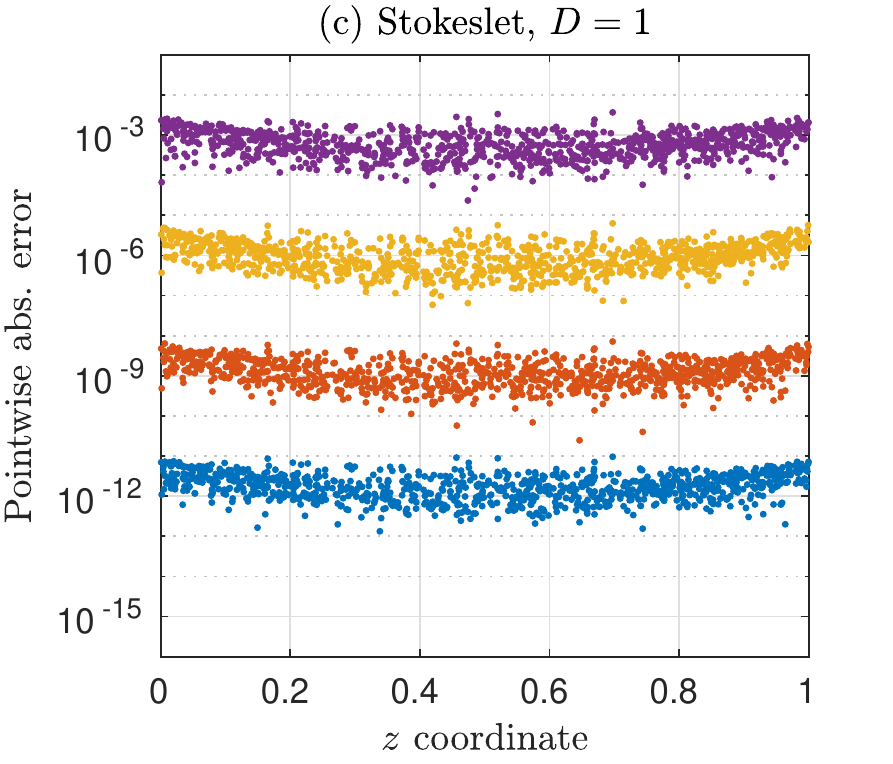}%
  \caption{Pointwise error. The Fourier-space part pointwise
  error for the stokeslet is plotted as a function of the $z$
  coordinate, for absolute rms error tolerances
  $\abstol=10^{-12}$ (blue), $10^{-9}$ (red), $10^{-6}$ (yellow),
  $10^{-3}$ (purple).
  (For $D=2$, the $z$ direction is the free direction, while for
  $D=1$, both the $y$ and $z$ directions are free.) For $D=0$,
  not shown here, the behaviour is qualitatively the same as in
  the $D=3$ case (also for $D=2,1$, the spatial dependence seems
  to come only from the $\bvec{k}^\per \neq \bvec{0}$ modes).
  Other kernels have similar behaviour as the
  stokeslet. The data consists of 10 different random
  systems with $N=100$ points each (source points are also
  evaluation points), superimposed. The primary cell is a cube with side
  length $L=1$; furthermore, $Q=1$, $\xi=10$. Other parameters are set
  from the error tolerances. The reference potential is computed
  using direct summation of the Fourier-space part Ewald sums
  \eqref{eq:fourier-potential-Dp}--\eqref{eq:1p-evaluated-integrals-k0},
  truncated such that the truncation error is $10^{-17}$
  ($k_\infty = 40\pi$).}
  \label{fig:pointwise-error}
\end{figure}

So far we have considered only the rms error \eqref{eq:abs-rms-error},
and the reader may wonder what can be said about the pointwise
error $\lvert \bvec{u}_h(\bvec{x}_m) -
\bvec{u}_\mathrm{ref}(\bvec{x}_m) \rvert$ at the evaluation
points $\bvec{x}_m$. For the harmonic kernel in $D=2$,
\cite{Lindbo2012} showed that the pointwise error has a weak
dependence on the coordinate in the free direction. The same
phenomenon can be seen for the kernels considered here, as shown
in Figure~\ref{fig:pointwise-error} for the stokeslet
(Fourier-space part). For $D=3$, the error has the same size in
the whole primary cell, but for $D=2,1$, the error is slightly larger at
the edges of the cell in the free directions, than at the middle
of the cell. However, the variation is not large, and the rms
error is still within one order of magnitude of the tolerance (as
already shown in section~\ref{sec:estimates-summary}).

\subsection{Timing results}
\label{sec:results-timing}

We will here study how the computational time of the method as a
whole (Fourier-space and real-space parts, i.e.\ the full
potential) varies with the number of points $N$ in the primary
cell and the error tolerance $\abstol$. We fix the box side length $L=1$ and source
strength quantity $Q=1$, since these do not affect the
computational time. The potential is evaluated at the source
points, i.e.\ $N_\mathrm{t} = N$. The precomputation step for the
$D=0$ case, described in section~\ref{sec:se-precomp}, is
excluded from the timing since it can be reused for multiple
solves.

First, we fix $\abstol=10^{-8}$ and vary the number of points
$N$, i.e.\ we measure the time complexity of the method. This is
done for two different types of random systems, shown in
\begin{figure}[b]
  \centering
  \includegraphics[width=0.4\textwidth]{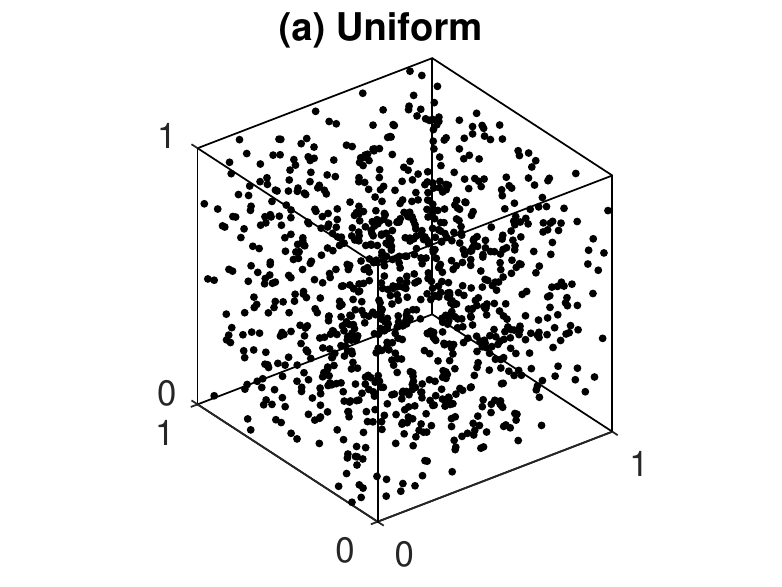}%
  \includegraphics[width=0.4\textwidth]{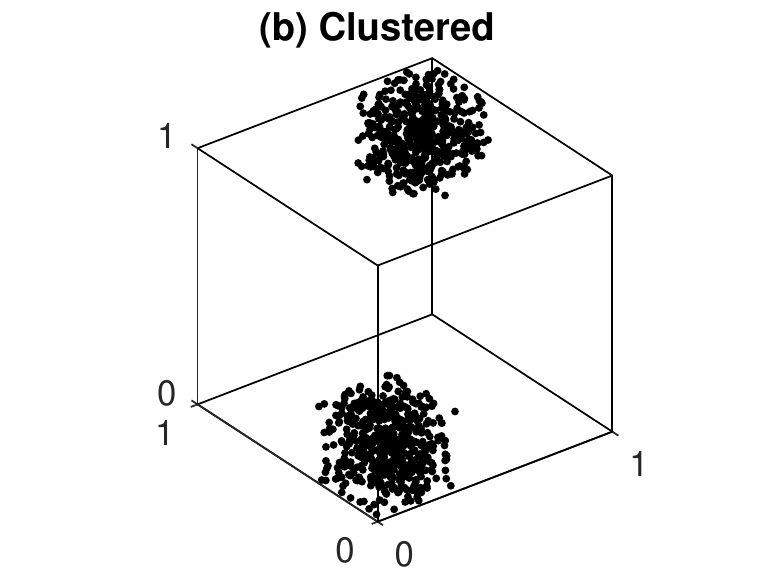}%
  \caption{(a) Uniform random system; (b) Clustered random
  system.}
  \label{fig:random-systems}
\end{figure}
Figure~\ref{fig:random-systems}: the uniform system that we have
already used above, and a clustered system consisting of denser
point clouds (with points only present in the subdomains
$[0,1/3)^3$ and $[2/3,1)^3$).
As $N$ increases, the point concentration $N/L^3$ increases, and
the cut-off radius $\rc$ must be adjusted to prevent the real-space
part from scaling like $O(N^2)$. We set $\rc$ such that the
expected average number of points within a ball of radius $\rc$,
i.e.\ $N_{\rc} := \frac{4}{3}\pi \rc^3 N/L^3$, is kept fixed as $N$
increases. This way, the time of the real-space part will scale
like $O(N)$. Given $\rc$, we determine $\xi$ from the truncation
error estimates
\eqref{eq:rs-trunc-est-stokeslet}--\eqref{eq:rs-trunc-est-rotlet},
with $\Etrunc^\rp=\abstol=10^{-8}$, and all remaining parameters
are then set according to the parameter selection procedure.
Roughly, we will have $\xi \sim \rc^{-1} \sim N^{1/3}$, and
furthermore the grid size $M \sim k_\infty \sim \xi$. The
gridding and gathering steps scale like $O(N)$, while the FFTs
scale like $O(M^3 \log M^3)$, and the scaling step like $O(M^3)$.
Since $M \sim N^{1/3}$, the Fourier-space part, and therefore the
method as a whole, will scale like $O(N \log N)$.

The expected average number $N_{\rc}$ of particles within a ball
of radius $\rc$ (from which $\xi$ is computed) is set to balance
the real-space part and Fourier-space part run times, so that the
total computation time is close to minimal. The
optimal value depends on the particle system (uniform or clustered), as
well as the periodicity and kernel (and the computer), and is
given (approximately) in Table~\ref{tab:complexity-optimal-pts}. Note that
$N_{\rc}$ is simply the volume of an $\rc$-ball times the
average particle concentration $N/L^3$ in the whole primary cell
(including also the empty regions for the clustered system;
$N/L^3$ has the same value for the uniform and clustered systems).
The fact that $N_{\rc}$ is selected smaller for the clustered
system than for the uniform one reflects that $\rc$ should be
smaller, since the particle clouds are denser. The reason that
$N_{\rc}$ is selected larger for lower periodicities~$D$ is that
the Fourier-space part becomes more expensive, so more
computational effort is shifted into the real-space part.

\begin{table}[htp]
  \centering
  \caption{Expected average number of points within a ball
  of radius $\rc$, i.e.\ $N_{\rc} = \frac{4}{3} \pi \rc^3 N/L^3$.}
  \label{tab:complexity-optimal-pts}
  \begin{tabular}{l|cccc|cccc}
    \toprule
    & \multicolumn{4}{c|}{\textbf{Uniform system}} &
    \multicolumn{4}{c}{\textbf{Clustered system}} \\
    & $D=0$ & $D=1$ & $D=2$ & $D=3$ &
    $D=0$ & $D=1$ & $D=2$ & $D=3$ \\
    \midrule
    Stokeslet, Rotlet & 2500 & 950 & 450 & 400 &
    800 & 300 & 170 & 120 \\
    Stresslet & 6000 & 1600 & 800 & 800 &
    3000 & 1200 & 370 & 300 \\
    \bottomrule
  \end{tabular}
\end{table}

\begin{figure}[htbp]
  \centering
  \includegraphics[width=0.5\textwidth]{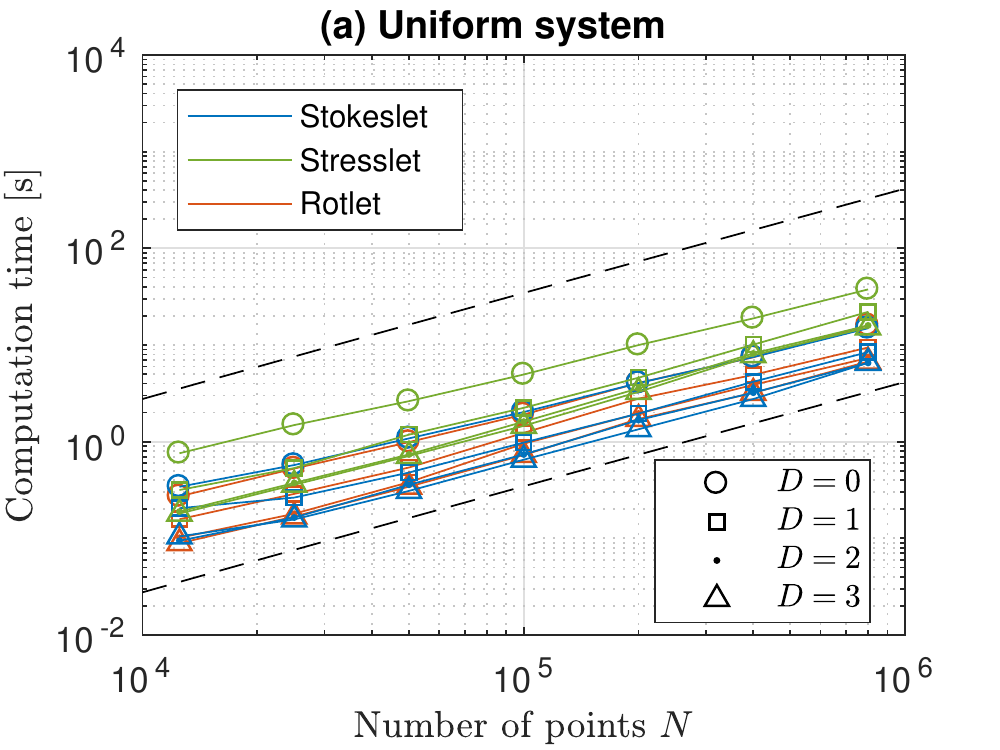}%
  \includegraphics[width=0.5\textwidth]{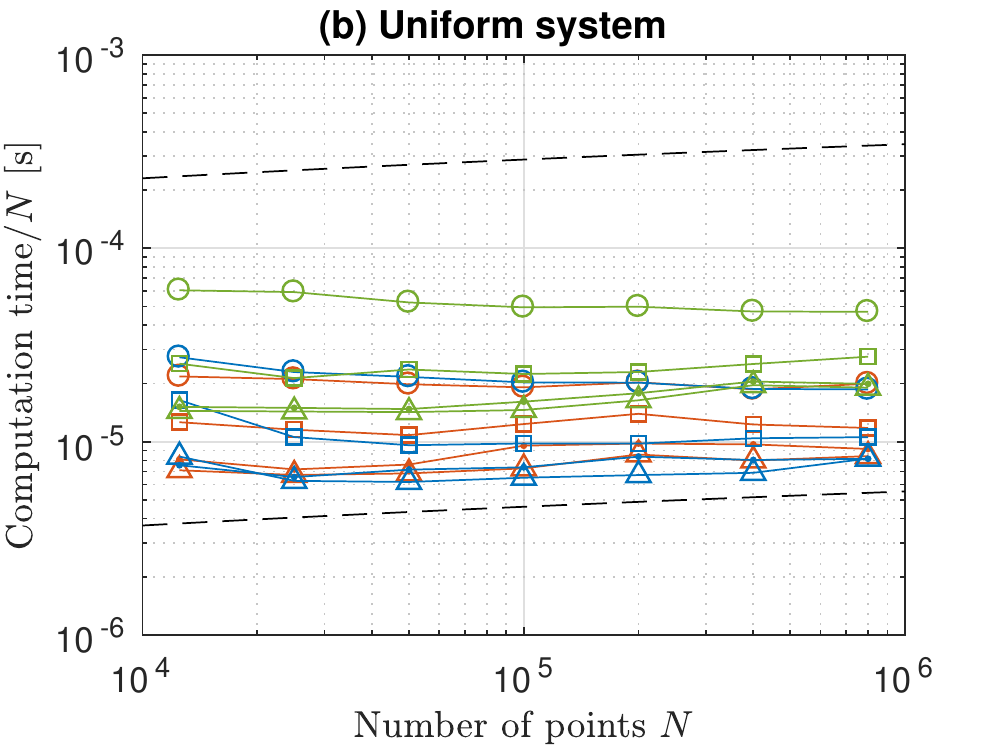}%
  \caption{Time complexity test, uniform system. The number of
  points $N$ varies ($N=2^j \times 12\,500$, $j=0,\ldots,6$)
  while $N_{\rc}$ is fixed to values given in
  Table~\ref{tab:complexity-optimal-pts}; other parameters are
  selected as stated in the main text. Total computation time
  excluding precomputation is shown in (a), and in (b) divided by
  $N$. Each data point is the average of 5 different point
  configurations. The dashed black lines indicate the expected
  complexity $O(N \log N)$.}
  \label{fig:complexity-uniform}
\end{figure}

\begin{figure}[htbp]
  \centering
  \includegraphics[width=0.5\textwidth]{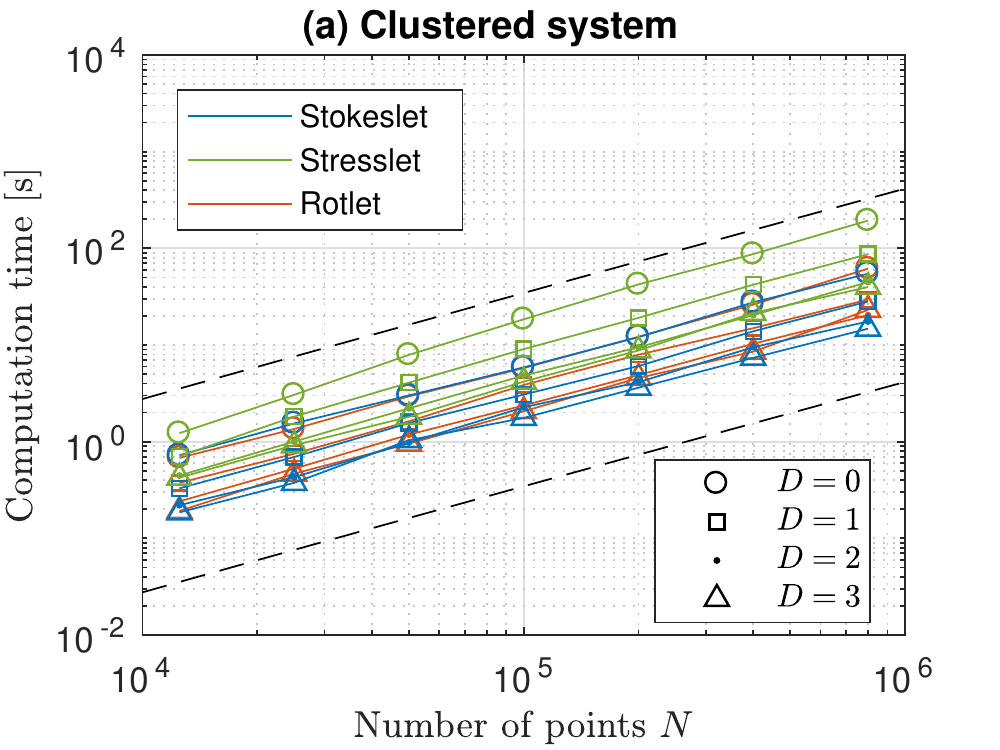}%
  \includegraphics[width=0.5\textwidth]{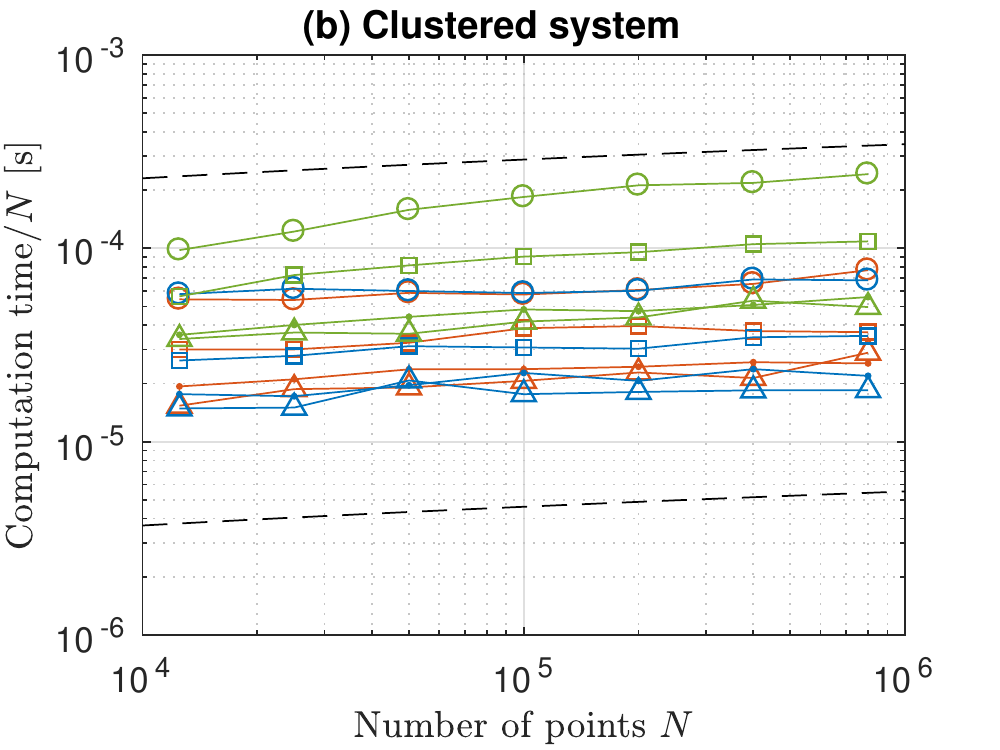}%
  \caption{Time complexity test, clustered system. For
  explanation, see Figure~\ref{fig:complexity-uniform}.}
  \label{fig:complexity-clustered}
\end{figure}

The computation time when increasing the number of points $N$,
while keeping $N_{\rc}$ fixed, is shown in
Figure~\ref{fig:complexity-uniform} for the uniform system, and
in Figure~\ref{fig:complexity-clustered} for the clustered
system. It can be noted that the computation time follows the
expected scaling $O(N \log N)$ for both systems, and that the
computation time is a bit larger for the clustered system than
for the uniform one. The uniform particle distribution is the
optimal one for the SE method, since there is no spatial
adaptivity in the Fourier-space part method.

As a second test, we fix $N = 10^5$ and measure the computation
time when varying the error tolerance $\abstol$ between $10^{-2}$
and $10^{-14}$. In this test, we use a fixed decomposition
parameter $\xi$, corresponding to the value selected (from
$N_{\rc}$) in the time complexity test for $N=10^5$ and
$\abstol=10^{-8}$; the value of $\xi$ is given per kernel,
periodicity and particle system in Table~\ref{tab:runtime-optimal-xi}.
Keeping $\xi$ fixed as the tolerance $\abstol$ is varied allows the run
time of both the real-space and Fourier-space parts to be
adjusted (through $\rc$ and $M$), so that the computational
balance between the two parts is kept approximately the same.
Parameters are set according to the parameter selection
procedure, as in the complexity test.

\begin{figure}[htp]
  \centering
  \includegraphics[width=0.33\textwidth]{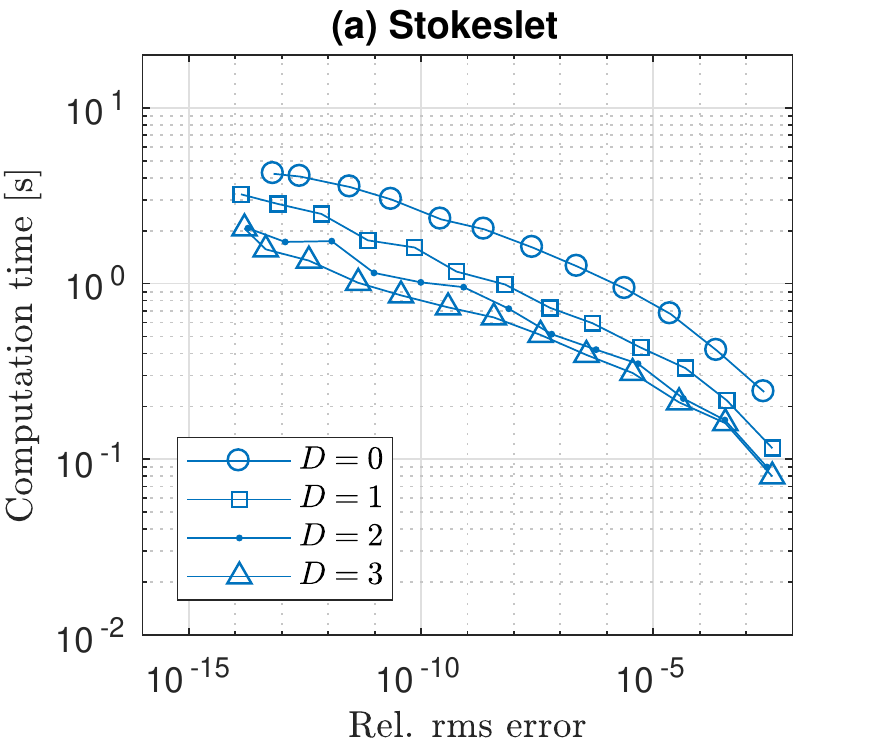}%
  \includegraphics[width=0.33\textwidth]{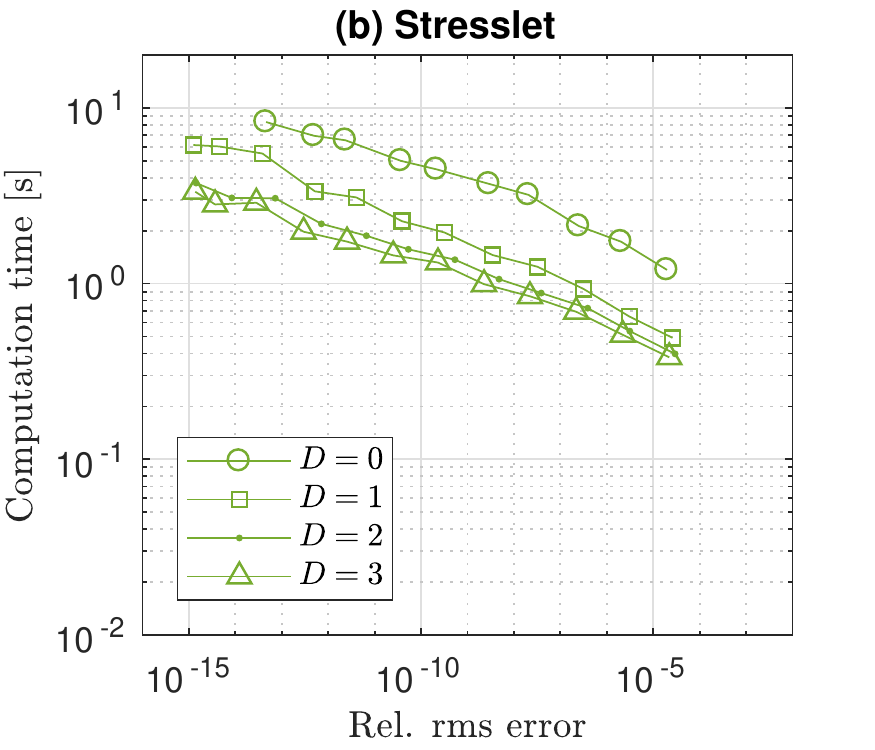}%
  \includegraphics[width=0.33\textwidth]{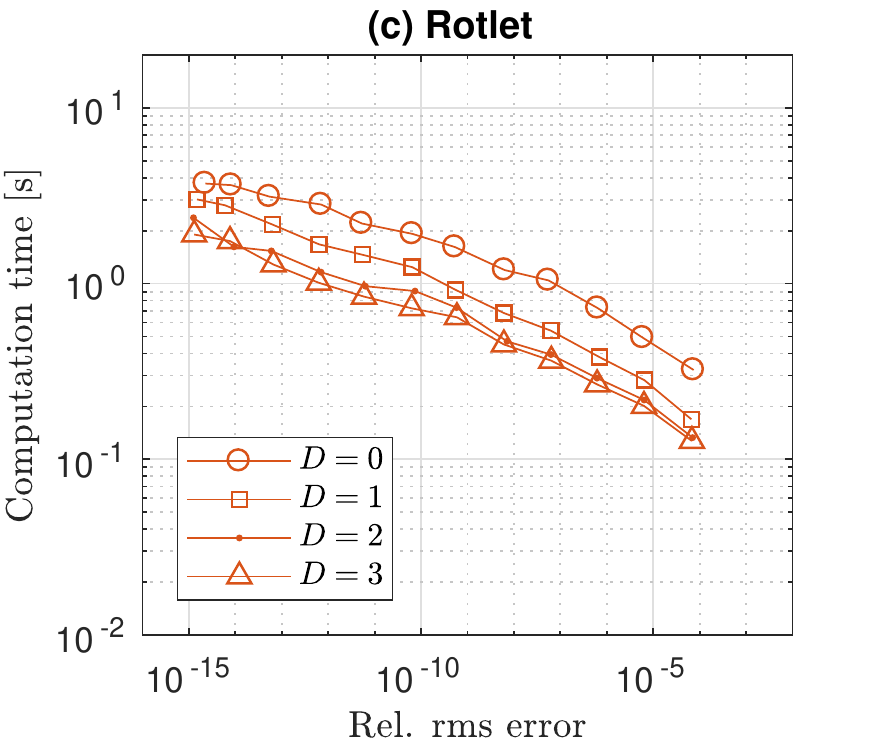}%
  \caption{Computation time versus error, uniform system.
  Total computation time excluding precomputation is plotted
  versus the actual relative rms error (of the full potential),
  for different kernels and periodicities~$D$.
  The value of $\xi$ for each case is given by
  Table~\ref{tab:runtime-optimal-xi}; furthermore, $N=10^5$,
  $L=1$, $Q=1$, while all other parameters are set from the
  absolute error tolerance $\abstol = 10^{-2}, 10^{-3}, \ldots,
  10^{-p}$ (stopping at $p \leq 14$ where the error starts to
  flatten out).
  Each data point is the average of 10 different point configurations.
  The absolute and relative errors are related by the rms value
  of the full potential, which is around 3 for the stokeslet, 500
  for the stresslet, and 200 for the rotlet. The reference
  potential is computed using the same method with
  $\abstol=10^{-16}$.}
  \label{fig:runtime}
\end{figure}

The result is shown in
Figure~\ref{fig:runtime} for the uniform system. (The clustered
system has a very similar behaviour, but with a factor 3 to 4
larger computational times, as already seen in the complexity
test. It is therefore omitted here.)
For each kernel, it can be seen that the computation time for
$D=2$ is very similar to $D=3$, while the time for $D=1$ is
at most twice as large; furthermore, the time for $D=0$ is around
three times as large as $D=3$. It can also be seen
that the stresslet requires more time than the other kernels
(up to twice as much) for the same relative error; this is
because the stresslet source strengths have more components, which for
example leads to more expensive gridding and more FFTs.

\begin{table}[htp]
  \centering
  \caption{Decomposition parameter $\xi$ used in the computation
  time versus error test. Computed from the values in
  Table~\ref{tab:complexity-optimal-pts} for $N=10^5$ and
  $\abstol=10^{-8}$.}
  \label{tab:runtime-optimal-xi}
  \begin{tabular}{l|cccc|cccc}
    \toprule
    & \multicolumn{4}{c|}{\textbf{Uniform system}} &
    \multicolumn{4}{c}{\textbf{Clustered system}} \\
    & $D=0$ & $D=1$ & $D=2$ & $D=3$ &
    $D=0$ & $D=1$ & $D=2$ & $D=3$ \\
    \midrule
    Stokeslet & 23.5579 & 32.3804 & 41.3955 & 43.0296 &
    34.2623 & 47.2969 & 57.0041 & 63.9175 \\
    Stresslet & 19.9813 & 31.1944 & 39.4021 & 39.4021 &
    25.2392 & 34.3700 & 51.0933 & 54.8344 \\
    Rotlet & 24.5091 & 33.9754 & 43.7207 & 45.4936 &
    36.0043 & 50.1320 & 60.7237 & 68.2974 \\
    \bottomrule
  \end{tabular}
\end{table}

\subsection{Window function comparison}
\label{sec:results-window}

Finally, we compare the PKB window function used in this paper
with the truncated Gaussian window traditionally used in the SE
method. The truncated Gaussian (TG) is given by (cf.\
\eqref{eq:KB-window})
\begin{equation}
  \label{eq:TG-window}
  \window_{0,\text{TG}}(r) =
  \begin{cases}
    \displaystyle
    \E^{-\alpha(r/a_\window)^2}, & \lvert r \rvert \leq a_\window, \\
    0, & \lvert r \rvert > a_\window,
  \end{cases}
\end{equation}
where $a_\window = hP/2$, and $\alpha$ is a shape parameter which
has the value $\alpha = 0.91 (\pi/2) P$. In the scaling step, the
Fourier transform of the window function is needed.
Rather than using the exact Fourier transform
$\wh{\window}_{0,\text{TG}}$ of the truncated Gaussian, which
involves the error function $\erf(\cdot)$, we prefer to use the
Fourier transform of the untruncated Gaussian
$\window_{0,\text{UG}}(r) = \E^{-\alpha(r/a_\window)^2}$, i.e.
\begin{equation}
  \wh{\window}_{0,\text{UG}}(k)
  =
  \sqrt{\frac{\pi}{\alpha}} a_\window \E^{-k^2 a_w^2/(4\alpha)}.
\end{equation}
The reason is that $\wh{\window}_{0,\text{UG}}$ is cheaper to
evaluate than $\wh{\window}_{0,\text{TG}}$, and the difference
between them is guaranteed to be around the selected error level
due to how parameters are selected. Thus, we use
$\window_{0,\text{TG}}$ in the gridding and gathering steps, but
$\wh{\window}_{0,\text{UG}}$ in the scaling step.

\begin{figure}[htp]
  \centering
  \includegraphics[width=0.5\textwidth]{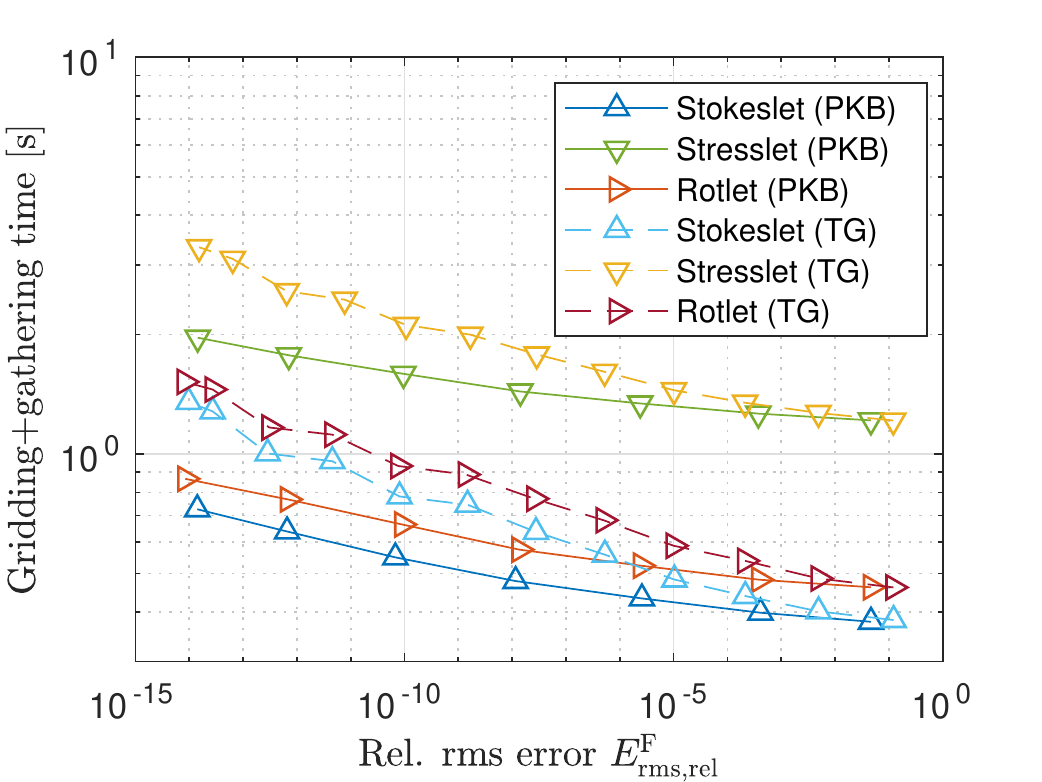}%
  \caption{Window function comparison, uniform system.
  Computation time of the gridding and gathering steps is shown
  versus the actual relative rms error of the Fourier-space part
  of the potential. Here, $N=10^5$, $Q=1$, $L=1$, and $\xi$ is
  given by Table~\ref{tab:runtime-optimal-xi}. Only the fully
  periodic ($D=3$) case is included here, but gridding and
  gathering times do not depend strongly on periodicity. The grid
  size is $M=196$ for the stokeslet, $M=192$ for the stresslet,
  and $M=208$ for the rotlet, which corresponds to a negligible
  absolute truncation error around $10^{-20}$. The data points
  are obtained by letting the window size $P$ vary from 2 to 14
  for the PKB window, and from 2 to 24 for the TG window (right
  to left), computing the error and the time for each $P$. Each data point
  is the average of 10 different point configurations. The rms
  value of the Fourier-space part of the potential is around 2
  for the stokeslet, 50 for the stresslet, and 20 for the rotlet.
  The reference potential is computed using the same method, with
  the PKB window and $P=20$.}
  \label{fig:window-comparison}
\end{figure}

To compare the PKB and TG windows, we consider for each kernel a
system of $N=10^5$ uniformly random sources, with $Q=1$, in a
cubic box of side length $L=1$. The parameter $\xi$ is selected
as in Table~\ref{tab:runtime-optimal-xi} (uniform system). The
parts of the SE method that depend most strongly on the window
function are the gridding and gathering steps, which have run
time $O(P^3)$ with respect to the window size $P$.
Ref.~\cite{Shamshirgar2021} showed, for the harmonic kernel, that $P$
can be reduced about 40~\% for the PKB window compared to the TG
window, for the same error. Thus, the gridding and gathering time
will be smaller for the PKB window than for the TG window. In
Figure~\ref{fig:window-comparison} we show that this is indeed
true also for the three kernels that we consider here: The
gridding and gathering time for the PKB window is significantly
reduced compared to the TG window, since the window size $P$ can
be selected smaller. This effect becomes more pronounced for
stricter error tolerances. The stokeslet and rotlet have
identical gridding and gathering steps; the only reason that the
run time is larger for the rotlet than for the stokeslet in
Figure~\ref{fig:window-comparison} is that $\xi$, and therefore
$M$, is selected larger for the rotlet, in accordance with
Table~\ref{tab:runtime-optimal-xi}.

Here, we have focused on the run time of the gridding and gathering
steps of the SE method. The FFT and IFFT steps do not depend on
the window function, while the scaling step has only a weak
dependence on the window (since the Fourier transform
$\wh{\window}_0$ must be evaluated). Thus, the PKB window will
typically lead to smaller computation times than the TG window,
also when the total computation time of the method is considered.
The percentual total time gain of using the PKB window will
depend on how much time is spent in the gridding and gathering
steps compared to the other steps, which depends on different
factors such as the number of points $N$, decomposition
parameter $\xi$, periodicity $D$, and the error tolerance
$\abstol$ (which sets $M$ and $P$). In the example considered
here, the gridding+gathering time makes up from 40~\% (for $P=2$)
to 70~\% (for $P=24$) of the total Fourier-space part run time.

We would like to point out that the SE method has been shown to
be competitive with other fast summation methods, such as the
fast multipole method (FMM); for instance, this has been done for
the harmonic kernel by \cite{Shamshirgar2021}, and for the
stokeslet by \cite{afKlinteberg2017}. In the latter case, the
comparison between SE and FMM was done in free space ($D=0$),
which is the most expensive periodic case for the SE method, but
the cheapest for the FMM. The comparison was done with the TG
window function in the SE method; as we have shown here, the PKB
window can further reduce the run time of the SE method.

\FloatBarrier
\section{Conclusions}
\label{sec:conclusions}

In this paper, we have presented the Spectral Ewald (SE) method,
a fast Ewald summation method, for three kernels (stokeslet,
stresslet, rotlet) of three-dimensional Stokes flow. The SE
method computes the Fourier-space part of the periodic potential;
it is based on the fast Fourier transform (FFT) and scales like $O(N
\log N)$ for $N$ particles in the primary cell. Arbitrary
periodicity ($D=3,2,1,0$ periodic directions) is supported, all
within the same framework. This paper marks the completion and
unification of the SE method for Stokes flow, uniting all kernels
and periodic cases, thus allowing to perform efficient
simulations of Stokes flow with arbitrary periodicity.

The modified kernels used to circumvent the singular behaviour of
the original kernels form a crucial part of the method for
$D=0,1,2$. We have improved the convergence of the modified
stokeslet and stresslet in the $D=0$ case, and derived new
modified kernels for the $D=1,2$ cases. With the new modified
kernels, the truncation error is independent of periodicity,
i.e.\ there is no penalty to the truncation error in the
$D=0,1,2$ cases compared to the $D=3$ case.
New improved truncation error estimates have been derived for the
stokeslet and stresslet, valid in all periodic cases.
We have also derived analytical formulas that can be used for
validation for the $D=1,2$ cases.

Our version of the SE method uses the polynomial Kaiser--Bessel
(PKB) window function introduced by \cite{Shamshirgar2021}, which
reduces the support, and thus the computational time, needed to
reach a given error compared with the truncated Gaussian (TG) window
traditionally used in the SE method. We have showed numerically
that the PKB window is indeed faster. The polynomial approximation
of the window is selected such that it introduces no further
error compared to the exact Kaiser--Bessel window. In the $D=1,2$
cases, an adaptive Fourier transform (AFT), first introduced by
\cite{Shamshirgar2017}, is used to reduce the cost of the FFTs.
In the $D=0$ case, a precomputation step is used to reduce the
computational time.

An automated procedure for selecting the parameters of the method
has been presented and tested. The error is within one order of
magnitude of the tolerance in the vast majority of cases. In the
numerical results, we have furthermore showed that the pointwise
error is well behaved, and that the computational time scales as
expected. The method is fastest in the $D=3$ case; the
computational time for $D=2$ is only slightly larger, while for
$D=1$ it is at most twice as large, and for $D=0$ around three
times as large as in the $D=3$ case.
The stresslet has more components than the other kernels and
therefore requires around twice the computational time, due to
more costly gridding and FFTs.

The SE method will be most efficient when applied to particles
which are uniformly distributed in a primary cell with low aspect
ratio (close to cubic). The computational time will increase if
points are very unevenly distributed, or if the cell has a high
aspect ratio. While the parameter selection procedure in
principle works for any rectangular cuboid cell, it has
been mainly tested for cells with low-to-moderate aspect ratios.

A drawback of the current parameter selection procedure is that
only an absolute error tolerance, not a relative tolerance, can
be given for the full potential (and real-space part). (The
relative error can of course be computed a posteriori, but not a
priori.) To be able to specify a relative tolerance for the full
potential (and real-space part), one would need an estimate of
its rms value. We have chosen not to construct such an estimate in
this paper, since it would necessarily have to incorporate the
distances from evaluation points to source points, making it more
complicated than the Fourier-space part estimate that we constructed.

Previous iterations of the SE method have been successfully
applied to simulations of Stokes flow based on boundary integral
equations, in the $D=3,0$ cases. In the future, we will apply
also the current version of the SE method to such simulations,
with arbitrary periodicity.
We believe that the techniques described in this paper would also
be applicable to other kernels, such as the
Rotne--Prager--Yamakawa tensor, assuming that a modified kernel
can be derived (needed for reduced periodicity).
Another possible future undertaking would be to examine the
so-called ``error pollution'' and aliasing errors further and try
to find a more refined rule, based on analytical
formulas, to adjust the uniform grid spacing $h$ and window size $P$.

The unified Spectral Ewald package with arbitrary periodicity
includes the three kernels of Stokes flow as well as the harmonic
kernel, and will be made available on GitHub \citep{SE_unified}.

\section*{Acknowledgements}

JB thanks Dr.\ Davood Saffar Shamshirgar for helpful discussions
on the implementation of the SE method during the initial stage
of this work. The authors gratefully acknowledge support from the
Swedish Research Council under Grant No.\ 2019-05206.

\appendix
\section{Analytical evaluation of Fourier integrals in doubly
periodic case}
\label{app:2p-integrals}

\subsection{Nonsingular case}
\label{app:2p-integrals-nonzero}

The goal here is to evaluate the integral
\begin{equation}
  \bmat{Q}^{2\per}(k_1,k_2,r_3;\xi)
  :=
  \frac{1}{2\pi} \int_{\mathbb{R}}
  \wh{\kernel}{}^\fp(k_1,k_2,\kappa_3)
  \E^{\I \kappa_3 r_3}
  \, \mathrm{d}\kappa_3,
  \qquad
  (k_1,k_2) \neq (0,0),
\end{equation}
that appears in \eqref{eq:2p-evaluated-integrals-nonzero},
obtained by combining \eqref{eq:fourier-potential-Dp} and
\eqref{eq:inverse-mixed-ft} for $D=2$. Here, the kernel $\kernel$
is the stokeslet, rotlet or stresslet. As noted below
\eqref{eq:diffop-rel-fourier-1}--\eqref{eq:diffop-rel-fourier-3},
we can write $\wh{\kernel}{}^\fp = \wh{\Kopt} \wh{A} \wh{\gamma}$,
where $A$ is either the harmonic kernel (for the rotlet) or
biharmonic kernel (for the stokeslet and stresslet). To be able
to convert $\wh{\Kopt}$ to $\Kopt$, we introduce a new quantity
\begin{equation}
  \label{eq:Qtilde-2p-def}
  \widetilde{\bmat{Q}}^{2\per}(k_1,k_2,r_1,r_2,r_3;\xi) :=
  \bmat{Q}^{2\per}(k_1,k_2,r_3;\xi) \E^{\I k_1 r_1} \E^{\I k_2 r_2}
  =
  \frac{1}{2\pi} \int_{\mathbb{R}}
  \wh{\kernel}{}^\fp(k_1,k_2,\kappa_3)
  \E^{\I (k_1,k_2,\kappa_3) \cdot (r_1,r_2,r_3)}
  \, \mathrm{d}\kappa_3.
\end{equation}
Inserting the relation $\wh{\kernel}{}^\fp = \wh{\Kopt} \wh{A}
\wh{\gamma}$ above, we find that
\begin{equation}
  \widetilde{\bmat{Q}}^{2\per}(k_1,k_2,r_1,r_2,r_3;\xi) =
  \Kopt \left(
  \E^{\I k_1 r_1} \E^{\I k_2 r_2}
  \frac{1}{2\pi} \int_{\mathbb{R}}
  \wh{A}(k_1,k_2,\kappa_3) \wh{\gamma}(k_1,k_2,\kappa_3;\xi)
  \E^{\I \kappa_3 r_3}
  \, \mathrm{d}\kappa_3
  \right).
\end{equation}
This shows that in order to compute $\bmat{Q}^{2\per}$, we can
first compute the integral
\begin{equation}
  \label{eq:QA-2p-def}
  Q^{A,2\per}(k_1,k_2,r_3;\xi) :=
  \frac{1}{2\pi} \int_{\mathbb{R}}
  \wh{A}(k_1,k_2,\kappa_3) \wh{\gamma}(k_1,k_2,\kappa_3;\xi)
  \E^{\I \kappa_3 r_3}
  \, \mathrm{d}\kappa_3,
\end{equation}
and then apply the relation
\begin{equation}
  \label{eq:2p-nonzero-mode-relation}
  \bmat{Q}^{2\per}(k_1,k_2,r_3;\xi)
  =
  \E^{-\I k_1 r_1} \E^{-\I k_2 r_2}
  \Kopt \Big(\E^{\I k_1 r_1} \E^{\I k_2 r_2} Q^{A,2\per}(k_1,k_2,r_3;\xi)\Big),
\end{equation}
obtained from \eqref{eq:Qtilde-2p-def}--\eqref{eq:QA-2p-def}.
To simplify the expressions in the following, we first introduce
some auxiliary variables, namely
\begin{equation}
  \label{eq:2p-aux}
  \begin{array}{c}
    \displaystyle
    \alpha = \sqrt{k_1^2+k_2^2}, \qquad
    \widetilde{k}_1 = \frac{k_1}{\alpha}, \qquad
    \widetilde{k}_2 = \frac{k_2}{\alpha},
    \\[10pt]
    \displaystyle
    \theta^{+} = \E^{\alpha r_3} \erfc\left(\frac{\alpha}{2\xi} + \xi r_3 \right), \qquad
    \theta^{-} = \E^{-\alpha r_3} \erfc\left(\frac{\alpha}{2\xi} - \xi r_3 \right), \qquad
    \beta^{+} = \theta^{+} + \theta^{-}, \qquad
    \beta^{-} = \theta^{+} - \theta^{-},
    \\[10pt]
    \displaystyle
    \lambda = \frac{2}{\sqrt{\pi}\xi} \E^{-\alpha^2/(2\xi)^2}\E^{-\xi^2 r_3^2}, \qquad
    \Lambda = \lambda + \frac{1}{\alpha} \beta^{+} - r_3 \beta^{-}, \qquad
    \Psi = -\frac{1}{\alpha} \beta^{+} - r_3 \beta^{-}.
  \end{array}
\end{equation}

For the harmonic kernel $H(\bvec{r}) = 1/\lvert\bvec{r}\rvert \ftpair
\wh{H}(\bvec{k}) = 4\pi/\lvert\bvec{k}\rvert^2$ and the
Ewald screening function $\gE$ \eqref{eq:ewald-screening},
the integral $Q^{H,2\per}$, defined by \eqref{eq:QA-2p-def}, was
stated by \cite{Tornberg2016}. The value of the integral is found
from \cite[3.954 (2), p.~504]{Gradshteyn2007},
\begin{equation}
  \label{eq:2p-QH-integral}
  Q^{H,2\per}(k_1,k_2,r_3;\xi) =
  \frac{1}{2\pi} \int_{\mathbb{R}} \frac{4\pi}{\alpha^2+\kappa_3^2}
  \E^{-(\alpha^2+\kappa_3^2)/(2\xi)^2} \E^{\I \kappa_3 r_3}
  \, \mathrm{d}\kappa_3
  = \frac{\pi}{\alpha} \beta^{+},
\end{equation}
with $\alpha$ and $\beta^{+}$ as in \eqref{eq:2p-aux}.
Applying the relation \eqref{eq:2p-nonzero-mode-relation}
for the rotlet, with $A=H$ and the operator $\Kopt^\rotlet$ given
by \eqref{eq:diffop-rel-2}, we find that the tensor
$\bmat{Q}^{\rotlet,2\per}$ is antisymmetric
and given by
\begin{equation}
  \label{eq:2p-rotlet-Q}
  \bmat{Q}^{\rotlet,2\per}(k_1,k_2,r_3;\xi)
  =
  \pi
  \begin{bmatrix}
    0 & -\beta^{-} & \I \widetilde{k}_2 \beta^{+} \\
    \beta^{-} & 0 & -\I \widetilde{k}_1 \beta^{+} \\
    -\I \widetilde{k}_2 \beta^{+} & \I \widetilde{k}_1 \beta^{+} & 0 \\
  \end{bmatrix}.
\end{equation}
Here, we have used that $\partial \beta^{+} / \partial r_3 = \alpha \beta^{-}$.

For the biharmonic kernel $\biharmonic(\bvec{r}) = \lvert \bvec{r} \rvert
\ftpair \wh{\biharmonic}(\bvec{k}) = -8\pi/\lvert \bvec{k} \rvert^4$ and
the Hasimoto screening function $\gH$ \eqref{eq:hasimoto-screening},
the integral $Q^{\biharmonic,2\per}$ is given by
\begin{align}
  \notag
  Q^{\biharmonic,2\per}(k_1,k_2,r_3;\xi) &=
  -\frac{1}{2\pi} \int_{\mathbb{R}} \frac{8\pi}{(\alpha^2+\kappa_3^2)^2}
  \E^{-(\alpha^2+\kappa_3^2)/(2\xi)^2}
  \left( 1 + \frac{\alpha^2+\kappa_3^2}{(2\xi)^2} \right)
  \E^{\I \kappa_3 r_3}
  \, \mathrm{d}\kappa_3
  \\
  \label{eq:2p-QB-integral}
  &=
  -4 \int_{\mathbb{R}}
  \left( \frac{1}{(\alpha^2+\kappa_3^2)^2}
  + \frac{1}{4\xi^2}\frac{1}{\alpha^2+\kappa_3^2} \right)
  \E^{-(\alpha^2+\kappa_3^2)/(2\xi)^2}
  \E^{\I \kappa_3 r_3}
  \, \mathrm{d}\kappa_3
  .
\end{align}
From the relation
\begin{equation}
  \label{eq:2p-B-H-relation}
  \left( \frac{1}{(\alpha^2+\kappa_3^2)^2}
  + \frac{1}{4\xi^2}\frac{1}{\alpha^2+\kappa_3^2} \right)
  \E^{-(\alpha^2+\kappa_3^2)/(2\xi)^2}
  =
  -\frac{1}{2\alpha} \frac{\partial}{\partial\alpha}
  \left(
  \frac{1}{\alpha^2+\kappa_3^2} \E^{-(\alpha^2+\kappa_3^2)/(2\xi)^2}
  \right)
\end{equation}
between the integrands of \eqref{eq:2p-QH-integral} and
\eqref{eq:2p-QB-integral}, one obtains the relation
\begin{equation}
  Q^{\biharmonic,2\per}(k_1,k_2,r_3;\xi) = \frac{1}{\alpha}
  \frac{\partial}{\partial \alpha} Q^{H,2\per}(k_1,k_2,r_3;\xi).
\end{equation}
Using the result \eqref{eq:2p-QH-integral} as well as
$\partial \beta^{+} / \partial\alpha = r_3\beta^{-}-\lambda$, we
obtain
\begin{equation}
  \label{eq:2p-QB-result}
  Q^{\biharmonic,2\per}(k_1,k_2,r_3;\xi) = -\frac{\pi}{\alpha^2} \Lambda,
\end{equation}
with $\Lambda$ as in \eqref{eq:2p-aux}.

Applying the relation \eqref{eq:2p-nonzero-mode-relation} for the
stokeslet, with $A=\biharmonic$ and the operator
$\Kopt^\stokeslet$ given by \eqref{eq:diffop-rel-1}, one finds
that the tensor $\bmat{Q}^{\stokeslet,2\per}$ is symmetric
and given by
\begin{equation}
  \label{eq:2p-stokeslet-Q}
  \bmat{Q}^{\stokeslet,2\per}(k_1,k_2,r_3;\xi)
  =
  \pi
  \begin{bmatrix}
    \widetilde{k}_2^2 \Lambda - \Psi & -\widetilde{k}_1 \widetilde{k}_2 \Lambda & -\I \widetilde{k}_1 r_3 \beta^{+} \\
    -\widetilde{k}_1 \widetilde{k}_2 \Lambda & \widetilde{k}_1^2 \Lambda - \Psi & -\I \widetilde{k}_2 r_3 \beta^{+} \\
    -\I \widetilde{k}_1 r_3 \beta^{+} & -\I \widetilde{k}_2 r_3 \beta^{+} & \Lambda \\
  \end{bmatrix},
\end{equation}
with $\Psi$ as in \eqref{eq:2p-aux}. Here, we have used that
$\partial \beta^{-}/\partial r_3 = \alpha \beta^{+} - 2\xi^2 \lambda$
and $\partial \Lambda / \partial r_3 = -r_3 \alpha \beta^{+}$.
The result \eqref{eq:2p-stokeslet-Q} has also been derived by
\cite{Lindbo2011b} using a somewhat different method and notation.
For the stresslet, with $A=\biharmonic$
and $\Kopt^\stresslet$ given by \eqref{eq:diffop-rel-3}, one
similarly finds that $\bmat{Q}^{\stresslet,2\per}$ is symmetric
with entries given by
\begin{equation}
  \label{eq:2p-stresslet-Q}
  \begin{array}{r@{\:}c@{\:}l}
    Q^{\stresslet,2\per}_{111}(k_1,k_2,r_3;\xi)
      &=&
      \displaystyle
      \I \pi \widetilde{k}_1 \alpha \left(
        (2 \widetilde{k}_2^2 + 1) \Lambda - 3 \Psi
      \right),
    \\[5pt]
    Q^{\stresslet,2\per}_{112}(k_1,k_2,r_3;\xi)
      &=&
      \displaystyle
      \I \pi \widetilde{k}_2 \alpha \left(
        (2 \widetilde{k}_2^2 - 1) \Lambda - \Psi
      \right),
    \\[5pt]
    Q^{\stresslet,2\per}_{113}(k_1,k_2,r_3;\xi)
      &=&
      \displaystyle
      - 2\pi \left(
        \xi^2 r_3 \lambda - 2 \widetilde{k}_1^2 r_3 \alpha \beta^{+} - 2 \beta^{-}
      \right),
    \\[5pt]
    Q^{\stresslet,2\per}_{122}(k_1,k_2,r_3;\xi)
      &=&
      \displaystyle
      \I \pi \widetilde{k}_1 \alpha \left(
        (2 \widetilde{k}_1^2 - 1) \Lambda - \Psi
      \right),
    \\[5pt]
    Q^{\stresslet,2\per}_{123}(k_1,k_2,r_3;\xi)
      &=&
      \displaystyle
      2 \pi \widetilde{k}_1 \widetilde{k}_2 r_3 \alpha \beta^{+},
    \\[5pt]
    Q^{\stresslet,2\per}_{133}(k_1,k_2,r_3;\xi)
      &=&
      \displaystyle
      \I \pi \widetilde{k}_1 \alpha (\Lambda + \Psi),
    \\[5pt]
    Q^{\stresslet,2\per}_{222}(k_1,k_2,r_3;\xi)
      &=&
      \displaystyle
      \I \pi \widetilde{k}_2 \alpha \left(
        (2 \widetilde{k}_1^2 + 1) \Lambda - 3 \Psi
      \right),
    \\[5pt]
    Q^{\stresslet,2\per}_{223}(k_1,k_2,r_3;\xi)
      &=&
      \displaystyle
      - 2 \pi
      \left(
        \xi^2 r_3 \lambda - 2 \widetilde{k}_2^2 r_3 \alpha \beta^{+} - 2 \beta^{-}
      \right),
    \\[5pt]
    Q^{\stresslet,2\per}_{233}(k_1,k_2,r_3;\xi)
      &=&
      \displaystyle
      \I \pi \widetilde{k}_2 \alpha (\Lambda + \Psi),
    \\[5pt]
    Q^{\stresslet,2\per}_{333}(k_1,k_2,r_3;\xi)
      &=&
      \displaystyle
      -2\pi \left(
        \xi^2 r_3 \lambda + 2 r_3 \alpha \beta^{+} - 2\beta^{-}
      \right),
  \end{array}
\end{equation}
again with auxiliary variables as in \eqref{eq:2p-aux}.

\subsection{Singular case}
\label{app:2p-integrals-k0}

We now want to compute the integral
\begin{equation}
  \bmat{Q}^{2\per,(0)}(r_3;\xi)
  :=
  \frac{1}{2\pi} \int_{\mathbb{R}}
  \wh{\kernel}{}^\fp(0,0,\kappa_3)
  \E^{\I \kappa_3 r_3}
  \, \mathrm{d}\kappa_3,
\end{equation}
that appears in \eqref{eq:2p-evaluated-integrals-k0}, obtained
from \eqref{eq:fourier-potential-Dp} and
\eqref{eq:inverse-mixed-ft} for $D=2$. Note that the integrand is
now singular at $\kappa_3=0$. The integral does not exist in the
Lebesgue sense, but we can interpret $\bmat{Q}^{2\per,(0)}$ as
the (one-dimensional) inverse Fourier transform of
$\wh{\kernel}{}^\fp(0,0,\cdot)$, in the distributional sense.
(An introduction to the distributional Fourier transform can be
found e.g.\ in \cite[chapter~8]{Vretblad2005}.) In the same way
as in \ref{app:2p-integrals-nonzero}, we apply the relation
$\wh{\kernel}{}^\fp = \wh{\Kopt} \wh{A} \wh{\gamma}$, noting that
in this case, $\wh{\Kopt}(0,0,\kappa_3) \E^{\I \kappa_3 r_3} =
\Kopt \E^{\I \kappa_3 r_3}$. The result is that
\begin{equation}
  \label{eq:2p-zero-mode-relation}
  \bmat{Q}^{2\per,(0)}(r_3;\xi) =
  \Kopt Q^{A,2\per,(0)}(r_3;\xi),
\end{equation}
where
\begin{equation}
  Q^{A,2\per,(0)}(r_3;\xi) :=
  \frac{1}{2\pi} \int_{\mathbb{R}}
  \wh{A}(0,0,\kappa_3) \wh{\gamma}(0,0,\kappa_3;\xi)
  \E^{\I \kappa_3 r_3}
  \, \mathrm{d}\kappa_3
  ,
\end{equation}
where $A$ is the harmonic kernel (for the rotlet) or biharmonic
kernel (for the stokeslet and stresslet). Note that $\wh{A}
\wh{\gamma} = \wh{A}^\fp$, so that
\begin{equation}
  Q^{A,2\per,(0)}(r_3;\xi) = \fourier^{-1}\{\wh{A}^\fp(0,0,\cdot;\xi)\}(r_3;\xi)
  =: \pft{A}^\fp(0,0,r_3;\xi).
\end{equation}
One might try to compute $Q^{A,2\per,(0)}(r_3;\xi)$ as the limit of
$Q^{A,2\per}(k_1,k_2,r_3;\xi) =: \pft{A}^\fp(k_1,k_2,r_3;\xi)$ from
\ref{app:2p-integrals-nonzero} as $\alpha = \sqrt{k_1^2+k_2^2} \to 0$. However, this
limit exists neither for the harmonic nor for the biharmonic.
Instead, we follow \cite{Tornberg2016} and use the Ewald
decomposition $A = A^\rp + A^\fp$ to write
\begin{equation}
  \label{eq:2p-zero-mode-split}
  \pft{A}^\fp(0,0,r_3;\xi) = \pft{A}(0,0,r_3) - \pft{A}^\rp(0,0,r_3;\xi),
\end{equation}
where
\begin{equation}
  \label{eq:2p-zero-mode-R}
  \pft{A}^\rp(0,0,r_3;\xi) = \lim_{\alpha \to 0}\pft{A}^\rp(k_1,k_2,r_3;\xi)
  = \lim_{\alpha \to 0} \Big(
    \pft{A}(k_1,k_2,r_3) - \pft{A}^\fp(k_1,k_2,r_3;\xi)
  \Big),
\end{equation}
and this limit exists since the real-space part of the kernel has
such a rapid decay in real space. It remains to compute
$\pft{A}(k_1,k_2,r_3)$ and $\pft{A}(0,0,r_3)$, but this is typically
an easier problem since there is no screening function involved.

For the harmonic, the derivation was done by \cite{Tornberg2016}.
The result is
\begin{equation}
  Q^{\harmonic,2\per,(0)}(r_3;\xi) = -2\pi \left(
  r_3 \erf(\xi r_3)
  +
  \frac{\E^{-\xi^2 r_3^2}}{\sqrt{\pi}\xi}
  \right)
  .
\end{equation}
For the rotlet, applying \eqref{eq:2p-zero-mode-relation} with
$A=\harmonic$ and $\Kopt^\rotlet$ given by \eqref{eq:diffop-rel-2},
we find that the tensor $\bmat{Q}^{\rotlet,2\per,(0)}$ is
antisymmetric and given by
\begin{equation}
  \label{eq:2p-rotlet-Q0}
  \bmat{Q}^{\rotlet,2\per,(0)}(r_3;\xi)
  =
  2\pi \erf(\xi r_3)
  \begin{bmatrix}
    0 & 1 & 0 \\
    -1 & 0 & 0 \\
    0 & 0 & 0 \\
  \end{bmatrix}
  .
\end{equation}

We now turn to the biharmonic kernel ($A=\biharmonic$). Our first goal will be to
compute $\pft{\biharmonic}(0,0,r_3)$ and $\pft{\biharmonic}(k_1,k_2,r_3)$, which are
needed in \eqref{eq:2p-zero-mode-split}--\eqref{eq:2p-zero-mode-R}.
These are defined by taking the one-dimensional inverse Fourier
transform in the free direction of $\wh{\biharmonic}(k_1,k_2,\kappa_3) =
-8\pi/(\alpha^2+\kappa_3^2)^2$, where again $\alpha=\sqrt{k_1^2+k_2^2}$.
For $\alpha > 0$, the inverse Fourier transform exists in the
classical sense, and (using e.g.\ \citep[13.2, F41b, p.~319]{RaadeWestergren2004} and
the convolution theorem)
\begin{equation}
  \label{eq:2p-Bbar-nonzero}
  \pft{\biharmonic}(k_1,k_2,r_3) = \frac{1}{2\pi} \int_{-\infty}^\infty
  \left( - \frac{8\pi}{(\alpha^2 + \kappa_3^2)^2} \right)
  \E^{\I \kappa_3 r_3} \, \mathrm{d}\kappa_3
  = -\frac{2\pi}{\alpha^3} \E^{-\alpha \lvert r_3 \rvert}
  (1 + \alpha \lvert r_3 \rvert), \qquad (k_1,k_2) \neq (0,0).
\end{equation}
For $\alpha=0$, the inverse Fourier transform has to be
interpreted in the distributional sense, the result being (the
integral here is only symbolic) \citep[13.2, F20, p.~317]{RaadeWestergren2004}
\begin{equation}
  \label{eq:2p-Bbar-zero}
  \pft{\biharmonic}(0,0,r_3) = \frac{1}{2\pi} \int_{-\infty}^\infty
  \left( - \frac{8\pi}{\kappa_3^4} \right)
  \E^{\I \kappa_3 r_3} \, \mathrm{d}\kappa_3
  = -\frac{2\pi}{3} \lvert r_3 \rvert^3.
\end{equation}
(The results in \eqref{eq:2p-Bbar-nonzero} and
\eqref{eq:2p-Bbar-zero} may also be viewed as solutions to the
PDE $-(-\alpha^2+\nabla_3^2)^2 \pft{\biharmonic} = 8\pi \delta(r_3)$ for
nonzero and zero $\alpha$, respectively. We mention this since it
corresponds to the view taken by \cite{Tornberg2016} for the
harmonic kernel.)

We can now compute the limit in \eqref{eq:2p-zero-mode-R} to find
$\pft{\biharmonic}^\rp(0,0,r_3;\xi)$. With
$\pft{\biharmonic}(k_1,k_2,r_3)$ taken from \eqref{eq:2p-Bbar-nonzero}
and $\pft{\biharmonic}^\fp(k_1,k_2,r_3;\xi)$ taken from \eqref{eq:2p-QB-result},
we find that we want to compute the limit
\begin{equation}
  \pft{\biharmonic}^\rp(0,0,r_3;\xi)
  =
  \lim_{\alpha\to0} \left(
    -\frac{2\pi \E^{-\alpha \lvert r_3 \rvert}
    (1 + \alpha \lvert r_3 \rvert)
    - \pi \alpha \Lambda}{\alpha^3}
  \right)
  .
\end{equation}
Both the numerator and denominator have zero as their limits, so
we apply L'H\^{o}pital's rule. Differentiating the numerator with
respect to $\alpha$, we get (using that $\partial \Lambda /
\partial \alpha = -\Lambda/\alpha - r_3^2 \beta^{+} - \alpha\lambda/(2 \xi^2)$)
\begin{equation}
  \pi \alpha \left( -2 r_3^2 \E^{-\alpha \lvert r_3 \rvert} + r_3^2 \beta^{+}
  + \frac{\alpha \lambda}{2 \xi^2} \right),
\end{equation}
and the denominator becomes $3\alpha^2$. Thus, the limit becomes
\begin{equation}
  \pft{\biharmonic}^\rp(0,0,r_3;\xi)
  =
  -\frac{\pi}{3}
  \lim_{\alpha\to0}
    \frac{-2 r_3^2 \E^{-\alpha \lvert r_3 \rvert} + r_3^2 \beta^{+}
  + \alpha \lambda/(2 \xi^2)}{\alpha}
  .
\end{equation}
Again, both numerator and denominator have zero limits, so we
apply L'H\^{o}pital's rule once more. Differentiating the
numerator now yields (using $\partial \beta^{+}/\partial \alpha =
r_3 \beta^{-} - \lambda$ and $\partial \lambda / \partial \alpha =
-\alpha\lambda/(2\xi^2)$)
\begin{equation}
  \label{eq:2p-lhopital-numerator}
  2 \lvert r_3 \rvert^3 \E^{-\alpha \lvert r_3 \rvert}
  + r_3^2 (r_3 \beta^{-} - \lambda)
  + \frac{\lambda}{2 \xi^2}
  - \frac{\alpha^2 \lambda}{4 \xi^4}
  ,
\end{equation}
while the denominator becomes $1$. The limit of
\eqref{eq:2p-lhopital-numerator} as $\alpha\to0$ can be computed,
and we get
\begin{equation}
  \label{eq:2p-Bbar-R-zero}
  \pft{\biharmonic}^\rp(0,0,r_3;\xi)
  =
  -\frac{2\pi}{3}
  \left(
  \lvert r_3 \rvert^3
  - r_3^3 \erf(\xi r_3)
  + \frac{\E^{-\xi^2 r_3^2}}{\sqrt{\pi}\xi} \left(
    \frac{1}{2\xi^2} - r_3^2
  \right)
  \right)
  .
\end{equation}
Finally, we use \eqref{eq:2p-zero-mode-split}, i.e.\ we subtract
$\pft{\biharmonic}^\rp(0,0,r_3;\xi)$ as given by
\eqref{eq:2p-Bbar-R-zero} from $\pft{\biharmonic}(0,0,r_3)$ as
given by \eqref{eq:2p-Bbar-zero}, which results in
\begin{equation}
  \label{eq:2p-Bbar-F-zero}
  Q^{\biharmonic,2\per,(0)}(r_3;\xi)
  =
  \pft{\biharmonic}^\fp(0,0,r_3;\xi)
  =
  -
  \frac{2\pi}{3}
  \left(
  r_3^3 \erf(\xi r_3)
  + \frac{\E^{-\xi^2 r_3^2}}{\sqrt{\pi}\xi} \left(
    r_3^2 - \frac{1}{2\xi^2}
  \right)
  \right)
  .
\end{equation}
This concludes the derivation for the biharmonic.

(In this case, i.e.\ $D=2$, it is also possible to compute
\eqref{eq:2p-Bbar-F-zero} directly by taking the one-dimensional
inverse Fourier transform in the $r_3$ direction of
$\wh{\biharmonic}(0,0,\kappa_3)\hgH(0,0,\kappa_3;\xi)$,
in the distributional sense (e.g.\ using Wolfram Mathematica~12). 
However, the derivation presented above is more straightforward
to apply in the $D=1$ case.)

We now apply \eqref{eq:2p-zero-mode-relation} for the stokeslet,
with $A=B$ and $\Kopt^\stokeslet$ given by
\eqref{eq:diffop-rel-1}. The resulting tensor
$\bmat{Q}^{\stokeslet,2\per,(0)}$ is symmetric
and given by
\begin{equation}
  \label{eq:2p-stokeslet-Q0}
  \bmat{Q}^{\stokeslet,2\per,(0)}(r_3;\xi)
  =
  -2\pi \left(
  2r_3 \erf(\xi r_3)
  +
  \frac{\E^{-\xi^2 r_3^2}}{\sqrt{\pi}\xi}
  \right)
  \begin{bmatrix}
    1 & 0 & 0 \\
    0 & 1 & 0 \\
    0 & 0 & 0 \\
  \end{bmatrix}
  .
\end{equation}
This result agrees with the expression for the $\bvec{k}^\per=\bvec{0}$ mode
derived by \cite{Lindbo2011b} using a different method.
For the stresslet, with $A=B$ and $\Kopt^\stresslet$ given by
\eqref{eq:diffop-rel-3}, we find that
$\bmat{Q}^{\stresslet,2\per,(0)}$ is given by
\begin{equation}
  \label{eq:2p-stresslet-Q0}
  \bmat{Q}^{\stresslet,2\per,(0)}(r_3;\xi)
  =
  -4\pi
  \left(
  \erf(\xi r_3)
  +
  \frac{\xi r_3}{\sqrt{\pi}} \E^{-\xi^2 r_3^2}
  \right)
  \bmat{C}^{\stresslet,2\per,(0)},
\end{equation}
where the constant tensor $\bmat{C}^{\stresslet,2\per,(0)}$ is
symmetric with entries given by
\begin{equation}
  C^{\stresslet,2\per,(0)}_{1lm}
  =
  \begin{bmatrix}
    0 & 0 & 1 \\
    0 & 0 & 0 \\
    1 & 0 & 0 \\
  \end{bmatrix}_{lm}
  ,
  \quad
  C^{\stresslet,2\per,(0)}_{2lm}
  =
  \begin{bmatrix}
    0 & 0 & 0 \\
    0 & 0 & 1 \\
    0 & 1 & 0 \\
  \end{bmatrix}_{lm}
  ,
  \quad
  C^{\stresslet,2\per,(0)}_{3lm}
  =
  \begin{bmatrix}
    1 & 0 & 0 \\
    0 & 1 & 0 \\
    0 & 0 & 1 \\
  \end{bmatrix}_{lm}
  .
\end{equation}

\section{Analytical evaluation of Fourier integrals in singly
periodic case}
\label{app:1p-integrals}

The structure of this section closely follows that of
\ref{app:2p-integrals} for the doubly periodic case.

\subsection{Nonsingular case}
\label{app:1p-integrals-nonzero}

We wish to evaluate the integral
\begin{equation}
  \bmat{Q}^{1\per}(k_1,r_2,r_3;\xi)
  :=
  \frac{1}{(2\pi)^2} \int_{\mathbb{R}^2}
  \wh{\kernel}{}^\fp(k_1,\kappa_2,\kappa_3)
  \E^{\I \kappa_2 r_2}
  \E^{\I \kappa_3 r_3}
  \, \mathrm{d}\kappa_2
  \, \mathrm{d}\kappa_3,
  \qquad
  k_1 \neq 0,
\end{equation}
that appears in \eqref{eq:1p-evaluated-integrals-nonzero}.
The kernel $\kernel$ may be the stokeslet, rotlet or stresslet.
Introducing the quantity
\begin{equation}
  \label{eq:Qtilde-1p-def}
  \widetilde{\bmat{Q}}^{1\per}(k_1,r_1,r_2,r_3;\xi) :=
  \bmat{Q}^{1\per}(k_1,r_2,r_3;\xi) \E^{\I k_1 r_1},
\end{equation}
and inserting $\wh{\kernel}{}^\fp = \wh{\Kopt} \wh{A}
\wh{\gamma}$, using that $\wh{\Kopt}\E^{\I \bvec{k} \cdot
\bvec{r}} = \Kopt\E^{\I \bvec{k} \cdot \bvec{r}}$, we find that
\begin{equation}
  \widetilde{\bmat{Q}}^{1\per}(k_1,r_1,r_2,r_3;\xi) =
  \Kopt \Big(
  \E^{\I k_1 r_1}
  Q^{A,1\per}(k_1,r_2,r_3;\xi)
  \Big),
\end{equation}
where
\begin{equation}
  \label{eq:QA-1p-def}
  Q^{A,1\per}(k_1,r_2,r_3;\xi) :=
  \frac{1}{(2\pi)^2} \int_{\mathbb{R}^2}
  \wh{A}(k_1,\kappa_2,\kappa_3) \wh{\gamma}(k_1,\kappa_2,\kappa_3;\xi)
  \E^{\I \kappa_2 r_2} \E^{\I \kappa_3 r_3}
  \, \mathrm{d}\kappa_2 \, \mathrm{d}\kappa_3.
\end{equation}
Thus,
\begin{equation}
  \label{eq:1p-nonzero-mode-relation}
  \bmat{Q}^{1\per}(k_1,r_2,r_3;\xi)
  =
  \E^{-\I k_1 r_1}
  \Kopt \Big(\E^{\I k_1 r_1} Q^{A,1\per}(k_1,r_2,r_3;\xi)\Big).
\end{equation}
To simplify the expressions in the following, we introduce the
auxiliary variables
\begin{equation}
  \label{eq:1p-aux}
  \rho = \sqrt{r_2^2+r_3^2}, \qquad
  U = \frac{k_1^2}{4\xi^2}, \qquad
  V = \xi^2 \rho^2.
\end{equation}

For the harmonic kernel $H(\bvec{r}) = 1/\lvert\bvec{r}\rvert$
and the Ewald screening function $\gE$
\eqref{eq:ewald-screening}, the integral $Q^{\harmonic,1\per}$,
defined by \eqref{eq:QA-1p-def}, was computed by
\cite{Tornberg2016}. The result is
\begin{equation}
  \label{eq:1p-QH-integral}
  Q^{H,1\per}(k_1,r_2,r_3;\xi) =
  \frac{1}{(2\pi)^2} \int_{\mathbb{R}^2} \frac{4\pi}{k_1^2 + \kappa_2^2 + \kappa_3^2}
  \E^{-(k_1^2 + \kappa_2^2 + \kappa_3^2)/(2\xi)^2}
  \E^{\I \kappa_2 r_2} \E^{\I \kappa_3 r_3}
  \, \mathrm{d}\kappa_2 \, \mathrm{d}\kappa_3
  = K_0(U,V),
\end{equation}
where $U$ and $V$ are given in \eqref{eq:1p-aux}, and
$K_\nu(\cdot,\cdot)$ is the incomplete modified Bessel
function of the second kind and order $\nu$ \citep{Harris2008},
defined by
\begin{equation}
  \label{eq:Knu-def}
  K_\nu(a,b) = \int_1^\infty \frac{\E^{-at-b/t}}{t^{\nu+1}} \, \mathrm{d}t.
\end{equation}
Now applying \eqref{eq:1p-nonzero-mode-relation} for the rotlet,
with $A=\harmonic$ and $\Kopt^\rotlet$ given by
\eqref{eq:diffop-rel-2}, we find that $\bmat{Q}^{\rotlet,1\per}$ is
antisymmetric
and given by
\begin{equation}
  \label{eq:1p-rotlet-Q}
  \bmat{Q}^{\rotlet,1\per}(k_1,r_2,r_3;\xi)
  =
  \begin{bmatrix}
    0 & 2 \xi^2 r_3 K_1(U,V) & -2 \xi^2 r_2 K_1(U,V) \\
    -2 \xi^2 r_3 K_1(U,V) & 0 & -\I k_1 K_0(U,V) \\
    2 \xi^2 r_2 K_1(U,V) & \I k_1 K_0(U,V) & 0 \\
  \end{bmatrix}
  .
\end{equation}
Here, we have used the relation
\begin{equation}
  \frac{\partial}{\partial b} K_\nu(a,b) = -K_{\nu+1}(a,b).
\end{equation}
(Note that if \eqref{eq:1p-evaluated-integrals-k0} is to be
evaluated at $(y,z) = (y_n,z_n)$, we must evaluate
\eqref{eq:1p-rotlet-Q} at $r_2=r_3=0$. This is not a problem, since
$K_\nu(a,0) = E_{\nu+1}(a)$, where $E_{\nu+1}(\cdot)$ denotes the
exponential integral of order $\nu+1$; note that $a=U>0$ here.)

For the biharmonic kernel $B(\bvec{r}) = \lvert \bvec{r} \rvert$
and the Hasimoto screening function $\gH$
\eqref{eq:hasimoto-screening}, the integral we want to compute is
\begin{align}
  \notag
  Q^{\biharmonic,1\per}(k_1,r_2,r_3;\xi) &=
  -\frac{1}{(2\pi)^2} \int_{\mathbb{R}^2}
  \frac{8\pi}{(k_1^2 + \kappa_2^2 + \kappa_3^2)^2}
  \E^{-(k_1^2 + \kappa_2^2 + \kappa_3^2) / (2\xi)^2}
  \left( 1 + \frac{k_1^2 + \kappa_2^2 + \kappa_3^2}{(2\xi)^2} \right)
  \E^{\I \kappa_2 r_2} \E^{\I \kappa_3 r_3}
  \, \mathrm{d}\kappa_2 \, \mathrm{d}\kappa_3
  \\
  \label{eq:1p-QB-integral-def}
  &=
  -\frac{2}{\pi} \int_{\mathbb{R}^2}
  \left( \frac{1}{(\alpha^2 + \kappa^2)^2}
  + \frac{1}{4\xi^2} \frac{1}{\alpha^2 + \kappa^2} \right)
  \E^{-(\alpha^2 + \kappa^2) / (2\xi)^2}
  \E^{\I \kappa_2 r_2} \E^{\I \kappa_3 r_3}
  \, \mathrm{d}\kappa_2 \, \mathrm{d}\kappa_3,
\end{align}
where we have introduced $\alpha := k_1$ and $\kappa^2 :=
\kappa_2^2 + \kappa_3^2$ to make the connection to the $D=2$ case
more apparent. Applying the relation \eqref{eq:2p-B-H-relation}
to the integrand of \eqref{eq:1p-QB-integral-def} and comparing with that of
\eqref{eq:1p-QH-integral}, we find that
\begin{equation}
  Q^{\biharmonic,1\per}(\alpha,r_2,r_3;\xi)
  =
  \frac{1}{\alpha}
  \frac{\partial}{\partial\alpha}
  Q^{H,1\per}(\alpha,r_2,r_3;\xi).
\end{equation}
Using the result of \eqref{eq:1p-QH-integral} and the relation
\begin{equation}
  \frac{\partial}{\partial a} K_\nu(a,b) = -K_{\nu-1}(a,b),
\end{equation}
one arrives that
\begin{equation}
  \label{eq:1p-QB-nonzero}
  Q^{\biharmonic,1\per}(k_1,r_2,r_3;\xi)
  =
  - \frac{1}{2 \xi^2} K_{-1}(U,V) .
\end{equation}
Note that the function $K_{-1}$ can be related to $K_1$ which
appeared in \eqref{eq:1p-rotlet-Q}, using the relation
\citep{Harris2008}
\begin{equation}
  a K_{-1}(a,b) - b K_1(a,b) = \E^{-a-b}.
\end{equation}

Applying \eqref{eq:1p-nonzero-mode-relation} for the stokeslet,
with $A=\biharmonic$ and $\Kopt^\stokeslet$ given by
\eqref{eq:diffop-rel-1}, one finds that
$\bmat{Q}^{\stokeslet,1\per}$ is symmetric
and given by
\begin{equation}
  \label{eq:1p-stokeslet-Q}
  \bmat{Q}^{\stokeslet,1\per}(k_1,r_2,r_3;\xi)
  =
  \begin{bmatrix}
    2K_0 - 2V K_1 & -\I k_1 r_2 K_0 & -\I k_1 r_3 K_0 \\
    -\I k_1 r_2 K_0 & 2 U K_{-1} + K_0 - 2 \xi^2 r_3^2 K_1 & 2 \xi^2 r_2 r_3 K_1 \\
    -\I k_1 r_3 K_0 & 2 \xi^2 r_2 r_3 K_1 & 2 U K_{-1} + K_0 - 2 \xi^2 r_2^2 K_1 \\
  \end{bmatrix}
  ,
\end{equation}
where $K_\nu$ ($\nu=-1,0,1$) is used as a shorthand for
$K_\nu(U,V)$, and $\rho$, $U$ and $V$ are as in
\eqref{eq:1p-aux}. For the stresslet, with $A=\biharmonic$ and
$\Kopt^\stresslet$ as in \eqref{eq:diffop-rel-3}, one finds that
$\bmat{Q}^{\stresslet,1\per}$ is symmetric with entries given by
\begin{equation}
  \label{eq:1p-stresslet-Q}
  \begin{array}{r@{\:}c@{\:}l}
    Q^{\stresslet,1\per}_{111}(k_1,r_2,r_3;\xi)
      &=&
      \displaystyle
      2 \I k_1 ( U K_{-1} + 3 K_0 - 3 V K_1 ),
    \\[5pt]
    Q^{\stresslet,1\per}_{112}(k_1,r_2,r_3;\xi)
      &=&
      \displaystyle
      4 \xi^2 r_2 ( U K_0 - 2 K_1 + V K_2 ),
    \\[5pt]
    Q^{\stresslet,1\per}_{113}(k_1,r_2,r_3;\xi)
      &=&
      \displaystyle
      4 \xi^2 r_3 ( U K_0 - 2 K_1 + V K_2 ),
    \\[5pt]
    Q^{\stresslet,1\per}_{122}(k_1,r_2,r_3;\xi)
      &=&
      \displaystyle
      2 \I k_1 ( U K_{-1} + \xi^2 (r_2^2 - r_3^2) K_1 ),
    \\[5pt]
    Q^{\stresslet,1\per}_{123}(k_1,r_2,r_3;\xi)
      &=&
      \displaystyle
      4\xi^2 \I k_1 r_2 r_3 K_1,
    \\[5pt]
    Q^{\stresslet,1\per}_{133}(k_1,r_2,r_3;\xi)
      &=&
      \displaystyle
      2 \I k_1 ( U K_{-1} + \xi^2 (r_3^2 - r_2^2) K_1 ),
    \\[5pt]
    Q^{\stresslet,1\per}_{222}(k_1,r_2,r_3;\xi)
      &=&
      \displaystyle
      -4\xi^2 r_2 ( 3 U K_0 + 3 K_1 - \xi^2 (r_2^2 + 3 r_3^2) K_2 ),
    \\[5pt]
    Q^{\stresslet,1\per}_{223}(k_1,r_2,r_3;\xi)
      &=&
      \displaystyle
      -4\xi^2 r_3 ( U K_0 + K_1 + \xi^2 (r_2^2 - r_3^2) K_2 ),
    \\[5pt]
    Q^{\stresslet,1\per}_{233}(k_1,r_2,r_3;\xi)
      &=&
      \displaystyle
      -4\xi^2 r_2 ( U K_0 + K_1 + \xi^2 (r_3^2 - r_2^2) K_2 ),
    \\[5pt]
    Q^{\stresslet,1\per}_{333}(k_1,r_2,r_3;\xi)
      &=&
      \displaystyle
      -4\xi^2 r_3 ( 3 U K_0 + 3 K_1 - \xi^2 (3r_2^2 + r_3^2) K_2 ),
  \end{array}
\end{equation}
again with auxiliary variables as in \eqref{eq:1p-aux} and $K_\nu
= K_\nu(U,V)$, $\nu=-1,0,1,2$.

\subsection{Singular case}
\label{app:1p-integrals-k0}

We wish to compute the integral
\begin{equation}
  \label{eq:1p-k0-integral-def}
  \bmat{Q}^{1\per,(0)}(r_2,r_3;\xi)
  :=
  \frac{1}{(2\pi)^2} \int_{\mathbb{R}^2}
  \wh{\kernel}{}^\fp(0,\kappa_2,\kappa_3)
  \E^{\I \kappa_2 r_2}
  \E^{\I \kappa_3 r_3}
  \, \mathrm{d}\kappa_2
  \, \mathrm{d}\kappa_3,
\end{equation}
that appears in \eqref{eq:1p-evaluated-integrals-k0}. The
integrand is now singular at the point $(\kappa_2,
\kappa_3)=(0,0)$. For the rotlet and stresslet kernels, the
singularity is of type $1/\kappa$, and can be removed by going to
polar coordinates. For the stokeslet, the singularity is of type
$1/\kappa^2$, and \eqref{eq:1p-k0-integral-def} must be
interpreted in the distributional sense. In the same way as in
\ref{app:2p-integrals-k0}, we have the relation
\begin{equation}
  \label{eq:1p-zero-mode-relation}
  \bmat{Q}^{1\per,(0)}(r_2,r_3;\xi) =
  \Kopt Q^{A,1\per,(0)}(r_2,r_3;\xi),
\end{equation}
with
\begin{equation}
  Q^{A,1\per,(0)}(r_2,r_3;\xi) :=
  \frac{1}{(2\pi)^2} \int_{\mathbb{R}^2}
  \wh{A}(0,\kappa_2,\kappa_3) \wh{\gamma}(0,\kappa_2,\kappa_3;\xi)
  \E^{\I \kappa_2 r_2} \E^{\I \kappa_3 r_3}
  \, \mathrm{d}\kappa_2 \, \mathrm{d}\kappa_3
  ,
\end{equation}
where $A$ is the harmonic kernel (for the rotlet) or biharmonic
kernel (for the stokeslet and stresslet). Also, we define the
two-dimensional inverse Fourier transform in the free directions
\begin{equation}
  \label{eq:1p-2d-generic-kernel}
  \pft{A}(0,r_2,r_3) := \fourier^{-1}\{ \wh{A}(0,\cdot,\cdot)\}(r_2,r_3),
\end{equation}
and use the relation
\begin{equation}
  \label{eq:1p-zero-mode-split}
  Q^{A,1\per,(0)}(r_2,r_3;\xi) =
  \pft{A}^\fp(0,r_2,r_3;\xi) = \pft{A}(0,r_2,r_3) - \lim_{k_1\to0}
  \Big(
  \pft{A}(k_1,r_2,r_3) - \pft{A}^\fp(k_1,r_2,r_3;\xi)
  \Big)
  .
\end{equation}
To simplify the expressions below, we introduce the auxiliary
variables
\begin{equation}
  \label{eq:1p-aux-zero}
  \widetilde{r}_2 = \frac{r_2}{\rho}, \qquad
  \widetilde{r}_3 = \frac{r_3}{\rho}, \qquad
  \varphi_0 = E_1(V) + \log(\rho^2) + 1, \qquad
  \varphi_1 = 1 - \E^{-V}, \qquad
  \varphi_2 = \varphi_1 - V\E^{-V},
\end{equation}
where $\rho$ and $V$ are as in \eqref{eq:1p-aux}, and
$E_\nu(\cdot)$ denotes the exponential integral of order $\nu$,
i.e.\ $E_\nu(a) = K_{\nu-1}(a,0)$.

For the harmonic ($A=H$), the kernel
$\pft{\harmonic}(0,r_2,r_3)$, cf.\ \eqref{eq:1p-2d-generic-kernel},
is the two-dimensional harmonic kernel. As in
\eqref{eq:2d-harmonic-kernel}, we introduce this kernel with an
arbitrary positive gauge constant $\ell_\harmonic$, i.e.\
\begin{equation}
  \pft{\harmonic}(0,r_2,r_3) = -2 \log(\rho/\ell_\harmonic).
\end{equation}
The derivation of $Q^{\harmonic,1\per,(0)}$ was done by
\cite{Tornberg2016} for the case $\ell_\harmonic=1$. For an
arbitrary $\ell_\harmonic$, the result becomes
\begin{equation}
  Q^{\harmonic,1\per,(0)}(r_2,r_3;\xi) =
  - \log(\rho^2/\ell_\harmonic^2) - E_1(\xi^2 \rho^2).
\end{equation}
Applying \eqref{eq:1p-zero-mode-relation} for the rotlet, with
$A=\harmonic$ and $\Kopt^\rotlet$ as in \eqref{eq:diffop-rel-2},
we find that $\bmat{Q}^{\rotlet,1\per,(0)}$ is antisymmetric and
given by
\begin{equation}
  \label{eq:1p-rotlet-Q0}
  \bmat{Q}^{\rotlet,1\per,(0)}(r_2,r_3;\xi)
  =
  \frac{2 \varphi_1}{\rho}
  \begin{bmatrix}
    0 & \widetilde{r}_3 & -\widetilde{r}_2 \\
    -\widetilde{r}_3 & 0 & 0 \\
    \widetilde{r}_2 & 0 & 0 \\
  \end{bmatrix}
  ,
\end{equation}
with auxiliary variables as in \eqref{eq:1p-aux-zero} and
\eqref{eq:1p-aux}. The limit as $(r_2,r_3) \to (0,0)$ is
$\bmat{Q}^{\rotlet,1\per,(0)}(0,0;\xi) = \bmatu{0}$.
Notably, $\bmat{Q}^{\rotlet,1\per,(0)}$ for the rotlet does not
depend on the gauge constant $\ell_\harmonic$.

For the biharmonic ($A=\biharmonic$), we will use
\eqref{eq:1p-zero-mode-split}, and we start by computing
$\pft{\biharmonic}(0,r_2,r_3)$ and $\pft{\biharmonic}(k_1,r_2,r_3)$.
These are defined by taking the two-dimensional inverse Fourier
transform in the free ($r_2$ and $r_3$) directions of
$\wh{\biharmonic}(k_1,\kappa_2,\kappa_3) =
-8\pi/(k_1^2+\kappa_2^2+\kappa_3^2)^2$. For $k_1>0$, the inverse
Fourier transform exists in the classical sense; it can be
computed by noting that $\wh{\biharmonic}$ is radial in
$(\kappa_2,\kappa_3)$, and using the formula
\begin{equation}
  \label{eq:2d-radial-IFT}
  f(r_2,r_3) = f(\rho) = \frac{1}{2\pi} \int_0^\infty
  J_0(\kappa \rho) \wh{f}(\kappa) \kappa \, \mathrm{d}\kappa
\end{equation}
for the inverse Fourier transform of a two-dimensional radial
function $\wh{f}(\kappa_2,\kappa_3) = \wh{f}(\kappa)$. Here,
$\rho = \sqrt{r_2^2+r_3^2}$, $\kappa=\sqrt{\kappa_2^2+\kappa_3^2}$,
and $J_0(\cdot)$ is the Bessel function of the first kind and order~0.
Applying \eqref{eq:2d-radial-IFT} to $\wh{B}$ yields
\begin{equation}
  \label{eq:1p-Bbar-tot-nonzero}
  \pft{\biharmonic}(k_1,r_2,r_3) = -\frac{2 \rho}{\lvert k_1 \rvert}
  K_1(\lvert k_1 \rvert \rho), \qquad k_1 \neq 0,
\end{equation}
where $K_1(\cdot)$ is the modified Bessel function of the second
kind and order~1. Note that $K_\nu(\cdot)$ of one argument is not
the same function as $K_\nu(\cdot,\cdot)$ of two arguments
introduced in \eqref{eq:Knu-def}.

For $k_1=0$, the inverse Fourier transform must be interpreted in
the distributional sense. By noting that the Fourier transform of
\begin{equation}
  - (\nabla_2^2 + \nabla_3^2)^2 \pft{\biharmonic}(0,r_2,r_3) =
  8\pi\delta(r_2)\delta(r_3)
\end{equation}
is precisely $\wh{\biharmonic}(0,\kappa_2,\kappa_3) =
-8\pi/(\kappa_2^2 + \kappa_3^2)^2$, we see that $\pft{\biharmonic}(0,r_2,r_3)$
is the fundamental solution of the two-dimensional biharmonic
equation, already introduced in \eqref{eq:2d-biharmonic-kernel},
\begin{equation}
  \label{eq:1p-Bbar-tot-zero}
  \pft{\biharmonic}(0,r_2,r_3) = - \rho^2 \log(\rho/\ell_\biharmonic) + \cB,
\end{equation}
where $\ell_\biharmonic$ is an arbitrary positive gauge constant,
and we will here set $\cB=0$.

We now compute the limit found in \eqref{eq:1p-zero-mode-split},
with $\pft{\biharmonic}(k_1,r_2,r_3)$ taken from \eqref{eq:1p-Bbar-tot-nonzero}
and $\pft{\biharmonic}^\fp(k_1,r_2,r_3;\xi)$ taken from
\eqref{eq:1p-QB-nonzero}. The limit is
\begin{equation}
  \label{eq:1p-Bbar-R-limit-1}
  \pft{\biharmonic}^\rp(0,r_2,r_3;\xi) =
  \lim_{k_1\to0} \left(
    -\frac{2 \rho}{\lvert k_1 \rvert} K_1(\lvert k_1 \rvert \rho)
    +
    \frac{1}{2 \xi^2} K_{-1}(U,V)
  \right).
\end{equation}
Recall that $k_1$ is found also in $U$, see \eqref{eq:1p-aux}.
Using the relation \citep{Harris2008}
\begin{equation}
  K_\nu(a,b) + K_{-\nu}(b,a) = 2 (a/b)^{\nu/2} K_\nu(2\sqrt{ab}),
\end{equation}
with $\nu=1$, we can rewrite \eqref{eq:1p-Bbar-R-limit-1} as (note that $2
\sqrt{UV} = \lvert k_1 \rvert \rho$, and $\sqrt{V/U} = 2 \xi^2
\rho / \lvert k_1 \rvert$)
\begin{equation}
  \label{eq:1p-Bbar-R-limit-2}
  \pft{\biharmonic}^\rp(0,r_2,r_3;\xi) =
  - \frac{1}{2 \xi^2}
  \lim_{k_1\to0} K_1(V,U)
  =
  - \frac{1}{2 \xi^2} K_1(\xi^2 \rho^2, 0)
  =
  - \frac{1}{2 \xi^2} E_2(\xi^2 \rho^2)
  .
\end{equation}
Inserting this, and $\pft{\biharmonic}(0,r_2,r_3)$ from \eqref{eq:1p-Bbar-tot-zero},
into \eqref{eq:1p-zero-mode-split}, we get
\begin{equation}
  \label{eq:1p-QB-F-result}
  Q^{\biharmonic,1\per,(0)}(r_2,r_3;\xi) =
  \pft{\biharmonic}^\fp(0,r_2,r_3;\xi)
  =
  - \rho^2 \log(\rho/\ell_\biharmonic)
  + \frac{1}{2 \xi^2} E_2(\xi^2 \rho^2)
  .
\end{equation}
At this point, we note that $-\rho^2 \log(\rho/\ell_\biharmonic) =
-\rho^2 \log(\rho) + \rho^2 \log(\ell_\biharmonic)$ and that the contribution to 
$\bmat{Q}^{\stokeslet,1\per,(0)}$ from the term $\rho^2
\log(\ell_\biharmonic)$, when applying $\Kopt^\stokeslet$ for the stokeslet,
given by \eqref{eq:diffop-rel-1}, is
\begin{equation}
  \Kopt^\stokeslet (\rho^2 \log(\ell_\biharmonic)) = \log(\ell_\biharmonic)
  \begin{bmatrix}
    4 & 0 & 0 \\
    0 & 2 & 0 \\
    0 & 0 & 2 \\
  \end{bmatrix},
\end{equation}
which upon insertion into \eqref{eq:1p-evaluated-integrals-k0}
leads to the constant contribution
\begin{equation}
  \label{eq:constant-stokeslet-field-contribution-1p}
  \bvec{u}^{\stokeslet,1\per,\fp,\bvec{k}^\per=\bvec{0},\mathrm{extra}}
  =
  \log(\ell_\biharmonic)
  \begin{bmatrix}
    4 & 0 & 0 \\
    0 & 2 & 0 \\
    0 & 0 & 2 \\
  \end{bmatrix}
  \frac{1}{L_1} \sum_{n=1}^N
  \bvec{f}(\bvec{x}_n)
\end{equation}
to the flow field $\bvec{u}^{\stokeslet,1\per,\fp,\bvec{k}^\per=\bvec{0}}$.
(This is the same contribution as
\eqref{eq:1p-stokeslet-extra-contribution}, with equality if
$\ell_\biharmonic/\ell_\harmonic=\E$.)
The constant contribution \eqref{eq:constant-stokeslet-field-contribution-1p}
can always be added to the solution in a separate step, if
desired, so we will not include it in the following results.
For the stresslet, the contribution from $\ell_\biharmonic$ is
always zero, i.e.\ $\Kopt^\stresslet (\rho^2
\log(\ell_\biharmonic)) = \bmatu{0}$.

Thus, setting $\ell_\biharmonic=1$ in \eqref{eq:1p-QB-F-result} and applying
\eqref{eq:1p-zero-mode-relation} for the stokeslet, we find that
$\bmat{Q}^{\stokeslet,1\per,(0)}$ is symmetric and given by
\begin{equation}
  \label{eq:1p-stokeslet-Q0}
  \bmat{Q}^{\stokeslet,1\per,(0)}(r_2,r_3;\xi)
  =
  -
  \begin{bmatrix}
    2(\varphi_0+\varphi_1) & 0 & 0 \\
    0 & \varphi_0 + 2\widetilde{r}_3^2 \varphi_1 & -2 \widetilde{r}_2 \hspace{0.5pt} \widetilde{r}_3 \varphi_1 \\
    0 & -2 \widetilde{r}_2 \hspace{0.5pt} \widetilde{r}_3 \varphi_1 & \varphi_0 + 2 \widetilde{r}_2^2 \varphi_1 \\
  \end{bmatrix}
  .
\end{equation}
Using the asymptotic expansion $E_1(t) = -\gamma - \log(t) +
O(t)$, we find that
\begin{equation}
  \lim_{\rho \to 0} \varphi_0 = -\gamma - \log(\xi^2) + 1,
\end{equation}
where $\gamma = 0.577\,215\,664\,9\ldots$ is the Euler--Mascheroni
constant. Using this and the fact that $\lim_{\rho \to 0}
\varphi_1 = 0$, we find that the limit of
\eqref{eq:1p-stokeslet-Q0} as $\rho \to 0$ is
\begin{equation}
  \label{eq:1p-stokeslet-Q0-limit}
  \bmat{Q}^{\stokeslet,1\per,(0)}(0,0;\xi)
  =
  \left( \gamma + \log(\xi^2) - 1 \right)
  \begin{bmatrix}
    2 & 0 & 0 \\
    0 & 1 & 0 \\
    0 & 0 & 1 \\
  \end{bmatrix}
  .
\end{equation}
For the stresslet, $\bmat{Q}^{\stresslet,1\per,(0)}$ is symmetric
with entries given by
\begin{equation}
  \label{eq:1p-stresslet-Q0}
  \begin{array}{r@{\:}c@{\:}l@{\hspace{4em}}r@{\:}c@{\:}l}
    Q^{\stresslet,1\per,(0)}_{111}(r_2,r_3;\xi)
      &=&
      0,
    &
    Q^{\stresslet,1\per,(0)}_{133}(r_2,r_3;\xi)
      &=&
      0,
    \\[5pt]
    Q^{\stresslet,1\per,(0)}_{112}(r_2,r_3;\xi)
      &=&
      R \widetilde{r}_2 (2 \varphi_1 - \varphi_2),
    &
    Q^{\stresslet,1\per,(0)}_{222}(r_2,r_3;\xi)
      &=&
      R \widetilde{r}_2 (3 \varphi_1 - (\widetilde{r}_2^2 + 3 \widetilde{r}_3^2) \varphi_2),
    \\[5pt]
    Q^{\stresslet,1\per,(0)}_{113}(r_2,r_3;\xi)
      &=&
      R \widetilde{r}_3 (2 \varphi_1 - \varphi_2),
    &
    Q^{\stresslet,1\per,(0)}_{223}(r_2,r_3;\xi)
      &=&
      R \widetilde{r}_3 (\varphi_1 - (\widetilde{r}_3^2 - \widetilde{r}_2^2) \varphi_2),
    \\[5pt]
    Q^{\stresslet,1\per,(0)}_{122}(r_2,r_3;\xi)
      &=&
      0,
    &
    Q^{\stresslet,1\per,(0)}_{233}(r_2,r_3;\xi)
      &=&
      R \widetilde{r}_2 (\varphi_1 - (\widetilde{r}_2^2 - \widetilde{r}_3^2) \varphi_2),
    \\[5pt]
    Q^{\stresslet,1\per,(0)}_{123}(r_2,r_3;\xi)
      &=&
      0,
    &
    Q^{\stresslet,1\per,(0)}_{333}(r_2,r_3;\xi)
      &=&
      R \widetilde{r}_3 (3 \varphi_1 - (3\widetilde{r}_2^2 + \widetilde{r}_3^2) \varphi_2),
  \end{array}
\end{equation}
with $R=-4/\rho$, and other auxiliary variables as in
\eqref{eq:1p-aux-zero}. The limit of \eqref{eq:1p-stresslet-Q0}
as $\rho \to 0$ is
\begin{equation}
  \label{eq:1p-stresslet-Q0-limit}
  \bmat{Q}^{\stresslet,1\per,(0)}(0,0;\xi) = \bmatu{0}.
\end{equation}

\section{Fourier-space truncation error estimates for stokeslet
and stresslet}
\label{app:truncation-estimates}

We here give the derivation of the truncation error
estimates \eqref{eq:fs-trunc-est-stokeslet} and \eqref{eq:fs-trunc-est-stresslet}
for the stokeslet and stresslet, respectively. The derivation
is based on the technique used by \cite{afKlinteberg2017},
adapted to the triply periodic case. The resulting estimates can,
however, be used in any periodicity.

\subsection{Stokeslet estimate}

For simplicity, the primary cell $\boxvar$ is assumed to be a cube
of side length $L$. We consider the Fourier-space Ewald sum
\eqref{eq:fourier-space-potential-3p} for the stokeslet in the
triply periodic case, i.e.\
\begin{equation}
  \label{eq:stokeslet-FS-3P-potential}
  \bvec{u}^{\stokeslet,3\per,\fp,\bvec{k}\neq\bvec{0}}(\bvec{x};\xi)
  = \frac{1}{L^3} \sum_{n=1}^N
  \sum_{\substack{\bvec{k} \in \wavenumset^3 \\ \bvec{k} \neq \bvec{0}}}
  \wh{\bmat{\stokeslet}}{}^\fp(\bvec{k};\xi)
  \bvec{f}(\bvec{x}_n)
  \E^{\I \bvec{k} \cdot (\bvec{x}-\bvec{x}_n)}
  ,
\end{equation}
where $\wavenumset^3$ is given by \eqref{eq:wavenumbers-3p}, and
$\wh{\bmat{\stokeslet}}{}^\fp$ is given by \eqref{eq:stokeslet-fourier}.
We will assume that \eqref{eq:stokeslet-FS-3P-potential} is
truncated outside a sphere of radius $\lvert \bvec{k} \rvert = k_\infty$
for some $k_\infty > 0$. (In the SE method, the sum is actually
truncated outside a cube of side length $2 k_\infty$, but that
can only lead to a smaller error, and the difference will in any
case not be large.)
The pointwise error from truncating
\eqref{eq:stokeslet-FS-3P-potential} is given by
\begin{equation}
  \label{eq:stokeslet-truncation-error}
  \bvec{e}^{\stokeslet,\fp}(\bvec{x};\xi)
  := \frac{1}{L^3} \sum_{n=1}^N
  \sum_{\substack{\bvec{k} \in \wavenumset^3 \\ \lvert\bvec{k}\rvert > k_\infty}}
  \wh{\bmat{\stokeslet}}{}^\fp(\bvec{k};\xi)
  \bvec{f}(\bvec{x}_n)
  \E^{\I \bvec{k} \cdot (\bvec{x}-\bvec{x}_n)}
  .
\end{equation}
To be able to follow \cite{afKlinteberg2017}, the sum over
$\bvec{k}$ is approximated by an integral. Multiplying and
dividing \eqref{eq:stokeslet-truncation-error} by $(\Delta k)^3$,
where $\Delta k = 2\pi/L$ is the wavenumber resolution, allows us
to make the approximation
\begin{equation}
  \label{eq:stokeslet-truncation-integral}
  \bvec{e}^{\stokeslet,\fp}(\bvec{x};\xi)
  \approx \frac{1}{(2\pi)^3} \sum_{n=1}^N
  \int_{\lvert \bvec{k} \rvert > k_\infty}
  \wh{\bmat{\stokeslet}}{}^\fp(\bvec{k};\xi)
  \bvec{f}(\bvec{x}_n)
  \E^{\I \bvec{k} \cdot (\bvec{x}-\bvec{x}_n)}
  \mathrm{d}\bvec{k}
  .
\end{equation}
Let us define
\begin{align}
  \notag
  E_{jl} (\bvec{r}) &:= \frac{1}{(2\pi)^3}
  \int_{\lvert\bvec{k}\rvert > k_\infty}
  \wh{S}{}^\fp_{jl}(\bvec{k};\xi)
  \E^{\I \bvec{k} \cdot \bvec{r}}
  \, \mathrm{d} \bvec{k}
  \\ &=
  \label{eq:stokeslet-Ejl-step1}
  \frac{1}{\pi^2}
  \int_{\lvert\bvec{k}\rvert > k_\infty}
  \frac{1}{\lvert \bvec{k} \rvert^2}
  \left(
    \delta_{jl} -
    \frac{k_j k_l}{\lvert \bvec{k} \rvert^2}
   \right)
  \E^{-\lvert \bvec{k} \rvert^2/(2\xi)^2}
  \left(
  1 + \frac{\lvert \bvec{k} \rvert^2}{(2\xi)^2}
  \right)
  \E^{\I \bvec{k} \cdot \bvec{r}}
  \, \mathrm{d} \bvec{k}
  ,
\end{align}
such that $\bvec{e}^{\stokeslet,\fp}(\bvec{x};\xi) \approx
\sum_{n=1}^N \bmat{E}(\bvec{x}-\bvec{x}_n) \bvec{f}(\bvec{x}_n)$.

As done by \cite{afKlinteberg2017}, we approximate the
directional component by its root mean square value (computed
using spherical coordinates),
\begin{equation}
  \delta_{jl} - \frac{k_j k_l}{\lvert \bvec{k} \rvert^2}
  \approx
  \sqrt{\frac{1}{9} \sum_{j,l=1}^3 \left(
  \delta_{jl} - \frac{k_j k_l}{\lvert \bvec{k} \rvert^2}
  \right)^2}
  =
  \frac{\sqrt{2}}{3},
\end{equation}
and keep only the highest-order term with respect to $\lvert
\bvec{k} \rvert$ in \eqref{eq:stokeslet-Ejl-step1}, since that
term will dominate the error for large $k_\infty$. Introducing
spherical coordinates $[k, \theta, \varphi]$, with the coordinate
system chosen such that $\bvec{k} \cdot \bvec{r} = kr \cos \theta$
(where $k = \lvert \bvec{k} \rvert$ and $r = \lvert \bvec{r} \rvert$),
we get the result
\begin{equation}
  E_{jl}(\bvec{r}) \approx
  2\pi
  \frac{\sqrt{2}}{12\pi^2\xi^2}
  \int_{k_\infty}^\infty
  \int_{0}^{\pi}
  \E^{-k^2/(2\xi)^2}
  \E^{\I kr \cos \theta}
  k^2 \sin \theta
  \, \mathrm{d} \theta \, \mathrm{d} k
\end{equation}
where the integral over $\varphi$ could be evaluated to $2\pi$ directly.
Computing also the integral over $\theta$ yields
\begin{equation}
  E_{jl}(\bvec{r}) \approx
  \frac{\sqrt{2}}{3\pi\xi^2}
  \int_{k_\infty}^\infty
  k^2 \frac{\sin(kr)}{kr}
  \E^{-k^2/(2\xi)^2}
  \, \mathrm{d} k.
\end{equation}
The exponential decay of $\E^{-k^2/(4\xi^2)}$ makes sure that the
dominant contribution comes from the beginning of the interval,
where $k \approx k_\infty$. This allows the approximation
\begin{equation}
  \label{eq:stokeslet-Ejl-step2}
  E_{jl}(\bvec{r}) \approx
  \frac{\sqrt{2}}{3\pi\xi^2}
  \frac{k_\infty}{r}
  \int_{k_\infty}^\infty
  \sin(kr) \E^{-k^2/(2\xi)^2}
  \, \mathrm{d} k.
\end{equation}
The remaining integral can be computed exactly in terms of the
error function, but to get a more manageable error estimate, we
will approximate it. Writing $\sin(kr) = \mathrm{Im} \{ \E^{\I kr} \}$,
and once again using $k \approx k_\infty$, we get
\begin{align}
  \notag
  \int_{k_\infty}^\infty
  \sin(kr) \E^{-k^2/(2\xi)^2}
  \, \mathrm{d} k
  &=
  \mathrm{Im} \left\{
  \int_{k_\infty}^\infty
  \E^{\I kr -k^2/(2\xi)^2}
  \, \mathrm{d}k
  \right\}
  =
  \mathrm{Im} \left\{
  \int_{k_\infty}^\infty
  \frac{\I r - k/(2\xi^2)}{\I r - k/(2\xi^2)}
  \E^{\I kr -k^2/(2\xi)^2}
  \, \mathrm{d}k
  \right\}
  \\ &\approx
  \notag
  \mathrm{Im}
  \left\{
    \frac{1}{\I r - k_\infty/(2\xi^2)}
    \int_{k_\infty}^\infty
    \Big(\I r - k/(2\xi^2)\Big)
    \E^{\I kr -k^2/(2\xi)^2}
    \, \mathrm{d}k
  \right\}
  \\ &=
  \label{eq:stokeslet-Ejl-part-integral}
  \frac{\E^{-k_\infty^2/(2\xi)^2}}{r^2 + k_\infty^2/(4\xi^4)}
  \left(
    r \cos(k_\infty r) + \frac{k_\infty}{2\xi^2} \sin(k_\infty r)
  \right).
\end{align}
Again assuming that $k_\infty$ is large (also compared to $r$),
we keep only the highest-order terms with respect to $k_\infty$.
Inserting \eqref{eq:stokeslet-Ejl-part-integral} into
\eqref{eq:stokeslet-Ejl-step2} then yields
\begin{equation}
  \label{eq:stokeslet-Ejl-step3}
  E_{jl}(\bvec{r}) \approx
  \frac{2\sqrt{2}}{3\pi}
  \frac{\sin(k_\infty r)}{r}
  \E^{-k_\infty^2/(2\xi)^2}
  =: E(\bvec{r})
  ,
\end{equation}
which is the final simplification of \eqref{eq:stokeslet-Ejl-step1}.
Note that this result no longer depends on the indices $j$ and $l$.

Let us denote the root mean square error of
$\bvec{e}^{\stokeslet,\fp}(\bvec{x};\xi)$
by
\begin{equation}
  \delta \bvec{e}^{\stokeslet,\fp} := \sqrt{
    \frac{1}{N} \sum_{m=1}^N \sum_{j=1}^3
    [e_j^{\stokeslet,\fp}(\bvec{x}_m;\xi)]^2
  }
  .
\end{equation}
Then, as $\bvec{e}^{\stokeslet,\fp}(\bvec{x};\xi) \approx
\sum_{n=1}^N \bmat{E}(\bvec{x}-\bvec{x}_n) \bvec{f}(\bvec{x}_n)$,
we use \cite[Lemma~1]{afKlinteberg2017} to approximate
$(\delta\bvec{e}^{\stokeslet,\fp})^2$ by
\begin{equation}
  (\delta \bvec{e}^{\stokeslet,\fp})^2 \approx
  \sum_{n=1}^N \sum_{j=1}^3
  \sum_{l=1}^3 [f_l(\bvec{x}_n)]^2
  \frac{1}{\lvert V \rvert}
  \int_V [E_{jl}(\bvec{r})]^2 \, \mathrm{d} \bvec{r}
  \approx
  \frac{3Q}{\lvert V \rvert}
  \int_V [E(\bvec{r})]^2 \, \mathrm{d} \bvec{r}
  ,
\end{equation}
where $V$ is a volume containing all vectors $\bvec{r}_{mn} =
\bvec{x}_m - \bvec{x}_n$, and $Q$ is defined as in
\eqref{eq:vector-source-quantity}. Selecting $V$ to be a sphere
of radius $\rho := \sqrt{3}L/2$, we get (inserting \eqref{eq:stokeslet-Ejl-step3}
and computing the integral over $V$ using spherical coordinates)
\begin{equation}
  (\delta \bvec{e}^{\stokeslet,\fp})^2
  \approx
  \frac{8Q}{\pi^2\rho^3}
  \E^{-k_\infty^2/(2\xi^2)}
  \left(
  \frac{\rho}{2} - \frac{\sin(2 k_\infty \rho)}{4 k_\infty}
  \right)
  .
\end{equation}
Since $k_\infty$ is large, the term involving $\sin(2 k_\infty
\rho)/k_\infty$ can be neglected. This leads to the estimate
\begin{equation}
  \delta \bvec{e}^{\stokeslet,\fp}
  \approx
  \frac{2 \sqrt{Q}}{\pi\rho}
  \E^{-k_\infty^2/(2\xi)^2}
  =
  \frac{4 \sqrt{Q}}{\sqrt{3} \pi L}
  \E^{-k_\infty^2/(2\xi)^2}
  ,
\end{equation}
which is \eqref{eq:fs-trunc-est-stokeslet}.

\subsection{Stresslet estimate}

The derivation for the stresslet \eqref{eq:stresslet-fourier}
proceeds in exactly the same way as for the stokeslet, with the
only difference being that the stresslet has an extra factor of
$\lvert \bvec{k} \rvert$ compared to the stokeslet, and that the
root mean square of the directional component is
\begin{equation}
  \sqrt{\frac{1}{27} \sum_{j,l,m=1}^3
  \left(
  \frac{\delta_{jl} k_m + \delta_{mj} k_l + \delta_{lm} k_j}{\lvert \bvec{k} \rvert}
  - 2 \frac{k_j k_l k_m}{\lvert \bvec{k} \rvert^3}
  \right)
  }
  = \sqrt{\frac{7}{27}},
\end{equation}
instead of $\sqrt{2}/3$. This means that the stresslet estimate
will differ from the stokeslet estimate by a factor of $k_\infty
\sqrt{7/6}$. Multiplying \eqref{eq:fs-trunc-est-stokeslet} by
this factor yields \eqref{eq:fs-trunc-est-stresslet}.

\bibliography{bibliography}

\end{document}